\documentclass[3p,a4paper,11pt]{elsarticle}

\usepackage{lineno,hyperref}
\usepackage{latexsym}
\usepackage{color}
\usepackage[table]{xcolor}
\usepackage{mathrsfs}
\usepackage{amssymb,amstext,amsfonts,amsmath,amsthm}
\usepackage{bm}
\usepackage{multirow}
\usepackage[ruled,linesnumbered]{algorithm2e}
\usepackage{booktabs}

\usepackage{jlcode}
\usepackage{listings}
\usepackage[utf8]{inputenc} 
\usepackage[english]{babel}
\usepackage{hyphenat}

\usepackage{calc}
\usepackage{graphicx,import}

\usepackage{soul}

\usepackage{caption}
\usepackage{subcaption}

\usepackage{arydshln}
\usepackage{natbib}
\usepackage[misc,geometry]{ifsym} 
\usepackage{empheq}
\usepackage{longtable}
\usepackage{comment}

\biboptions{authoryear}

\usepackage{tikz}

\def\bbA{{\mathbb A}} 
 
\def\bbC{{\mathbb C}}

\def\bbR{{\mathbb R}}

\def\tr[#1]{\textcolor{red}{#1}}
\def\tg[#1]{\textcolor{green}{#1}}
\def\tb[#1]{\textcolor{blue}{#1}}

\def\mcC{{\mathcal C}}

\def\mcF{{\mathcal F}}

\def\mcR{{\mathcal R}}

\def\mcT{{\mathcal T}}
\def\mcU{{\mathcal U}}
\def\mcV{{\mathcal V}}
\def\mcX{{\mathcal X}}

\def\dif{{\text{d}}}

\def\r{ \right) }

\newcommand{\bfA}{{\bf A}}

\newcommand{\bfB}{{\bf B}}
\newcommand{\bfc}{{\bf c}}
\newcommand{\bfC}{{\bf C}}
\newcommand{\bfd}{{\bf d}}
\newcommand{\bfD}{{\bf D}}
\newcommand{\bfe}{{\bf e}}

\newcommand{\bff}{{\bf f}}
\newcommand{\bfF}{{\bf F}}
\newcommand{\bfg}{{\bf g}}
\newcommand{\bfG}{{\bf G}}

\newcommand{\bfI}{{\bf I}}

\newcommand{\bfK}{{\bf K}}
\newcommand{\bfM}{{\bf M}}

\newcommand{\bfN}{{\bf N}}

\newcommand{\bfP}{{\bf P}}

\newcommand{\bfr}{{\bf r}}

\newcommand{\bfu}{{\bf u}}

\newcommand{\bfv}{{\bf v}}

\newcommand{\bfx}{{\bf x}}

\newcommand{\bfPhi}{\boldsymbol{\Phi}}

\newcommand{\bfxi}{\boldsymbol{\xi}}

\newcommand{\bszer}{\boldsymbol{0}}

\newcommand{\bstheta}{\boldsymbol{\theta}}

\def\R{{\mathbb R}}
\def\x{{\mathbf x}}

\DeclareMathOperator*{\Ass}{ \mathlarger{\mathlarger{\mathlarger{\boldsymbol{\bbA}}}} }

\renewcommand{\leq}{\leqslant}
\renewcommand{\geq}{\geqslant}

\newcommand{\Reach}[3]{\mathcal{R}}
\def\mcX{\mathcal{X}}
\def\mcR{\mathcal{R}}
\def\mcU{\mathcal{U}}
\def\mcF{\mathcal{F}}
\usepackage{lineno}
\newcommand*\patchAmsMathEnvironmentForLineno[1]{%
	\expandafter\let\csname old#1\expandafter\endcsname\csname #1\endcsname
	\expandafter\let\csname oldend#1\expandafter\endcsname\csname end#1\endcsname
	\renewenvironment{#1}%
	{\linenomath\csname old#1\endcsname}%
	{\csname oldend#1\endcsname\endlinenomath}}%
\newcommand*\patchBothAmsMathEnvironmentsForLineno[1]{%
	\patchAmsMathEnvironmentForLineno{#1}%
	\patchAmsMathEnvironmentForLineno{#1*}}%
\AtBeginDocument{%
	\patchBothAmsMathEnvironmentsForLineno{equation}%
	\patchBothAmsMathEnvironmentsForLineno{align}%
	\patchBothAmsMathEnvironmentsForLineno{flalign}%
	\patchBothAmsMathEnvironmentsForLineno{alignat}%
	\patchBothAmsMathEnvironmentsForLineno{gather}%
	\patchBothAmsMathEnvironmentsForLineno{multline}%
}

\usepackage[]{algorithm2e}

\usepackage{hyperref}

\usepackage{tikz}
\usetikzlibrary{shapes,arrows}

\usepackage{relsize}

\definecolor{darkcandyapplered}{rgb}{0.64, 0.0, 0.0}
\definecolor{darkred}{rgb}{0.55, 0.0, 0.0}

\def\tr[#1]{\textcolor{darkcandyapplered}{#1}}
\def\tb[#1]{\textcolor{blue}{#1}}

\journal{Computers \& Structures}

\usepackage{bbm}

\newcommand{\rev}{} 

\usepackage{ulem}
\usepackage{enumitem,amssymb}
\newlist{todolist}{itemize}{2}
\setlist[todolist]{label=$\square$}
\usepackage{pifont}

\newtheorem{prop}{Proposition}

\begin{document}

\begin{frontmatter}
	
  \title{Combining Set Propagation with Finite Element Methods for Time Integration in Transient Solid Mechanics Problems}

  \author[1]{Marcelo Forets}
  \author[2]{Daniel Freire Caporale}
  \author[3]{Jorge M. {P\'erez Zerpa}}

\address[1]{Departamento de Matemática y Aplicaciones, Centro Universitario Regional del Este, Universidad de la República, Maldonado, Uruguay}
\address[2]{Instituto de Física, Facultad de Ciencias, Universidad de la República, Montevideo, Uruguay}
\address[3]{Instituto de Estructuras y Transporte, Facultad de Ingeniería, Universidad de la República, Montevideo, Uruguay}

\begin{abstract}
	The Finite Element Method (FEM) is the gold standard for spatial discretization in numerical simulations for a wide spectrum of real-world engineering problems.
	Prototypical areas of interest include linear heat transfer and linear structural dynamics problems modeled with partial differential equations (PDEs).
	While different algorithms for direct integration of the equations of motion exist, exploring all feasible behaviors for varying loads, initial states and fluxes in models with large numbers of degrees of freedom remains a challenging task.
	In this article we propose a novel approach, based in set propagation methods and motivated by recent advances in the field of Reachability Analysis.
	Assuming a set of initial states and inputs, the proposed method consists in the construction of a union of sets (flowpipe) that enclose the infinite number of solutions of the spatially discretized PDE.
	We present the numerical results obtained in {\rev five} examples to illustrate the capabilities of our approach, and {\rev compare its performance against reference numerical integration methods}.
	We conclude that{\rev, for problems with single known initial conditions, the proposed method is accurate. For problems with uncertain initial conditions included in sets, the proposed method can compute all the solutions of the system more efficiently than numerical integration methods.
}
	%
\end{abstract}

\begin{keyword}
Reachability Analysis \sep
Finite Element Method \sep
Heat Transfer \sep
Structural Dynamics  \sep
Numerical Verification
  \end{keyword}	
\end{frontmatter}


\section{Introduction} \label{sec:introduction}

Transient problems are of great relevance in a large and diverse set of disciplines in Engineering.
{\rev In structural dynamics,} problems with time-varying loads are formulated in {\rev applications such as } wind turbine modeling \citep{Xue2020} or the study of nonlinear behavior of cables \citep{martinelli2001numerical}.
{\rev In} heat transfer analysis, {\rev transient problems arise in} applications ranging from designing massive concrete structures \citep{Wilson1974,Tahersima2017}, to simulating the effects {\rev of COVID-19} on the respiratory tract \citep{Grau-Bartual2020}.

{\rev
	The governing equations of the mathematical models arising in transient problems are numerically solved considering both a discretization in space and time.
	The Finite Element Method (FEM), emerged in the second half of the 20th century \citep{Zienkiewicz1972}, has become the gold-standard approach for spatial discretization of PDEs.
	Regarding time integration, methods keep being actively developed \citep{Bathe2012,Kim2018,Malakiyeh2019, Li2017a,kim2021}, highlighting the challenge this represents.
}

{\rev In addition, all mathematical models developed and applied to any real-world problem are subject to uncertainties.}
The material properties, the loads applied, the geometries, {\rev the initial conditions,} etc., must be provided as an input of the models, and are always submitted to errors.
Hence, the effect of {\rev uncertainties} over simulation results must be quantified.
In the last decades, several methods have been developed using different mathematical and numerical techniques, which can be classified into probabilistic and non-probabilistic approaches.
The main probabilistic approaches are Uncertainty Propagation and Monte Carlo methods. The former is based on the calculation of probability density functions for sampling statistics and has been successfully applied in structural systems under dynamic loads \citep{capillon2016,imholz2020,horn2021}. Monte Carlo methods rely on repeated numerical simulations starting from different initial conditions and parameter choices \citep{shinozuka1972monte}.
Non-probabilistic approaches include interval methods \citep{shu2001interval, muhanna2011interval}, and fuzzy logic methods {\rev \citep{muhanna1999formulation,Haiqing2010a}}.

The most common practice to solve transient problems relies on exploring different behaviors by simulation and testing.
However, if uncertainty is present and the application requires an exhaustive exploration of the state space, such approach becomes computationally intractable.
For instance, to guarantee that a certain maximum temperature in a massive concrete structure is not exceeded \citep{Tahersima2017} can be interpreted as a \textit{verification problem}.
A prominent framework used to solve verification problems consists in approximating the set of states that are reachable by a dynamical system, from all initial states and for all admissible inputs and parameters.
Such approach is called \textit{reachability analysis} and it has been recently applied to solve verification, control, path planning and stability problems in diverse domains \citep{althoff2020set}.
One of the main approaches to reachability analysis relies on \textit{set propagation}, where the solution of an ordinary differential equation (ODE) is expressed in terms of sets rather than numbers. 
Several algorithms have been devised for linear \citep{girard2006efficient, le2010reachability, althoff2016combining, bogomolov2018reach} and non-linear \citep{henzinger1998algorithmic, asarin2003reachability, althoff2008reachability, chen2013flow} differential equations.
{\rev These algorithms are implemented in numerical tools such as SpaceEx \citep{frehse2011spaceex}, Flowstar \citep{chen2013flow}, CORA \citep{althoff2015introduction} or JuliaReach \citep{juliareach}.} %
Despite the existence of numerous algorithms, their successful application requires expert knowledge about algorithm tuning \citep{wetzlinger2020adaptive}.

The main motivation of our work is to integrate reachability analysis and finite element methods in structural analysis applications.
{\rev Recent specific works have initiated research in such direction, as in \citep{scacchioli2014assessment}, where the} authors use ellipsoidal set representations{\rev. %
 This approach does not scale } for systems with many degrees of freedom. In \citep{BakTJ19} the authors develop a set propagation method based on linear programming that is applied to a heat diffusion problem with a large ($1000^3$) mesh obtained using finite differences. In comparison with our technique, theirs' returns reach-set approximations at time points, instead of reach-sets at time intervals as in our method.
In \citep{althoff2019reachability} the author extended the theoretical basis of set propagation to handle large linear systems with Krylov subspace approximations of the state transition matrix acting on zonotopes.
Their results are illustrated on the finite element model of a bridge subject to set-valued wind forces.

In this article we present a {\rev novel} framework for the resolution of linear transient dynamics problems in solid mechanics, that builds upon set propagation techniques.
Starting from the assembled FEM system of differential equations, our method is based on {\rev the construction of a set including the solutions for the first time step and propagating the set to obtain an enclosure of the solutions for posterior times.}
{\rev 
The methodology was implemented using JuliaReach\footnote{\href{http://juliareach.com}{juliareach.com}}, a state-of-the-art open source library for reachability analysis. Finite-element assembly was done using ONSAS\footnote{\href{https://github.com/ONSAS}{github.com/ONSAS}}.
}

This article is organized as follows.
In Section~\ref{sec:preliminaries} the basic concepts of the Finite Element Method and Reachability Analysis are described, including an illustrative example complementing the presentation of the latter.
In Section~\ref{sec:methodology} the key components of the proposed methodology are developed.
In Section~\ref{sec:results} the results obtained considering five numerical examples (including reference examples from the literature and a realistic Engineering problem) are shown. %
Finally, in Section \ref{sec:conclusions} the conclusions are exposed and additional proofs are available as part of the Appendix.

\section{Preliminaries} \label{sec:preliminaries}

%
%

%
%

{\rev 
\subsection{ Transient problems governing equations} \label{sec:transient_problems}
}

In this section we recall the basic equations of the heat transfer and structural dynamics problems, and omit their development which can be found in the literature \citep{Wilson1974,Hughes1987a,Bathe2014}.

For the heat transfer problem, considering an isotropic solid in the region $\Omega$ with conductivity $\kappa$, density $\rho$ and specific heat $c$, the governeing equations are
\begin{equation}\label{eqn:heatfem}
	\bfC_\theta \dot{\bstheta}(t) + \bfK_\theta \bstheta(t) = \bff_\theta(t),
\end{equation}
where the matrix $\bfK_\theta$ is the assembled convection matrix
\begin{equation}\label{eqn:heatfem_convection}
	\bfK_\theta = \Ass_{e=1}^{n_e} \bfK_\theta^e,  \qquad 
	\bfK_\theta^e = \int_{\Omega^e} \kappa \left( \bfB^e \right)^T \bfB^e \dif V + \int_{\Gamma_{R}^e} h \left(\bfN^e\right)^T \bfN^e \dif S,
\end{equation}
with $n_e$ being the number of elements and $\bfN$ and $\bfB$ the matrices of interpolation functions and their derivatives, respectively, the matrix $\bfC_\theta$ is the assembled diffusion matrix
\begin{equation}\label{eqn:heatfem_diffusion}
	\bfC_\theta = \Ass_{e=1}^{n_e} \bfC_\theta^e  \qquad 
	\bfC_\theta^e = \int_{\Omega^e} \rho c \left(\bfN^e\right)^T \bfN^e \dif V,
\end{equation}
and the vector of heat sources $\bff_\theta(t)$ is given by
\begin{equation}\label{eqn:heatfem_force}
	\bff_\theta(t)
	=
	\Ass_{e=1}^{n_e} \int_{\Omega^e} Q_{int}(t) \bfN^e \dif V
	+ \int_{\Gamma_N^e} q_{inp}(t)  \bfN^e \dif S 
	+ \int_{\Gamma_R^e} h \theta_\infty(t) \bfN^e \dif S ,
\end{equation} 
with $Q_{int}$ and $q_{inp}$ being volumetric and boundary external heat sources, $h$ being a convection constant and $\theta_\infty$ being the external temperature.

For the structural dynamics problem, considering a solid with density $\rho$ and linear elastic constitutive behavior given by the tensor $\bbC$, the governing equations can be written as follows
\begin{equation}\label{eqn:dynamicsfem}
	\bfM_\bfu \ddot{\bfu}(t)+ \bfC_\bfu \dot{\bfu}(t) + \bfK_\bfu \bfu(t) = \bff_{\bfu}(t),
\end{equation}
where $\bfK_\bfu$ is the stiffness matrix, given by:
\begin{equation}
	\bfK_\bfu = %
	\Ass_{e=1}^{n_e} \bfK_\bfu^e,  \qquad 
	\bfK_\bfu^e = \int_{\Omega^e} \left( \bfB_u^e \right)^T   \bbC \bfB_u^e \, \dif V,
\end{equation}
and $\bfM_\bfu$ is the mass matrix, given by:
\begin{equation}
	\bfM_\bfu = %
	\Ass_{e=1}^{n_e} \bfM_\bfu^e,  \qquad 
	\bfM_\bfu^e = \int_{\Omega^e} \rho  (\bfN_u^e)^T 
	\bfN_u^e \, \dif V,
\end{equation}
where the matrices $\bfN_u^e$ and $\bfB_u^e$ are the matrices given by the interpolation functions and gradients for the displacement fields, and $\bff_{\bfu}$ is the vector of applied loads. %
The viscosity (or damping) matrix $\bfC_\bfu$ is assumed to be known or estimated as a combination of the stiffness and the mass matrices \citep{clough1993dynamics}.

{\rev
\subsection{Numerical Integration Methods} \label{sec:numerical_integration}
}

In this section, the equations of three reference numerical time-integration methods \citep{Bathe2014} are recalled. %
A fixed time step $\Delta t$ is considered, and the time domain is discretized in values $t_k = \Delta t \cdot k$, $k \in \mathbb{N}$.

The first reference numerical integration technique is the Backward Euler Method, which is used to solve the heat transfer equations. %
Considering an initial temperature $\bstheta_0$, the method consists in solving the following linear system for each time step:
\begin{equation}
	\left( \bfK_\theta \Delta t + \bfC_\theta \right) \bstheta_{k+1} =
	\bff_{\theta,{k+1}} \Delta t + \bfC_\theta \bstheta_k.
\end{equation}

Regarding structural dynamics problems, the Newmark and the Bathe methods are presented.
Assuming that initial displacements and velocities are given, after computing the initial acceleration, each method provides a different procedure to obtain subsequent displacements and velocities. The Newmark method provides the following system:
\begin{equation}
	( b_0 \bfM_\bfu + b_1 \bfC_\bfu + \bfK_\bfu ) \bfu_{k+1} %
	= \bff_{\bfu,k+1} %
	+ \bfM_\bfu ( b_0 \bfu_k + b_2 \dot{\bfu}_k
	+ \ddot{\bfu}_k ) %
	+ \bfC_\bfu (b_1 \bfu_k + \dot{\bfu}_k ),
\end{equation}
where the auxiliary constants $b_0 = 4/\Delta t^2$, $b_1 = 2/\Delta t$ and $b_2 = 4/\Delta t$ are introduced. Update rules are provided for the velocities and accelerations.

In the Bathe method, the displacements are computed in two sub-steps, with corresponding update rules for the velocities. %
The linear systems corresponding to these sub-steps are the following:
\begin{eqnarray}
	\left( a_0 \bfM_\bfu + a_1 \bfC_\bfu + \bfK_\bfu \right) \bfu_{(k+1/2)} %
	&=& \bff_{\bfu,(k+1/2)} %
	+ \bfM_\bfu ( a_0 \bfu_k + a_4 \dot{\bfu}_k
	+ \ddot{\bfu}_k ) %
	+ \bfC_\bfu ( a_1 \bfu_k + \dot{\bfu}_k ), \\
	( a_2 \bfM_\bfu + a_3 \bfC_\bfu + \bfK_\bfu ) \bfu_{k+1} %
&=& \bff_{\bfu,k+1} + \bfM_\bfu ( a_5 \bfu_{(k+1/2)} + a_6 \bfu_k + a_1 \dot{\bfu}_{(k+1/2)}
+ a_7 \dot{\bfu}_k )  \nonumber\\
& \dots & + \,\, \bfC_u ( a_1 \bfu_{(k+1/2)} + a_7 \bfu_k ),
\end{eqnarray}
where the auxiliary constants $a_0 = 16/\Delta t^2$, $a_1 = 4/\Delta t$, $a_2 = 9/\Delta t^2$, $a_3 = 3 / \Delta t$, $a_4 = 8/\Delta t$, $a_5 = 12 / \Delta t^2$, $a_6 = -3 / \Delta t^2$ and $a_7 = -1 / \Delta t$ were introduced.

\subsection{Set propagation concepts} \label{sec:set_propagation_preliminaries}

In this section we consider set-based integration in the context of linear differential equations and introduce two essential concepts: reach-set and flowpipe. Then we define the representations used in our approach: hyperrectangles, zonotopes and support functions. The definitions are illustrated by means of a textbook example: the harmonic oscillator. More complex and larger problems are considered in Section~\ref{sec:results}.

\subsubsection{Illustrative example: harmonic oscillator}\label{sec:illustrative_example}

The set propagation definitions introduced in this section are illustrated using a single degree of freedom second order problem given by the following equation:
\begin{equation} \label{eq:harmonic_oscillator}
	\centering
	\ddot{u}(t) + \omega^2 u(t) = 0,
\end{equation}
where $\omega$ is a scalar and $u$ is the unknown. %
This problem can be associated with a spring-mass system, where $u(t)$ is the elongation of the spring at time $t$ and the mass and stiffness set a natural frequency $\omega$. In this case we consider $\omega = 4 \pi$.

Let us introduce the state variable $v(t) := \dot{u}(t)$ and define the vector $\bfx(t) = [ u(t), v(t) ]^T$. Then Eq.~\eqref{eq:harmonic_oscillator} can be written in the following first order form
\begin{equation} \label{eq:harmonic_oscillator_first_order}
	\dot{\bfx}(t) = 
	\left[
	\begin{matrix}
		0 & 1 \\
		-\omega^2 & 0
	\end{matrix}
	\right] \, \bfx(t).
\end{equation}

\begin{figure}[htbp]
	\centering
	\begin{subfigure}[b]{0.7\textwidth}
		\centering
		\includegraphics[width=\textwidth]{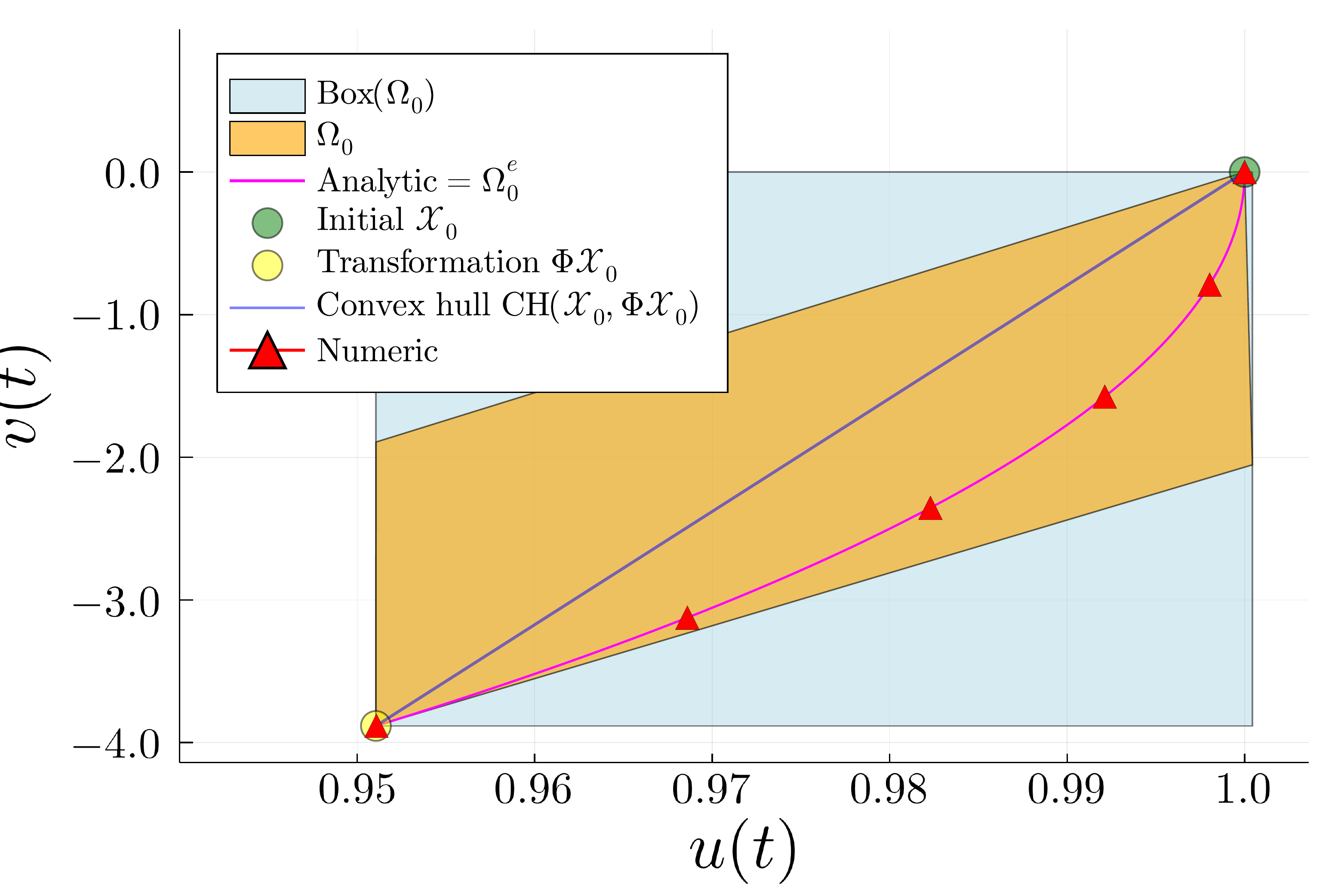}
		\caption{The set $\Omega_0$ (orange) encloses the true solution (magenta) at the endpoints and at any intermediate time between $0$ and $\delta = 0.025$.}
		\label{fig:sdof_singleton_discretization}
	\end{subfigure}
	~~
	\begin{subfigure}[b]{0.7\textwidth}
		\centering
		\includegraphics[width=\textwidth]{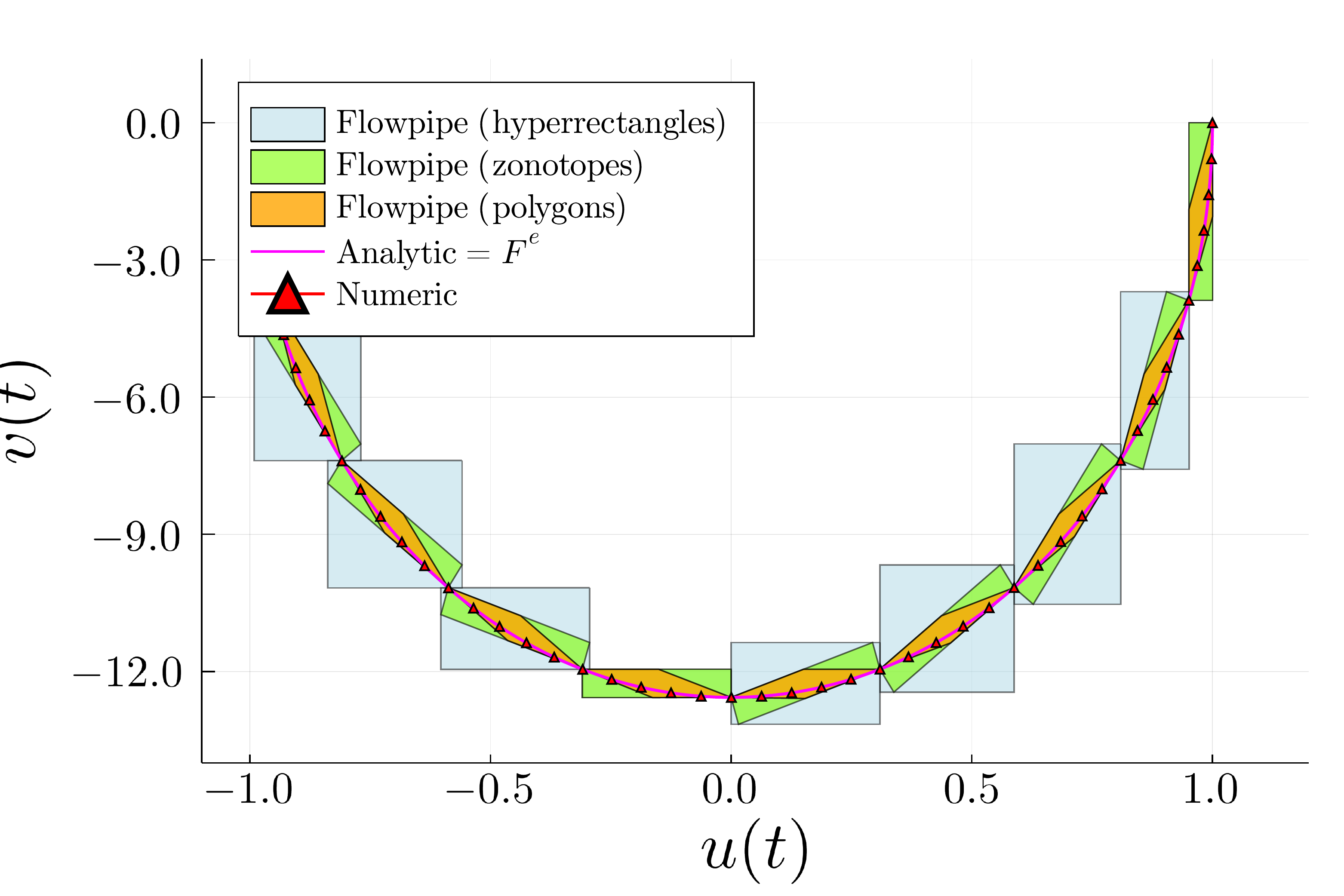}
		\caption{Set propagation using hyperrectangles (lightblue), zonotopes (green) and to assess the zonotopic overapproximation of $\Omega_0$, polygons (orange).}
		\label{fig:sdof_singleton_flowpipe}
	\end{subfigure}
	\caption{Illustration of set-based integration of the simple harmonic oscillator example with single initial condition. Shown are phase-space plots $u(t)$ vs. $v(t)$ for the first reach-set (Fig. \ref{fig:sdof_singleton_discretization}), and for the first 9 reach-sets (Fig. \ref{fig:sdof_singleton_flowpipe}).}
	\label{fig:sdof_singleton}
\end{figure}

Two cases are considered for the initial conditions, first a single initial condition given by $u(0)=1$ and $v(0)=0$, and a second case where the initial conditions belong to a set. %
This problem is solved using set propagation and the results obtained for the single initial condition are shown in Fig.~\ref{fig:sdof_singleton}.
The time step-size used is $\delta = 0.025$.
The set $\Omega_0$ from Fig.~\ref{fig:sdof_singleton_discretization} includes the analytic trajectory within $[0, \delta]$, and such set is propagated to cover the solution for further time intervals as in Fig.~\ref{fig:sdof_singleton_flowpipe}.
It is worth noting that even when the initial condition is a singleton, the set propagation technique returns a sequence of sets, and such sequence guarantees an enclosure of the analytic solution at \textit{any} time $t \in [0, T] \subset \mathbb{R}$.
For comparison, we also plot the analytic solution of Eq.~\eqref{eq:harmonic_oscillator_first_order} (magenta solid line) and a numerical solution obtained using the Newmark's method with time step $\delta/5$ (red triangle markers).
The figure shows other set representations and concepts that will be developed in the following paragraphs.

\begin{figure}[htbp]
	\centering
	\begin{subfigure}[b]{0.7\textwidth}
		\centering
		\includegraphics[width=\textwidth]{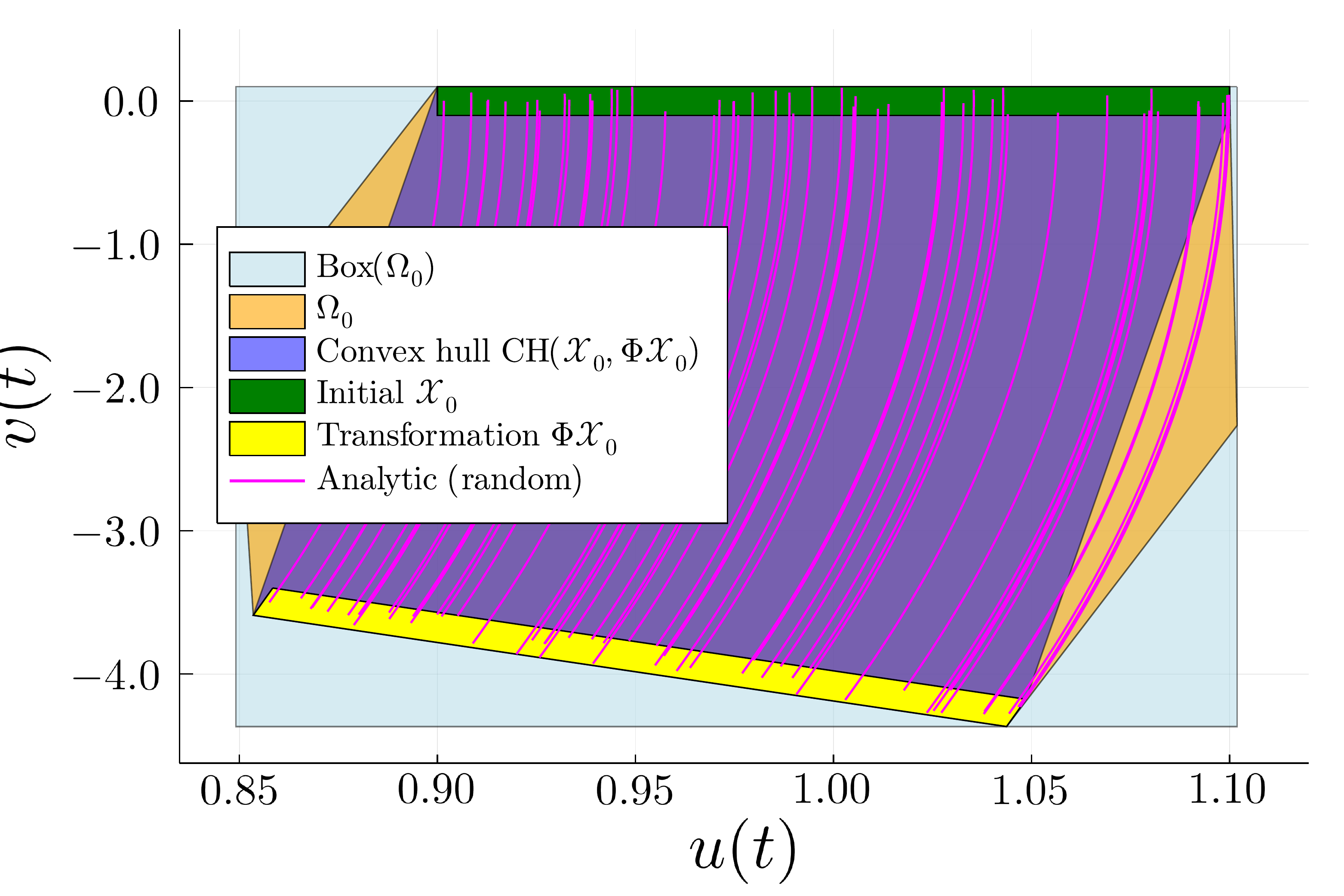}
		\caption{The set $\Omega_0$ (orange) covers the right-most trajectories for intermediate times, which naturally escape linear interpolations (violet).}
		\label{fig:sdof_distributed_discretization}
	\end{subfigure}
	~~
	\begin{subfigure}[b]{0.7\textwidth}
		\centering
		\includegraphics[width=\textwidth]{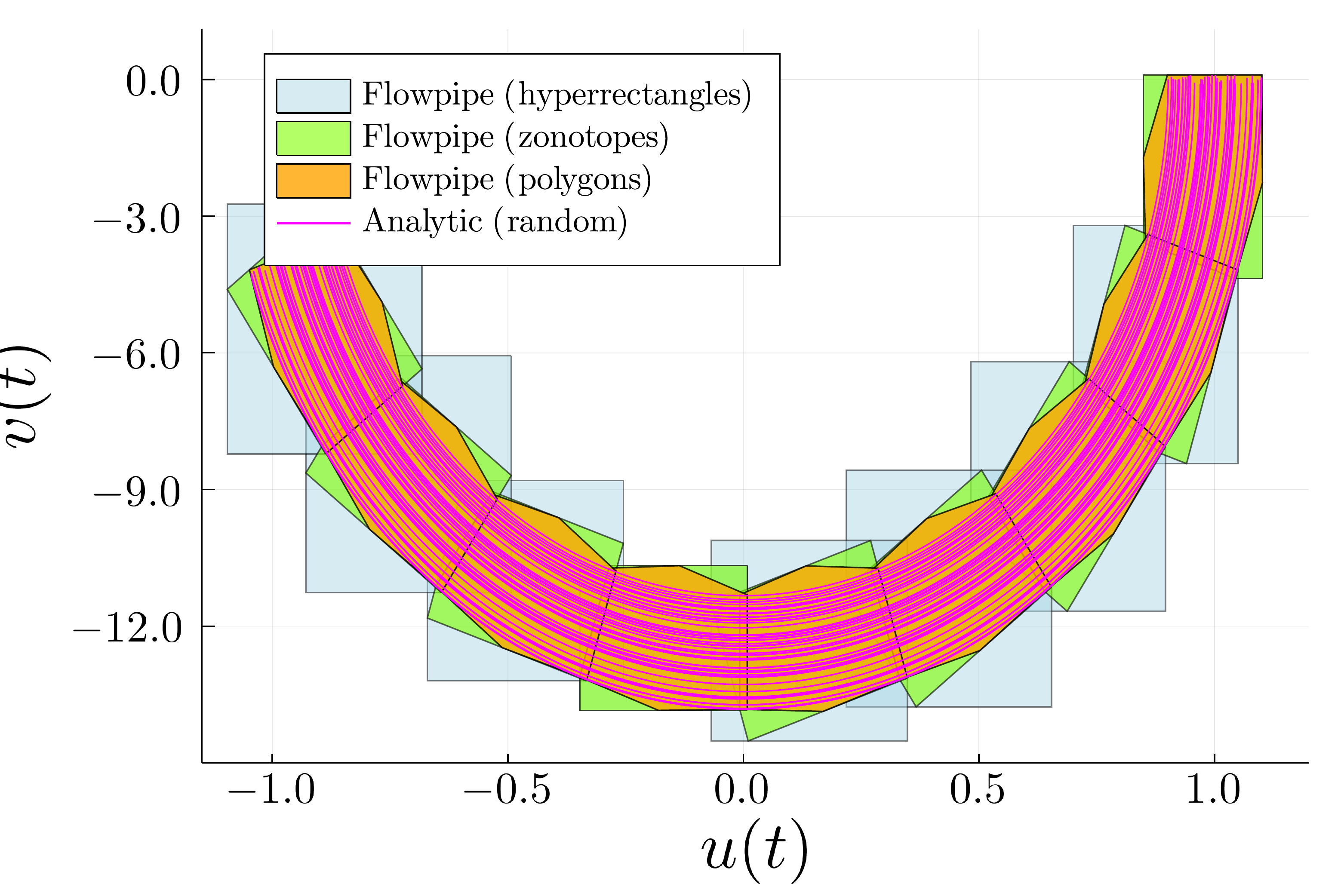}
		\caption{The flowpipe construction method for Eq.~\eqref{eq:discrete_recurrence} is wrapping-free: the area of the sets does not increase with time.}
		\label{fig:sdof_distributed_flowpipe}
	\end{subfigure}
	\caption{Illustration of set-based integration of the simple harmonic oscillator with distributed initial conditions. Note that a single set-based integration covers all dynamically feasible behaviors: to illustrate this fact we plot 50 analytic solutions (magenta) uniformly distributed over $\mathcal{X}_0$.}
	\label{fig:sdof_distributed}
\end{figure}

The result obtained when a set of initial conditions is considered is shown in Fig.~\ref{fig:sdof_distributed}, together with a few dozen trajectories with random initial conditions drawn from the initial set $\mathcal{X}_0 = \mathcal{U}_0 \times \mathcal{V}_0 = [0.9, 1.1] \times [-0.1, 0.1]$, where $\times$ denotes the Cartesian product.

In both cases, it is shown that a single set-based integration covers infinitely many trajectories \textit{in dense time}, i.e. for all intermediate times between $0$ and $T$, where $T > 0$ is the time horizon, there is a reach-set that covers all the exact solutions.

\subsubsection{Trajectories, reach-sets and flowpipes} \label{ssec:reachsets}

Let us consider an initial-value problem
\begin{equation}\label{eq:linearODE}
	\dot{\bfx}(t) = \bfA \bfx(t),\qquad \bfx(0) \in \mcX_0,~t \in [0, T],
\end{equation}
with state matrix $\bfA \in \mathbb{R}^{n\times n}$ and where $\bfx(t) \in \mathbb{R}^n$ is the state vector at time $t$.
The initial condition is $\bfx(0) \in \mcX_0\subset \mathbb{R}^n$ with $\mcX_0$ a closed and convex set. In this article we will refer to any state belonging to $\mcX_0$ as an \textit{admissible} state.
The guiding principle in set-based integration is to represent reachable states using sets in the Euclidean space that are \textit{propagated} according to the system's dynamics.
The obtained sequence of sets, by construction, covers the exact trajectories for given time intervals, for any admissible initial condition or external actions. For instance, the result for the harmonic oscillator presented in Section~\ref{sec:illustrative_example} using a single initial condition is shown in Fig.~\ref{fig:sdof_singleton} (green circle), and in Fig.~\ref{fig:sdof_distributed} (green box) the results for a box of initial conditions are shown.

Given the system of ODEs in Eq.~\eqref{eq:linearODE}, a \textit{trajectory} from an initial state $\bfx_0$ is the solution for a given initial condition $\bfx(0) = \bfx_0$, denoted $\varphi(\bfx_0, t) : \R^n \times [0, T] \rightarrow \R^n$, which we know is $\varphi(\bfx_0, t) = e^{\bfA t}\bfx_0$.
The \textit{reach-set} at a \textit{time point} $t \in [0, T]$ is
\begin{equation}\label{eq:reachset_point}
	\mcR^e(\mcX_0, t) := \bigl\{ \varphi(\bfx_0, t): \bfx_0 \in \mcX_0 \bigr\}.
\end{equation}
This is the set of states which are reachable at time $t$ starting from any admissible initial condition.
Generalizing Eq.~\eqref{eq:reachset_point} to time intervals is straightforward: the reach-set over $[t_0, t] \subseteq [0, T]$ is
\begin{equation}\label{eq:reachset_interval}
	\mcR^e(\mcX_0, [t_0, t]) := \bigl\{ \varphi(\bfx_0, s): \bfx_0 \in \mcX_0, \ \forall~s \in [t_0, t] \bigr\}.
\end{equation}
The difference between reach-sets for time points and for time intervals is that the former are evaluated at a particular time $t$, while the latter include \textit{all} states reachable for \textit{any} time between $t_0$ and $t$, hence $\mcR^e(\mcX_0, t) \subseteq \mcR^e(\mcX_0, [t_0, t])$. 
Finally, given a time step $\delta > 0$ and a collection of reach-sets $\left\{\mcR^e(\mcX_0, [\delta(k-1), \delta k)]\right\}_{k=1}^N$ with $N = T/\delta$, the \textit{flowpipe} of Eq.~\eqref{eq:linearODE} is the set union,
\begin{equation}\label{eq:flowpipe}
	\mcF^e(\mcX_0, [0, T]) := \bigcup_{k=1}^N \mcR^e(\mcX_0, [\delta(k-1), \delta k)]).
\end{equation}

In Eqs.~\eqref{eq:reachset_point}-\eqref{eq:flowpipe}, the superindex $e$ is used to remark that the definitions refer to the \textit{exact} or true sets. In practice, reachable sets can only be obtained approximately and several algorithms to numerically approximate reachable sets are known in the literature. Different methods crucially depend on how the sets are represented and the cost of operating with those representations.

\subsubsection{Set representations}

Three set representations commonly used in reachability analysis of linear systems are zonotopes, hyperrectangles and support functions. Let us briefly recall their definition.

\paragraph{Zonotopes}

Given a vector $\bfc \in \mathbb{R}^n$ that we call \textit{center} and a list of $p \ge n$ vectors that we call \textit{generators}, $\bfg_1, \ldots, \bfg_p$, $\bfg_i\in \mathbb{R}^n$, the associated zonotope is the set
\begin{equation} \label{eq:zonotope}
	Z = \left\{\bfx \in \mathbb{R}^n : \bfx = \bfc + \sum_{i=1}^p \xi_i \bfg_i,\quad \xi_i \in [-1, 1] \right\}.
\end{equation}
In compact notation we write $Z = \langle \bfc, \bfG\rangle_Z$, where each generator is stored as a column of the generators matrix $\bfG \in \mathbb{R}^{n \times p}$. Zonotopes can be interpreted in different ways, such as the Minkowski sum of line segments, or as the affine map of a unit ball in the infinity norm. The zonotope representation has been successfully used in reachability analysis because it offers a very compact representation relative to the number of vertices or faces, so it is well suited to represent high dimensional sets \citep{GLGM06, althoff2016combining}.

\paragraph{Hyperrectangles}

Given vectors $\bfc \in \mathbb{R}^n$ and $\bfr \in \mathbb{R}^n$ that we call \textit{center} and \textit{radius} respectively, the associated hyperrectangle is the set
\begin{equation} \label{def:hyperrectangle}
	H = \left\{\bfx \in \mathbb{R}^n : \vert x_i - c_i \vert \leq r_i, \quad i = 1, \ldots, n\right\}.
\end{equation}
In compact notation we write $H = \langle \bfc, \bfr \rangle_H$. Every hyperrectangle is a zonotope but the converse is not true. In fact, the linear map of a hyperrectangle can be (exactly) represented as a zonotope.
While computations with hyperrectangles can be performed very efficiently, using hyperrectangles for reachability can result in coarse overapproximations unless they are used in concurrence with other representations at intermediate parts of the algorithm. 

\paragraph{Support functions}

The support function of a set $X \subset \mathbb{R}^n$ along direction $\bfd \in \mathbb{R}^n$ is the scalar
\begin{equation}\label{eq:support_function}
	\rho(\bfd, X) = \max_{\bfx \in X} \bfd^T \bfx,
\end{equation}
where the superscript $T$ denotes transposition. 
The set of points in $X$ that are the maximizers of Eq.~\eqref{eq:support_function} are called the \textit{support vectors}. The intuition behind Eq.~\eqref{eq:support_function} is that support functions represent the farthest (signed) distance to the origin of the set $X$ along direction $\bfd$. Given that the direction can be chosen at will, support functions can be used to obtain the solution corresponding to a linear combination of the state variables at once and at a reduced cost~\citep{LeGuernic09}.


\section{Methodology} \label{sec:methodology}

%
%
%
%

{\rev
In this section the proposed method is presented.}
{\rev Given a heat transfer or structural dynamics system of linear differential equations, the proposed algorithm can be summarized in the following three stages:
\begin{itemize}
	\item[(i)]  \textbf{Homogeneization and transformation into a first order system.} Transform the system {\rev of differential equations  into a linear homogeneous system of first order,} as in Eq.~\eqref{eq:linearODE}, using the method described in Section~\ref{ssec:homogeneization}.
	\item[(ii)] \textbf{Conservative time discretization.} Build a set $\Omega_0 \subset \mathbb{R}^n$ that contains all exact trajectories between time $0$ and time $\delta > 0$, i.e. satisfying the property:
	\begin{equation}\label{eq:discretization}
		\Omega_0^e := \mcR^e(\mcX_0, [0, \delta]) \subseteq \Omega_0,
	\end{equation}
	{\rev 
		using the results presented in Section~\ref{ssec:discretization}.}
	\item[(iii)] \textbf{Set propagation.} Set the initial set $X_0 := \Omega_0$ and the state transition matrix $\bfPhi := e^{\bfA\delta}$, and {\rev compute the reachable sets for the subsequent time intervals, solving the} recurrence{\rev :}
	\begin{equation}\label{eq:discrete_recurrence}
		X_{k+1} = \bfPhi X_k,\qquad X_0 \subset \mathbb{R}^n,~k = 0, 1,\ldots, N-2,
	\end{equation}
	{\rev where $X_k$ is an overapproximation of the exact reachable sets for the time interval $[k\delta, (k+1) \delta]$. This is done using available set propagation strategies described in Section~\ref{ssec:set_propagation}.}
\end{itemize}

{\rev The proposed method consists in an explicit numerical approach, since Eq.~\eqref{eq:discrete_recurrence} is used to calculate the reachable set at a later time interval from the reachable set at the current time interval.
If the state dimension is large, Krylov methods are used to evaluate the action of the matrix exponential on specific directions, instead of computing the matrix exponential in full.
}


\subsection{Homogeneization and transformation into a first order system} \label{ssec:homogeneization}

We now describe {\rev the first step of the method, consisting in:} an algebraic transformation of the governing equations {\rev (heat transfer or structural dynamics)}, into a first order, homogeneous system. %
%
%

{\rev 
\subsubsection{Transformation into a first order form}

The first step of the method is to reformulate the governing equations as a system of the form of Eq.~\eqref{eq:linearODE}. %
In the case of heat transfer problems, using that the diffusion matrix $\bfC_\theta$ is a square positive definite (thus invertible) matrix, the system given by Eq.~\eqref{eqn:heatfem} can be written as
\begin{equation}\label{eqn:heatTrans}
	\dot{\bstheta}(t) = - \bfC_\theta^{-1}  \bfK_\theta \bstheta(t) + \bfC_\theta^{-1} \bff_\theta(t).
\end{equation}
For the structural dynamics equations, using that the matrix $\bfM_\bfu$ is positive definite (and invertible), the system given by Eq.~\eqref{eqn:dynamicsfem} can be written as:
\begin{equation}\label{eqn:dynTrans}
 \left\{
 \begin{array}{l}
\dot{\bfu} = \bfv \\
\dot{\bfv} = 	 -\bfM_\bfu^{-1} \bfK_\bfu \bfu   -\bfM_\bfu^{-1} \bfC_\bfu \bfv   +\bfM_\bfu^{-1} \bff_{\bfu}(t) 
 \end{array}
 \right.
\end{equation}
}

\subsubsection{Homogeneization}\label{sec:homogeneization}

{\rev In this stage of the method, the system of equations is converted to a homogeneous form. In order to do that, the inputs vector is integrated to the system of ODEs as a combination of auxiliary variables with a given dynamic. }

Let us consider that the inputs vector  $\bff(t)$ {\rev (representing loads $\bff_{\bfu}$ or heat fluxes $\bff_\theta$)} can be decomposed as follows:
\begin{equation}\label{eqn:decompEta}
	\bff(t) = \sum_{i=1}^{n_f} \bff_0^{(i)} \eta^{(i)}(t),
\end{equation}
where $n_f$ is the number of input components, each $\bff_0^{(i)} \in \mathbb{R}^m$ is a known vector and the functions $\eta^{(i)}(t):[0,T]\rightarrow \bbR$ are
functions that can be expressed in terms of the solution of a first order homogeneous system of differential equations.

The property of functions $\eta^{(i)}$ can be written as follows: for each function $\eta^{(i)}$ there exists a vector of $q^{(i)}$ functions $\bfxi^{(i)}(t) = [\xi_1^{(i)}(t), \ldots, \xi_{q^{(i)}}^{(i)}(t)]^T$, such that, for each $i = 1, \ldots, n_f$:
\begin{equation}\label{eqn:etaODE}
	\left\{
	\begin{array}{l}
		\eta^{(i)}(t) = \xi_1^{(i)}(t)\\
		\dot{\bfxi}^{(i)}(t) = \bfB^{(i)} \bfxi^{(i)}(t)
	\end{array}
	\right.
\end{equation}
where $\bfB^{(i)}$ is a $q^{(i)}\times q^{(i)}$ matrix. %
To exemplify definition \eqref{eqn:etaODE}, two special cases are mentioned: exponential functions (with $q = 1$) and trigonometric functions (with $q = 2$){\rev, which correspond to systems
\begin{equation} \label{eq:input_function_cases}
\left\{
\begin{array}{l}
\eta(t) = \xi_1(t)\\
\dot{\xi_1}(t) = \alpha \cdot \xi_1(t)\\
\end{array}
\right.
\quad
\text{and}
\quad
\left\{
\begin{array}{l}
	\eta(t) = \xi_1(t)\\
	\dot{\xi_1}(t) = \xi_2(t)\\
	\dot{\xi_2}(t) = -\omega^2 \xi_1(t)\\
\end{array}
\right.,
\end{equation}
respectively, with $\alpha$ and $\omega$ being real scalars. In particular, the case $\alpha = 0$ is used to model constant functions.}

Substituting Eq.~\eqref{eqn:etaODE} into Eq.~\eqref{eqn:decompEta}, we can reformulate Eq.~\eqref{eqn:heatTrans} as
\begin{equation}\label{eqn:homheat}
	\left\{
	\begin{array}{l}
		{\rev
		\dot{\bstheta}(t) = - \bfC_\theta^{-1}  \bfK_\theta \bstheta (t) +  \sum_{i=1}^{n_f} \bfC_\theta^{-1}  \bfF_0^{(i)} \bfxi^{(i)}(t) }\\
		\dot{\bfxi}^{(i)}(t) = \bfB^{(i)} \bfxi^{(i)}(t) \quad i=1,\dots,n_f
	\end{array}
	\right.
\end{equation}
{\rev
where 
$$
\bfF_0^{(i)} = \begin{bmatrix}
	\begin{array}{c|c}
		\bff_0^{(i)} & {\bf0}_{m\times (q^{(i)}-1)}
	\end{array}
\end{bmatrix}.
$$
}
Similarly, Eq.~\eqref{eqn:dynTrans} {\rev can be rewritten} as %
\begin{equation}\label{eqn:homdyn}
	\left\{
	\begin{array}{l}
		{\rev 
    	\dot{\bfu} = \bfv }\\
    	{\rev 
    	\dot{\bfv} = 	 -\bfM_\bfu^{-1} \bfK_\bfu \bfu   -\bfM_\bfu^{-1} \bfC_\bfu \bfv   + \sum_{i=1}^{n_f} \bfM_\bfu^{-1} \bfF_0^{(i)} \bfxi^{(i)}(t) } \\
		\dot{\bfxi}^{(i)}(t) = \bfB^{(i)} \bfxi^{(i)}(t) \quad i=1,\dots,n_f
	\end{array}
	\right.
\end{equation}
obtaining then an homogeneous system of differential equations with $m+\sum_{i=1}^{n_f} q^{(i)}$ unknowns. %

\subsubsection{{\rev Matrix Coupled formulation and initial conditions} } \label{sec:transformation_firstorder}

{\rev
The systems of differential equations in Eq.~\eqref{eqn:homheat} and Eq.~\eqref{eqn:homdyn},} can be written in a first order matrix coupled form as follows:
\begin{equation} \label{eqn:coupled_formulation}
	\dfrac{d}{dt}\begin{bmatrix}
		\bfx(t) \\
		\bfxi^{(1)}(t)\\
		\vdots \\
		\bfxi^{(n_f)}(t)
	\end{bmatrix} = \bfA
	\begin{bmatrix}
		\bfx(t) \\
		\bfxi^{(1)}(t)\\
		\vdots \\
		\bfxi^{(n_f)}(t)
	\end{bmatrix},
\end{equation}
where the {\rev vector of} unknowns $\bfx(t)$ {\rev contains the \textit{state variables} ($\bstheta$ for heat transfer or $[\bfu, \bfv]$ for structural dynamics)}, the {\rev unknowns in the vector}  $\bfxi^{(i)}(t)$ are called \textit{input variables} {\rev  and the matrix $\bfA$ is given by
	\begin{equation} \label{eq:homogeneized_matrix_heat}
	\bfA = 
	\left[
	\begin{array}{c|c;{2pt/2pt}c;{2pt/2pt}c}
		-\bfC_\theta^{-1}\bfK_\theta & \bfC_\theta^{-1}\bfF_0^{(1)} & \dots & \bfC_\theta^{-1}\bfF_0^{(n_f)} \\ \hline
		\bf0 & \bfB^{(1)} & \bf0 & \bf0 \\ \hdashline[2pt/2pt]
		\vdots & \bf0 & \ddots & \bf0 \\ \hdashline[2pt/2pt]
		\bf0 & \dots & \bf0 & \bfB^{(n_f)}
	\end{array}
	\right],
\end{equation}
for heat transfer problems or
\begin{equation} \label{eq:homogeneized_matrix_dynamics}
	\bfA = 
	\left[
	\begin{array}{c;{2pt/2pt}c|c;{2pt/2pt}c;{2pt/2pt}c}
		\bf0 & \bfI & \bf0 & \dots & \bf0 \\ \hdashline[2pt/2pt]
		-\bfM_\bfu^{-1}\bfK_\bfu & -\bfM_\bfu^{-1}\bfC_\bfu & \bfM_\bfu^{-1}\bfF_0^{(1)} & \dots & \bfM_\bfu^{-1}\bfF_0^{(n_f)} \\ \hline
		\bf0 &\bf0 & \bfB^{(1)} & \bf0 & \bf0 \\ \hdashline[2pt/2pt]
		\vdots & \vdots &\bf0 & \ddots & \bf0 \\ \hdashline[2pt/2pt]
		\bf0 & \bf0 & \dots & \bf0 & \bfB^{(n_f)}
	\end{array}
	\right],
\end{equation}
for structural dynamics problems.	
}

{\rev The initial conditions of the system given in Eq.~\eqref{eqn:coupled_formulation}  can be defined by  the Cartesian product $\mathcal{X}_0\times \mathcal{C}_0$, where $\mathcal{X}_0$ is the set of initial conditions of the state variables  and $\mathcal{C}_0$ is the set of initial conditions for the differential equations in Eq.~\eqref{eqn:etaODE}.
}

Concerning the initial condition for the input variables in Eq.~\eqref{eqn:coupled_formulation}, that we call $\mathcal{C}_0$, they are not restricted to be a single value; in fact they can be chosen to be e.g. a hyperrectangle, in which case we can model an uncertainty of the initial conditions around nominal values for a whole family of input functions, with a single flowpipe computation.

\begin{figure}[htb!]
	\centering
	\begin{subfigure}[b]{0.48\textwidth}
		\includegraphics[width=\textwidth]{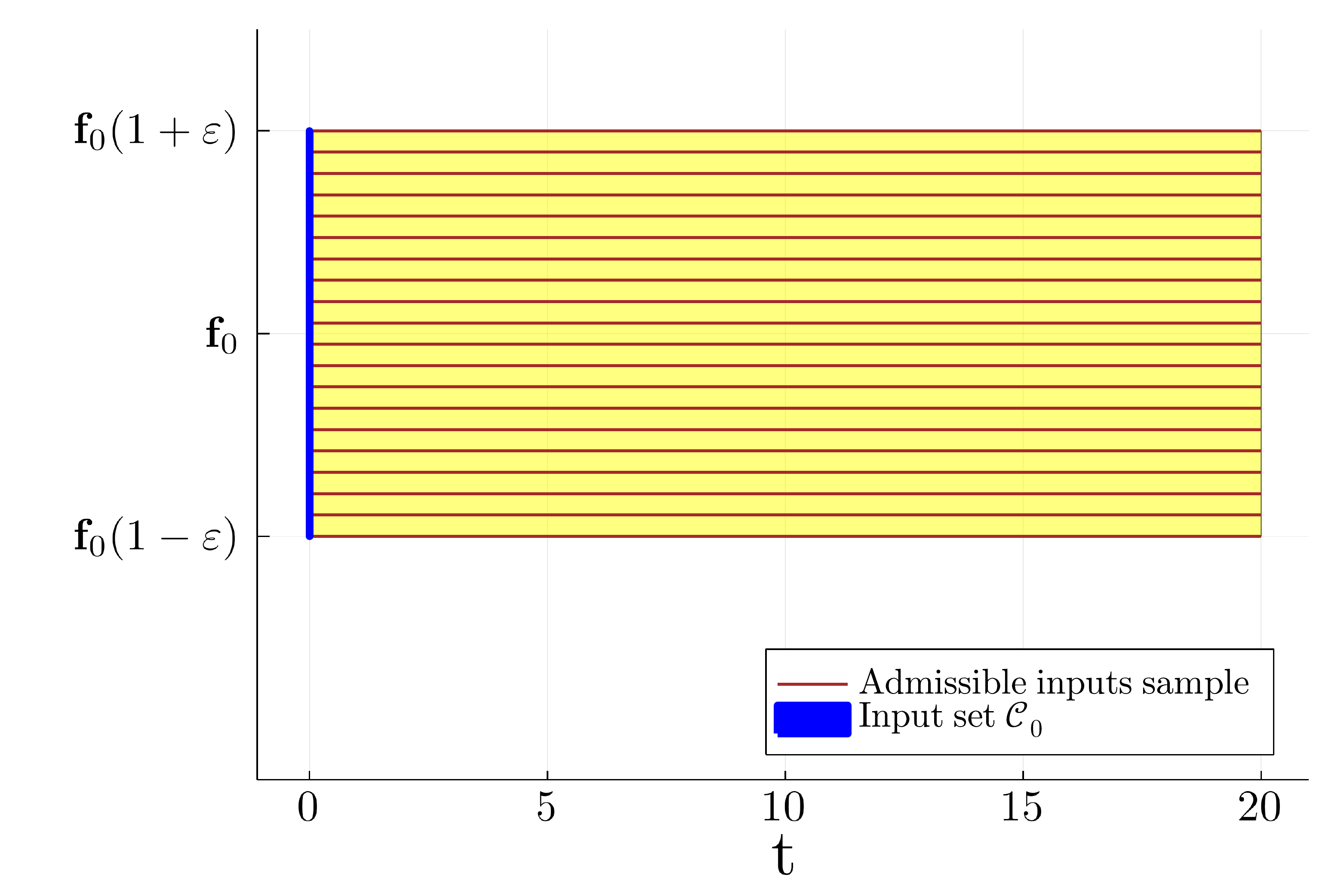}
		\caption{Example of constant input model.}
		\label{fig:constantforcesRA}
	\end{subfigure}
	\hfill
	\begin{subfigure}[b]{0.50\textwidth}
		\includegraphics[width=\textwidth]{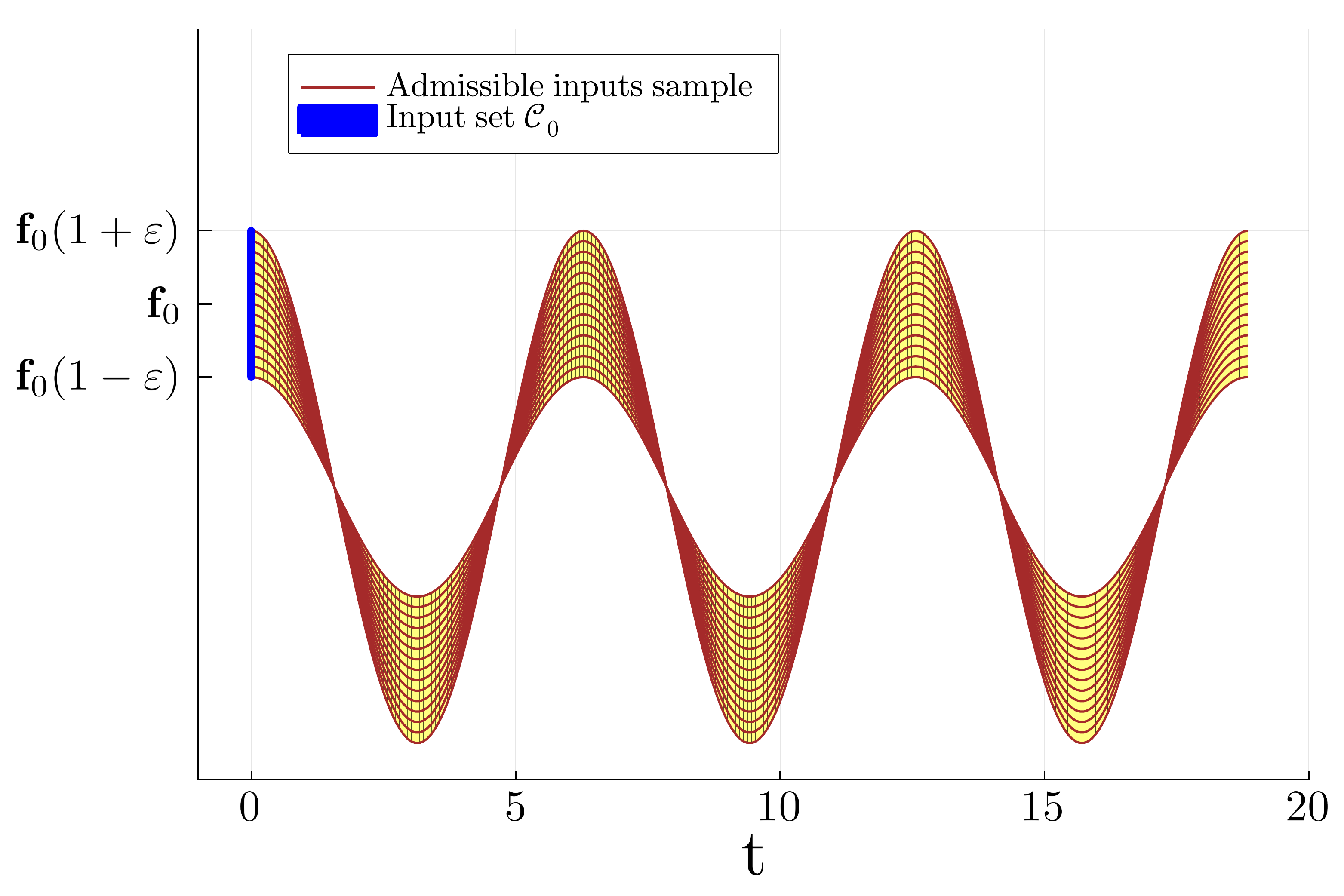}
		\caption{Example of time-varying input model.}
		\label{fig:sineforcesRA}
	\end{subfigure}
	\caption{Illustration of different input function assumptions. A set $\mathcal{C}_0$ is used to span all possible input functions with the specified behavior; a uniform sample is shown for each case.}
\end{figure}

{\rev To illustrate the formulation, assume that there is only one constant forcing term $\bff_0$ ($n_f = 1$, $\bfB^{(1)} = {\bf0}_{1\times 1}$), which corresponds to the case $\alpha = 0$ in Eq.~\eqref{eq:input_function_cases}. We take an interval $\mathcal{C}_0 = [1-\varepsilon, 1+\varepsilon] \subset \mathbb{R}$ for some $\varepsilon > 0$. With this initial condition, we can model all forcing functions that take any constant value within $\bff_0 \pm \varepsilon \bff_0$, as shown in Fig.~\ref{fig:constantforcesRA}.
On the other hand, in Fig.~\ref{fig:sineforcesRA} we show an example that corresponds to a time-varying (sinusoidal) input function with $\omega = 1$ in Eq.~\eqref{eq:input_function_cases}}.
The three-dimensional heat transfer problem in Example 4 of Section~\ref{sec:Example4_Heat3D} combines three different input terms.
}

\subsection{Conservative time discretization} \label{ssec:discretization}

In this section we explain how to perform a time discretization that is conservative in the sense of including all admissible behaviors in dense time, and continue with the description of different methods to propagate such initial set according to the system's state transition matrix.

We recall that $\mathcal{X}_0$ is the set of initial states and, since the analytic solution of Eq.~\eqref{eq:linearODE} at $t = \delta$ is $\varphi(\bfx_0, \delta) = \bfPhi \bfx_0$ for each $\bfx_0 \in \mathcal{X}_0$, then $\bfPhi \mathcal{X}_0$ is the set of transformed states after one time increment of size $\delta$, i.e. $\mathcal{R}^e(\mathcal{X}_0, \delta) = \bfPhi \mathcal{X}_0$.
The reach-set for the time interval $[0, \delta]$ is more difficult to obtain because there is no closed-form expression. 
Let $CH(X, Y)$ denote the \textit{convex hull} of the set union of two sets $X \subseteq \mathbb{R}^n$, $Y \subseteq \mathbb{R}^n$, i.e. the smallest convex set that contains $X$ and $Y$, defined by
\begin{equation}
	CH(X, Y) = \left\{ z \in \mathbb{R}^n : z = \lambda x + (1-\lambda) y,\quad \lambda \in [0, 1],~x \in X,~y \in Y \right\}.
\end{equation}
It is clear (compare e.g. with Fig.~\ref{fig:sdof_singleton_discretization}) that the analytic solution does not follow the straight line of the convex hull $CH(\mathcal{X}_0, \bfPhi \mathcal{X}_0)$ for intermediate times between $0$ and $\delta$, because trajectories curve due to the system's dynamics. Hence we need to enlarge, or bloat, the convex hull by a suitable amount such that all the exact reachable states $\Omega_0^e$ are covered.

To obtain a set $\Omega_0$ satisfying Eq.~\eqref{eq:discretization}, we propose a slight variation of the method presented in \citep{frehse2011spaceex}. The advantage of the new method is that it doesn't require to solve an optimization problem for each direction of interest.
The symbol $\boxdot(\mathcal{X})$ represents the symmetric interval hull of $\mathcal{X} \subset \mathbb{R}^n$, that is the smallest hyperrectangle centered at the origin that contains $\mathcal{X}$.
\begin{prop} \label{prop:step_intersect}
	Define the following pair of sets
	\begin{eqnarray}
		\Omega_0^+ &=& CH(\mathcal{X}_0, \bfPhi \mathcal{X}_0 \oplus E^+(\bfA, \mathcal{X}_0, \delta)) \label{eq:discretization_step_intersect_pos}\\
		\Omega_0^- &=& CH(\bfPhi \mathcal{X}_0, \mathcal{X}_0 \oplus E^+(\bfA, \bfPhi \mathcal{X}_0, \delta)) \label{eq:discretization_step_intersect_neg}
	\end{eqnarray}
	where $E^+(\bfA, \mathcal{X}_0, \delta) = \boxdot\left(\bfP(\vert \bfA \vert, \delta)\boxdot(\bfA^2 \mathcal{X}_0)\right)$
	and $\bfP(\bfA, \delta) = \sum_{i=0}^\infty \bfA^i \delta^{i+2}/(i+2)!$. Then, the set $\Omega_0 := \Omega^+_0 \cap \Omega_0^-$ encloses $\Omega_0^e$, i.e. $\Omega_0^e \subseteq \Omega_0$. Moreover, $\Omega_0 \to \Omega_0^e$ as $\delta \to 0$.
\end{prop}

The proof of Prop.~\ref{prop:step_intersect} is left to \ref{sec:prop1}.
The intuition behind Prop.~\ref{prop:step_intersect} is that first Eq.~\eqref{eq:discretization_step_intersect_pos} performs one step forward in time (from $t = 0$ to $t = \delta$) with a bloating set $E^+$ that is a small box around the origin ($E^+$ tends to zero as $\delta \to 0$), and second, Eq.~\eqref{eq:discretization_step_intersect_neg} performs one step in reverse time (from $t = \delta$ to $t = 0$) which are then intersected to remove spurious overapproximation errors.

The result in our running example is shown in  Fig.~\ref{fig:sdof_singleton_discretization}: note how the set $\Omega_0$ (orange) encloses the exact (curved) solution for all intermediate times between $0$ and $\delta$.
The case depicted in Fig.~\ref{fig:sdof_distributed_discretization} is more interesting: it corresponds to the case in which the set of initial states $\mathcal{X}_0$ is a hyperrectangle.
While the obtained set $\Omega_0$ (orange) may not be the tightest convex set containing the exact reachable set $\Omega_0^e$, it is nonetheless a very good approximation.
Finally, let us remark that $\Omega_0$ is in general a polygon (or a polytope in higher dimensions). In practice, either we approximate it with a hyperrectangle (see Fig.~\ref{fig:sdof_distributed_discretization} in lightblue), or use lazy (i.e.~on-demand) evaluations of the support function. While the former approach leads to faster runtimes for large time horizons, the latter does not introduce additional overapproximation error. The hyperrectangular approximation of $\Omega_0$ can be easily computed using support functions.

Having obtained the set $\Omega_0$ to fulfill the conservative time discretization requirement in Eq.~\eqref{eq:discretization}, we proceed in the next paragraphs to explain the set-based propagation of the sequence. To solve Eq.~\eqref{eq:discrete_recurrence}, three different approaches are presented: zonotopes, hyperrectangles and support functions. Each method offers a different efficiency vs accuracy tradeoff and scales differently with the system's dimension. For Section \ref{sec:results} we have selected the method that is more suitable for each problem.
A complete evaluation of the different methods is out of the scope of this article (see also \citep{althoff2020arch} for a comparison of results obtained by different software tools on a set of benchmark problems).

\subsection{Set propagation} \label{ssec:set_propagation}

\subsubsection{Set propagation using zonotopes}
\label{ssec:set_propagation_zonotopes}

To solve the set-based recurrence of Eq.~\eqref{eq:discrete_recurrence} using zonotopes, it suffices to observe that zonotopes are closed under linear maps: $\bfM Z = \langle \bfM \bfc, \bfM \bfG \rangle_Z$ is again a zonotope 
for any $\bfM \in \mathbb{R}^{m \times n}$, obtained by transforming the center and generators of $Z$ by the linear map $\bfM$.
Therefore if $Z_0 = \langle \bfc_0, \bfG_0 \rangle_Z$ is a zonotope, it holds that
\begin{equation}\label{eq:recurrence_sol_zonotope}
	Z_k = \langle \bfPhi^k \bfc_0, \bfPhi^k \bfG_0 \rangle_Z, \qquad k = 0, 1, \ldots, N-1.
\end{equation}
In general, the initial set $\Omega_0$ from Prop.~\ref{prop:step_intersect} is not a zonotope so we overapproximate it with a hyperrectangle as a pre-processing step.
Set propagation using zonotopes is illustrated in Fig.~\ref{fig:sdof_singleton_flowpipe} for singleton resp. in Fig.~\ref{fig:sdof_distributed_flowpipe} for distributed initial conditions (light-green).

Finally we remark that, in a simplified model of computation, the time complexity of Eq.~\eqref{eq:recurrence_sol_zonotope} is $\mathcal{O}(Nn^2 p)$, and the space complexity is $\mathcal{O}(Nnp)$.
Typically, $p = n$ hence storing all intermediate zonotopes can be expensive if the state-space dimension $n$ is sufficiently large. However, there are situations in which we are only interested in obtaining minimum and maximum bounds of the flowpipe and in that case it is more efficient and equally accurate to propagate axis-aligned hyperrectangles. This is the subject of the next subsection.

\subsubsection{Set propagation using hyperrectangles} \label{ssec:set_propagation_hyperrectangles}

Hyperrectangular approximations are particularly useful in problems where it is relevant to find the maximum or minimum values of all the state variables; for example in the longitudinal oscillations of a bar if we are interested in computing the maximum displacements and velocities (cf. Example~\ref{sec:Example2_Clamped}), or in a massive concrete structure model we may be interested in the maximum temperature achieved at any mesh point (cf. Example~\ref{sec:Example4_Heat3D}).

Given a zonotope $Z$, we say that $H$ is a \textit{tight} hyperrectangular approximation if $H$ is the smallest hyperrectangle that contains $Z$. The following proposition formalizes the idea of approximating the zonotope flowpipe with tight hyperrectangles.

\begin{prop} \label{prop:box_propagation}
	Let $H_0 = \langle \bfc_0, \bfr_0 \rangle_H$ be a hyperrectangle enclosing $\Omega_0^e$. For each $k \geq 1$, define the sequence
	$\bfc_{k} = \bfPhi^k \bfc_0$
	and $\bfr_{k} = \vert \bfPhi^k \vert \bfr_0$,  where the absolute values are taken component-wise. Then, the sequence of hyperrectangles $H_k = \langle \bfc_k, \bfr_k\rangle_H$ tightly overapproximates the solutions of the recurrence \eqref{eq:discrete_recurrence} for all $k = 0,\ldots, N-1$. 
\end{prop}

The proof of Prop.~\ref{prop:box_propagation} is left to \ref{sec:prop2}.
Concerning the computational cost, the time complexity of the method in Prop.~\ref{prop:box_propagation} is  $\mathcal{O}(Nn^2)$, and the space complexity is $\mathcal{O}(Nn)$. Note that computing with hyperrectangles is comparatively less expensive than zonotopes, although the latter have more expressive power than the former.

Set propagation using hyperrectangles is illustrated in Fig.~\ref{fig:sdof_singleton_flowpipe} for singleton initial conditions resp. in Fig.~\ref{fig:sdof_distributed_flowpipe} for distributed initial conditions (light-blue). For the 6th reach-set the hyperrectangle and the zonotope coincide; in general though, the hyperrectangle encloses the zonotope. However, for the $k$-th reach-set, the maximum and minimum of $Z_k$ are always the same to those obtained with $H_k$ as claimed in Prop.~\ref{prop:box_propagation}.

While solving Eq.~\eqref{eq:discrete_recurrence} using hyperrectangles scales very well with the dimension, this approach is limited to obtaining maximum (resp. minimum) values over solutions. We can interpret such computation as the maximum of the scalar product $\bfd^T \bfx(t)$ over all $\bfx(t) \in \mcR^e(\mcX_0, [\delta (k-1), \delta k])$, where $\bfd$ is any canonical direction of $\mathbb{R}^n$, i.e. of the form $\bfd = (0, \ldots, 1, \ldots, 0)^T$. The minimum can be obtained analogously. Using support functions, as explained in the next section, we can generalize the computation for an arbitrary direction $\bfd \in \mathbb{R}^n$.

\subsubsection{Set propagation using support functions} \label{ssec:set_propagation_support_functions}

The first systematic attempt to solve Eq.~\eqref{eq:discrete_recurrence} using support functions was presented in \citep{LeGuernic2010250}. The idea is to make use of the property $\rho(\bfd, \bfM X) = \rho(\bfM^T \bfd, X)$ for any matrix $\bfM \in \mathbb{R}^{n \times n}$. 
By successive application of this rule,
\begin{equation} \label{eq:recurrence_support_function}
	\rho(\bfd, X_k) = \rho((\bfPhi^T)^k \bfd, X_0),\qquad k = 0, 1,\ldots, N-1.
\end{equation}
The support function of common set representations (e.g. hyperrectangles and zonotopes) is known analytically so in those cases it can be computed very efficiently. In addition, the support function of an intersection satisfies that $\rho(\bfd, X\cap Y) \leq \min \{\rho(\bfd, X), \rho(\bfd, Y)\}$, while for the case of Minkowski sums, $\rho(\bfd, X \oplus Y) = \rho(\bfd, X) + \rho(\bfd, Y)$.

Using support functions we can compute linear combinations of the state variables, i.e. we can solve the recurrence \eqref{eq:discrete_recurrence} \textit{along a given direction} directly solving a recurrence for the support function, instead of computing the flowpipe for all variables and then extracting the relevant information.
A practical application of the solution method in Eq.~\eqref{eq:recurrence_support_function} to obtain the spatial gradient in a heat transfer problem is presented in Example \ref{sec:Example3_Heat}.

{\rev
\subsection{Example code}
}

The pseudo-code presented in Algorithm~\ref{algo:reach} illustrates the main steps of our approach.
To describe with more detail the implementation, a simplified version of the proposed method is presented for the harmonic oscillator from Section~\ref{sec:illustrative_example}. For numeric and symbolic operations with sets, we use the LazySets.jl\footnote{The library can be downloaded from \href{https://github.com/JuliaReach/LazySets.jl}{github.com/JuliaReach/LazySets.jl}} library, written in the Julia programming language \citep{bezanson2017julia}. The library implements all required set representations. We show the Julia commands (introduced after the \texttt{julia>} prompt) followed by the resulting output.

\begin{algorithm}[htbp!]
	\caption{Set propagation method pseudo-code.}
	\label{algo:reach}
	\KwIn{%
		$\bfM$, $\bfC$, $\bfK$: FEM assembled matrices,
		$\bfF$: vector of loads,
		$\mcX_0 \subseteq R^n$: initial states set
		\goodbreak
		$\delta$: time step increment, $N$: number of time steps\\
	}
	\vspace*{1mm}
	\KwOut{List of reachable sets}
	\vspace*{2mm}
	$\bfA, \mcC_0 = homogeneize(\bfM, \bfC, \bfK, \bfF)$ \tcp*{use Eq.\eqref{eq:homogeneized_matrix_heat} or \eqref{eq:homogeneized_matrix_dynamics}} \label{line:homogeneize}
	$\Omega_0 = discretize(\mcX_0 \times \mcC_0, \bfA, \delta)$ \tcp*{time discretization using Prop.\ref{prop:step_intersect}}\label{line:discretize}
	$\bfPhi = exp(\bfA \delta)$ \tcp*{matrix exponential} \label{line:expm}
	$X_0 = \Omega_0$ \tcp*{reachable states approximation for $[0, \delta]$} \label{line:setprop_init}
	\For{$k = 0$ \KwTo $N-2$}{%
		$ X_{k+1} = \bfPhi X_k$ \tcp*{set propagation; use either Eq.\eqref{eq:recurrence_sol_zonotope}, Prop.\ref{prop:box_propagation} or Eq.\eqref{eq:recurrence_support_function}}
		\label{line:setprop}
	} \label{line:setprop_end}
	\Return{$X_0,  X_1, \cdots, X_{N-1}$} \tcp*{flowpipe approximation}
\end{algorithm}

The first step of the method is the formulation of the initial-value problem as an homogeneous, first order system which corresponds to Line \ref{line:homogeneize} of the pseudo-code. In the present case, there are no external forces. The initial states are defined choosing a suitable set representation, e.g. a hyperrectangle in this case.

\begin{minipage}{\linewidth}
\begin{lstlisting}[label=lst:formulation, caption=Stage 1. Initial-value problem formulation.]
julia> using LazySets
julia> A = [0 1; -(4π)^2 0] # state matrix
2×2 Matrix{Float64}:
   0.0    1.0
-157.914  0.0
julia> X₀ = Hyperrectangle([1.0, 0.0], [0.1, 0.1]); # initial states (center, radii)
julia> δ = 0.025; # time step
\end{lstlisting}
\end{minipage}

The second step of the algorithm is the conservative time discretization, which corresponds to Line \ref{line:discretize} of the pseudo-code. Here we computed $\bfP$ and $\bfPhi$ in full; for larger problems, matrix functions of exponential type are used.
The operation of multiplication between matrices and sets, as well as the symmetric interval hull are already implemented in LazySets.jl. Symbolic Minkowski sums, convex hulls and intersections are also implemented in the library. Such specific features simplify writing set-based algorithms.
 
\begin{minipage}{\linewidth}
\begin{lstlisting}[label=lst:discretization, caption=Conservative time discretization computation.]
julia> P = sum(abs.(A)^i * δ^(i+2) / factorial(i+2) for i in 0:10) # cf. Prop. 1
2×2 Matrix{Float64}:
  0.00031508  0.00000262
  0.00041327  0.00031508
julia> E₊(X₀) = symmetric_interval_hull(P * symmetric_interval_hull(A^2 * X0)); E₊(X₀)
Hyperrectangle{Float64, Vector{Float64}, Vector{Float64}}([0.0 0.0], [0.05477208, 0.07676220])
julia> Φ = exp(A*δ)
2×2 Matrix{Float64}:
  0.95105652  0.02459079
 -3.88322208  0.95105652
julia> Ω₀ = CH(X₀, Φ*X₀ ⊕ E₊(X₀)) ∩ CH(Φ*X₀, X₀ ⊕ E₊(Φ*X₀)); # symbolic discretized states
julia> H₀ = box_approximation(Ω₀) # box enclosure
Hyperrectangle{Float64, Vector{Float64}, Vector{Float64}}([0.97471, -2.13332], [0.12868, 2.23332])
\end{lstlisting}
\end{minipage}

Here we have computed an overapproximation of $\Omega_0$, namely a hyperrectangular enclosure.\footnote{The set $\Omega_0$ may also be propagated symbolically, i.e. without overapproximation, using support functions.} Such set is centered in $[0.97, -2.13]$ and has radius $0.13$ and $2.23$ on the $u$ and $v$ coordinates respectively, as shown in light-blue in Fig.~\ref{fig:sdof_distributed_discretization}. 

The third step of the method consists in set propagation (see Line \ref{line:setprop} of the pseudo-code).
To illustrate the propagation using zonotopes, let us obtain the center and generators matrix of the sixth zonotope in Fig.~\ref{fig:sdof_distributed} (hence $k = 5$):
 
\begin{minipage}{\linewidth}
\begin{lstlisting}[label=lst:setprop_zono, caption=Set propagation using zonotopes.]
julia> Z6 = Zonotope(Φ^5 * H₀.center, Φ^5 * genmat(H₀));
Zonotope{Float64, Vector{Float64}, Matrix{Float64}}(
[-0.16976461, -12.24853154],  # center
[-0.00000000  0.17772235;     # generators matrix
 -1.61711795  0.00000000]) 
\end{lstlisting}
\end{minipage}
Observe that each generator is proportional to a canonical vector, which explains the axis-aligned shape of the sixth zonotope in Fig.~\ref{fig:sdof_distributed}. Similarly we can obtain the corresponding hyperrectangle using Prop.~\ref{prop:box_propagation}:

\begin{minipage}{\linewidth}
	\begin{lstlisting}[label=lst:setprop_hyper, caption=Set propagation using hyperrectangles.]
julia> H6 = Hyperrectangle(Φ^5 * H₀.center, abs.(Φ^5) * H₀.radius)
Hyperrectangle{Float64, Vector{Float64}, Vector{Float64}}(
[-0.16976, -12.24853],    # center
[ 0.17772,   1.61712])    # radius
	\end{lstlisting}
\end{minipage}
Finally we note that set propagation using support functions is available in the libary by use of the \texttt{support\_function} method (Unicode alias: $\rho$).


%
%

\section{Numerical results} \label{sec:results}

In order to obtain comparative results about the performance and accuracy of the methods discussed, five numerical examples are solved.
Each example introduces a new motivation and level of complexity in the comparative process of the methods.

The reachability problems are solved using the JuliaReach library \citep{juliareach, bogomolov2019juliareach}, implemented in the Julia language \citep{bezanson2017julia}.
The computations associated to the Finite Element Method are done using source code from the ONSAS library \citep{onsas} implemented in GNU-Octave \citep{eaton1997gnu}.
The visualization is done using JuliaPlots \citep{tom_breloff_2021_4725318} and Paraview \citep{ahrens2005paraview}, {\rev and unstructured meshes are created using GMSH \citep{Geuzaine2009a}}.
For matrix functions of exponential type, the implementation in \citep{DifferentialEquations.jl-2017} is used.
All the results are executed in a laptop computer running Linux OS, with an Intel(R) Core(TM) i7-8705G CPU, 3.10GHz processor and 16GB RAM using Julia v.1.6.1. {\rev The results presented can be obtained using scripts publicly-available\footnote{Available at: \href{https://github.com/JuliaReach/SetPropagation-FEM-Examples}{github.com/JuliaReach/SetPropagation-FEM-Examples}}}.

\subsection{Example 1: accuracy analysis of the harmonic oscillator}

In this example, the one-dimensional harmonic oscillator presented in Section~\ref{sec:illustrative_example} is considered. We compare the set propagation results in terms of accuracy with respect to reference numerical integration methods and the analytic solution.

\subsubsection{Problem definition}

The governing equation is given in Eq.~\eqref{eq:harmonic_oscillator} and the parameters considered are frequency $\omega = 4\pi$ (i.e.~period $0.5$ s), and initial conditions $u(0)=1$, $\dot{u}(0)=v(0) = 0$.
Different time steps are considered: $\delta = 0.025$ s and $\delta = 0.05$ s.

\subsubsection{Numerical resolution}

\begin{figure}[htb]
	\centering
	\begin{subfigure}[b]{0.45\textwidth}
		\centering
		\includegraphics[width=\textwidth]{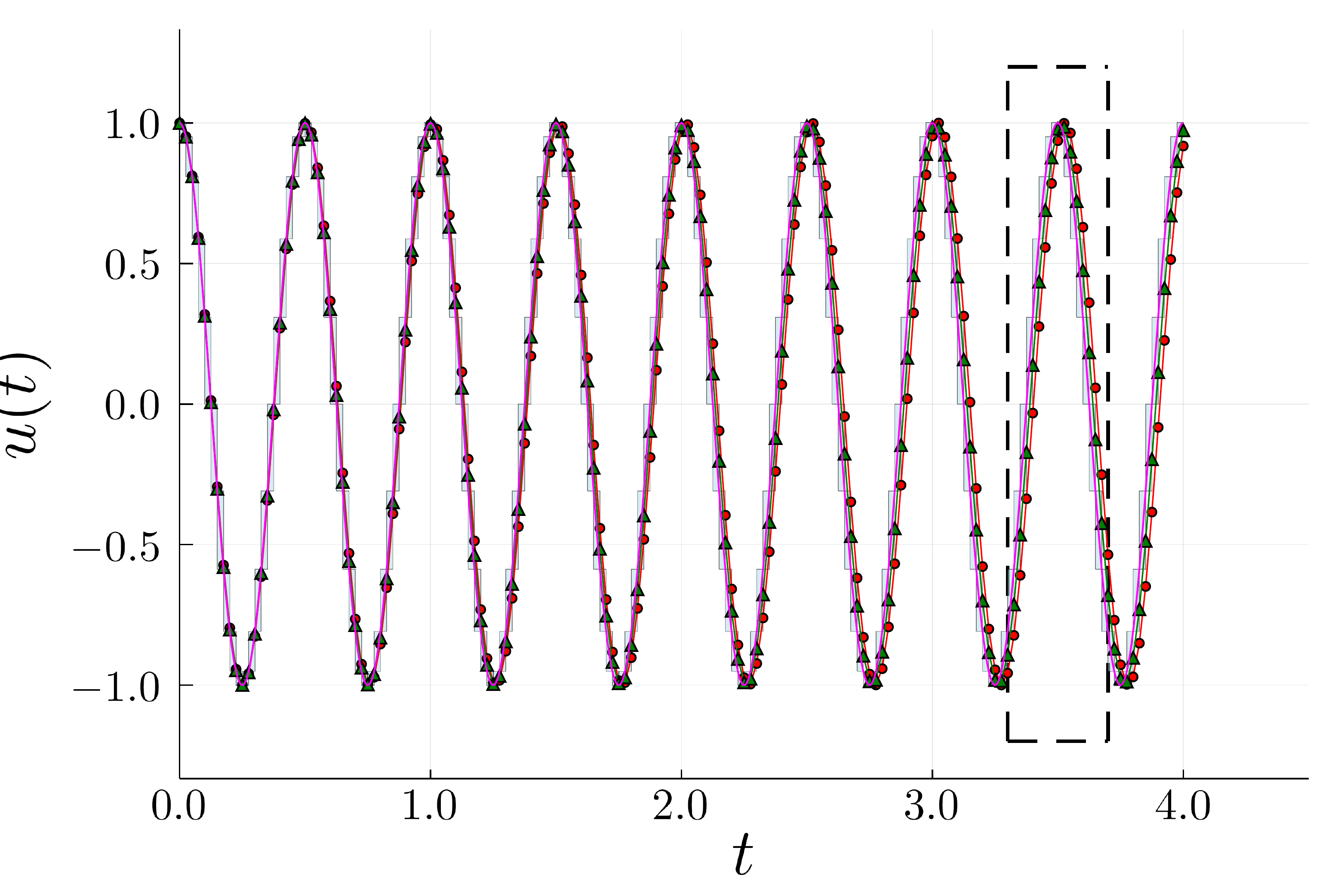}
		\caption{Displacements for $t \in [0, 4]$.}
		\label{fig:SDOF_x_vs_t_alpha_small_all}
	\end{subfigure}
	~~
	\begin{subfigure}[b]{0.45\textwidth}
		\centering
		\includegraphics[width=\textwidth]{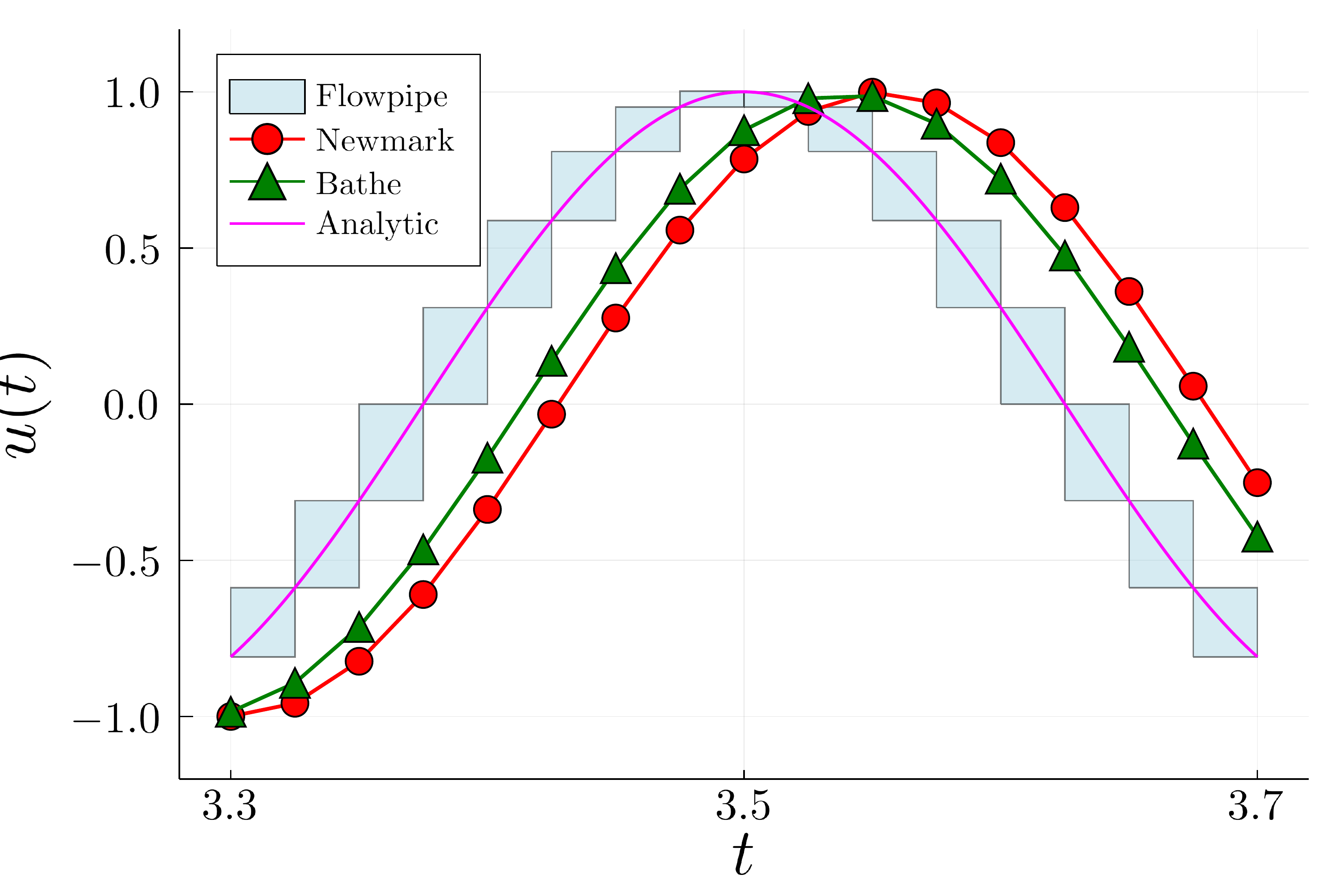}
		\caption{Zoom for $t \in [3.3, 3.7]$.}
		\label{fig:SDOF_x_vs_t_alpha_small_zoom}
	\end{subfigure}
	\caption{Example 1: displacements obtained using $\alpha = \delta / T = 0.05$. The analytic solution (magenta) is shown superposed with the numerical solution obtained with Newmark's method (red circles) and Bathe's method (green triangles).}
	\label{fig:SDOF_x_vs_t_alpha_small}
\end{figure}

The results obtained using time step $\delta = 0.025$ s are shown in Figure~\ref{fig:SDOF_x_vs_t_alpha_small}. The solution obtained using the Newmark method is shown with red circles, the solution using the Bathe method is represented by green triangles and the analytic solution is represented by a solid magenta line.
Regarding the solution obtained by set propagation, we represent the flowpipe using polygons and show their projection onto the $u(t)$ axis using light-blue boxes.

In Figure~\ref{fig:SDOF_x_vs_t_alpha_small_zoom} it can be observed that the analytic solution is contained in the flowpipe. On the other hand, the numerical methods show a considerable period elongation (after 7 periods of simulation). Execution times are reported in Table~\ref{tab:runtimes_sdof}.

\begin{table}[htb]
	\centering
	{\rev
	\begin{tabular}{@{}clll@{}}
		\cmidrule(l){2-4}
		& Set Propagation & Newmark & Bathe \\ \midrule
		$\delta = 0.025$ & 0.29 ms (140 kB)  & 0.44 ms (321 kB)  & 0.56 ms (561 kB) \\ \midrule
		$\delta = 0.05$ & 0.17 ms (91 kB)  & 0.18 ms (161 kB)  & 0.36 ms (282 kB)
	\end{tabular}
	\caption{Execution times and memory usage for the harmonic oscillator example.}
\label{tab:runtimes_sdof}
}
\end{table}

\begin{figure}[htbp]
	\centering
	\begin{subfigure}[b]{0.45\textwidth}
		\centering
		\includegraphics[width=\textwidth]{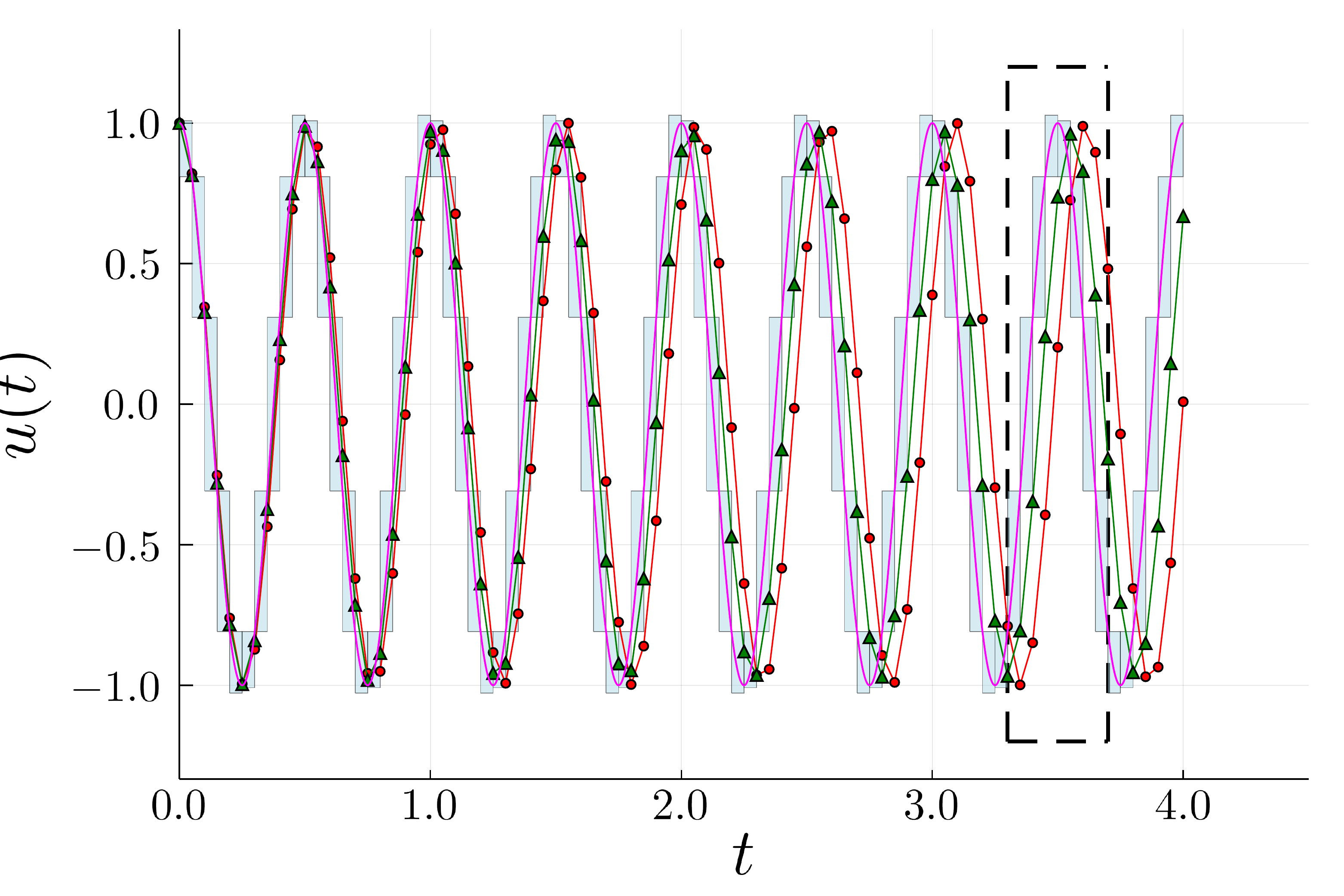} 
		\caption{Displacements for $t \in [0, 4]$.}
		\label{fig:SDOF_x_vs_t_alpha_big_all}
	\end{subfigure}
	~~
	\begin{subfigure}[b]{0.45\textwidth}
		\centering
		\includegraphics[width=\textwidth]{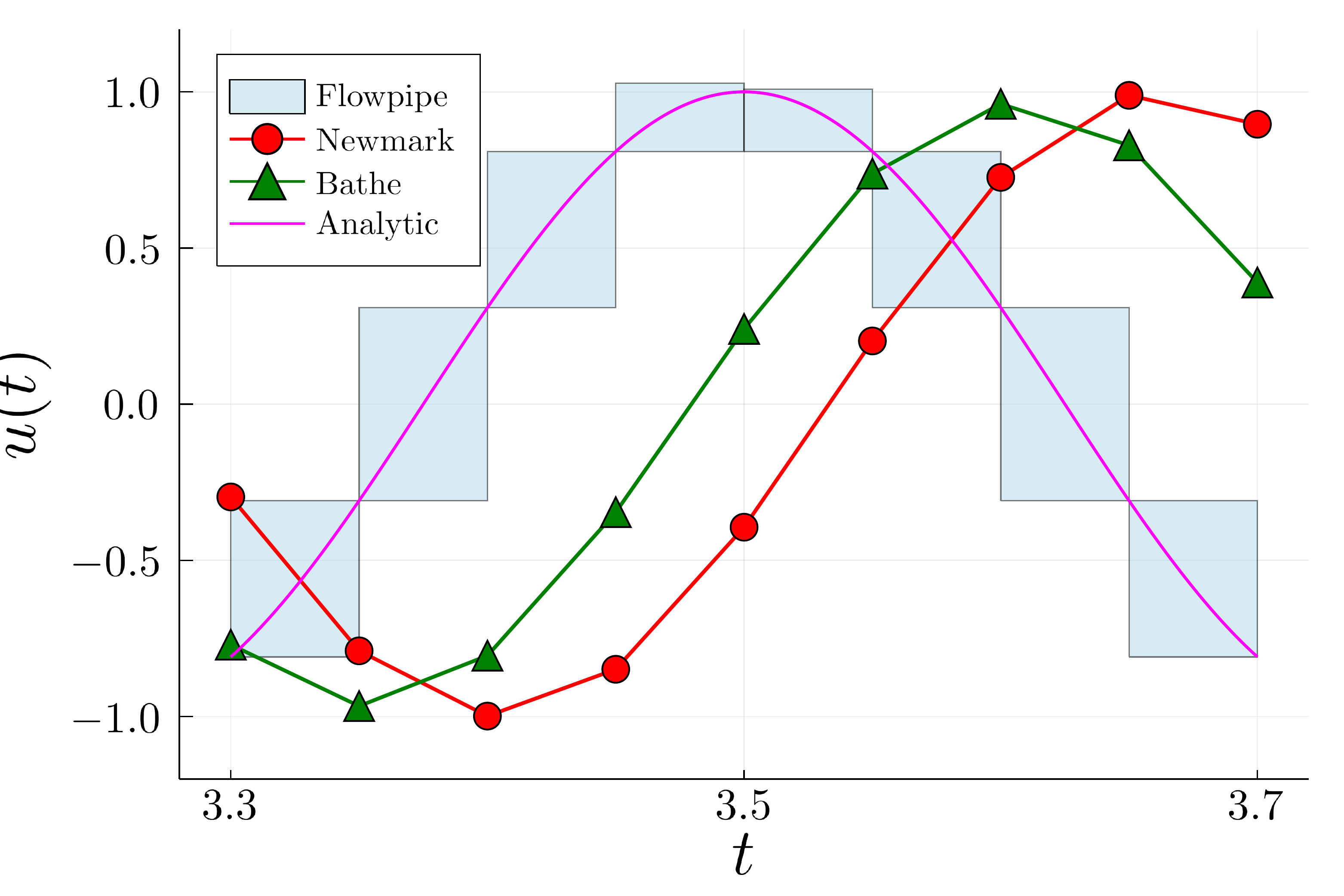}
		\caption{Zoom for $t \in [3.3, 3.7]$.}
		\label{fig:SDOF_x_vs_t_alpha_big_zoom}
	\end{subfigure}
	\caption{Example 1: displacements using $\alpha = \delta / T = 0.1$. The analytic solution (magenta) is shown superposed with the numerical solution obtained with Newmark's method (red circles) and Bathe's method (green triangles).}
	\label{fig:SDOF_x_vs_t_alpha_big}
\end{figure}

The results obtained using a larger time step $\delta = 0.05$ s are shown in Figure~\ref{fig:SDOF_x_vs_t_alpha_big}. %
In this case, in Figure~\ref{fig:SDOF_x_vs_t_alpha_big_zoom} it can be observed a greater period elongation for the numerical methods. The flowpipe still includes the analytic solution but the width of the intervals on the $u(t)$ axis is larger than in the previous case.

\subsubsection{Accuracy analysis}

Accuracy analysis of numerical integration methods can be done using two tools: Period Elongation ($PE$) and Amplitude Decay ($AD$) \citep{Bathe2012, Bathe2014}. %
In this section we compute these magnitudes for the numerical results obtained and draw a comparison with respect to the set propagation method.
The numerical data used to compute these magnitudes corresponds to a simulation for 50 periods, i.e. 25 seconds.

The PE can be computed using the following expression:
\begin{equation}
	\centering
	PE=\dfrac{T_{num}-T_{nat}}{T_{nat}},
	\label{eq:PE}
\end{equation}
where $T_{nat}$ is the natural period of the oscillator and $T_{num}$ is the period of the numerical solution. %
For the numerical integration methods, we compute $T_{num}$ using standard spectral analysis methods. On the other hand, to our knowledge there are no methods for the computation of PE estimates of flowpipes, thus we devised a procedure to compute $T_{num}$ in the set-based setting. The notions used in the method are presented in \ref{sec:apePEAD}, were both PE and AD estimations are described.

\begin{figure}[hbtp]
	\centering
	\begin{subfigure}[b]{0.45\textwidth}
		\centering
		\includegraphics[width=\textwidth]{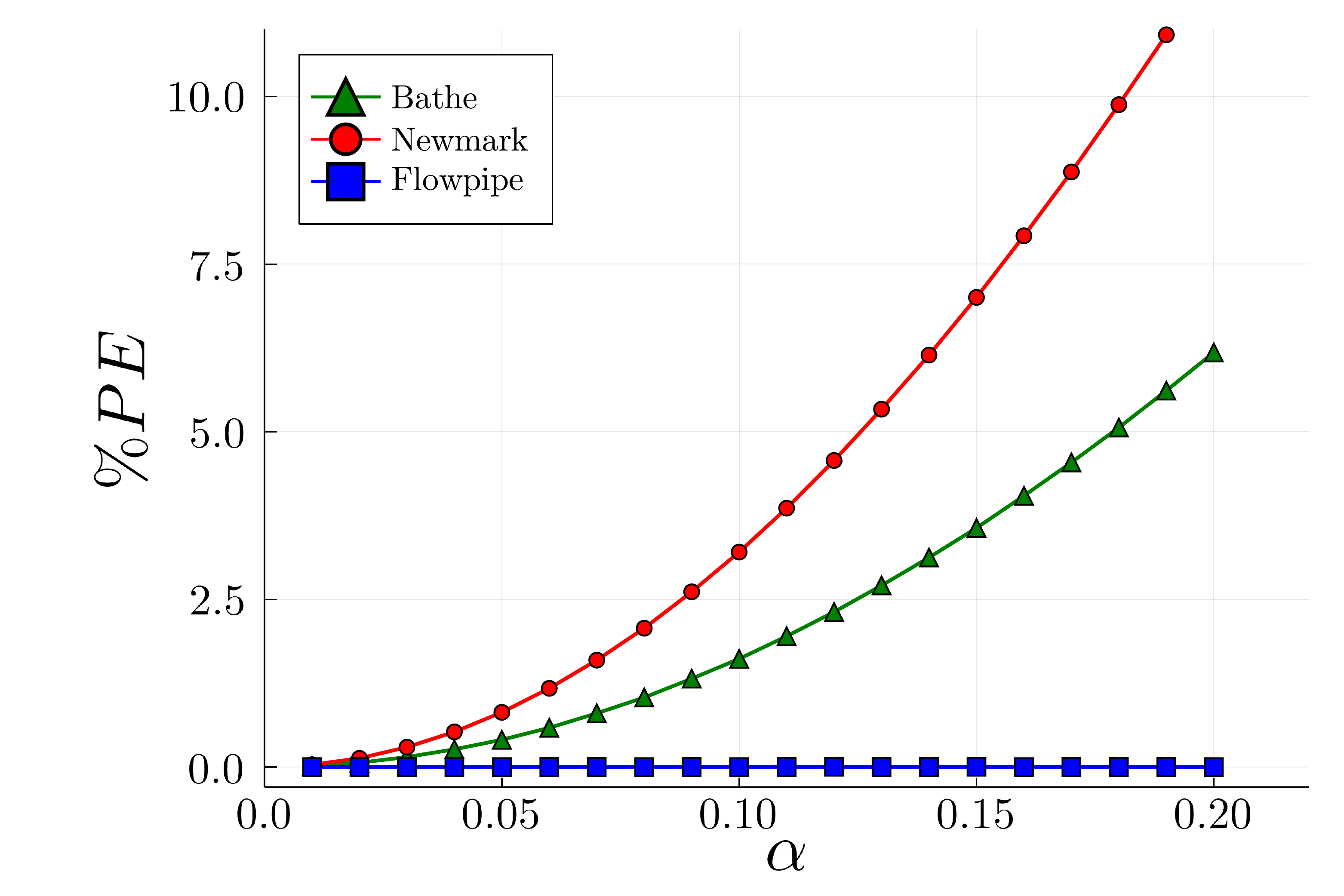}
		\caption{Percentage period elongation ($PE$).}
		\label{fig:PE_all}
	\end{subfigure}
	~~
	\begin{subfigure}[b]{0.45\textwidth}
		\centering
		\includegraphics[width=\textwidth]{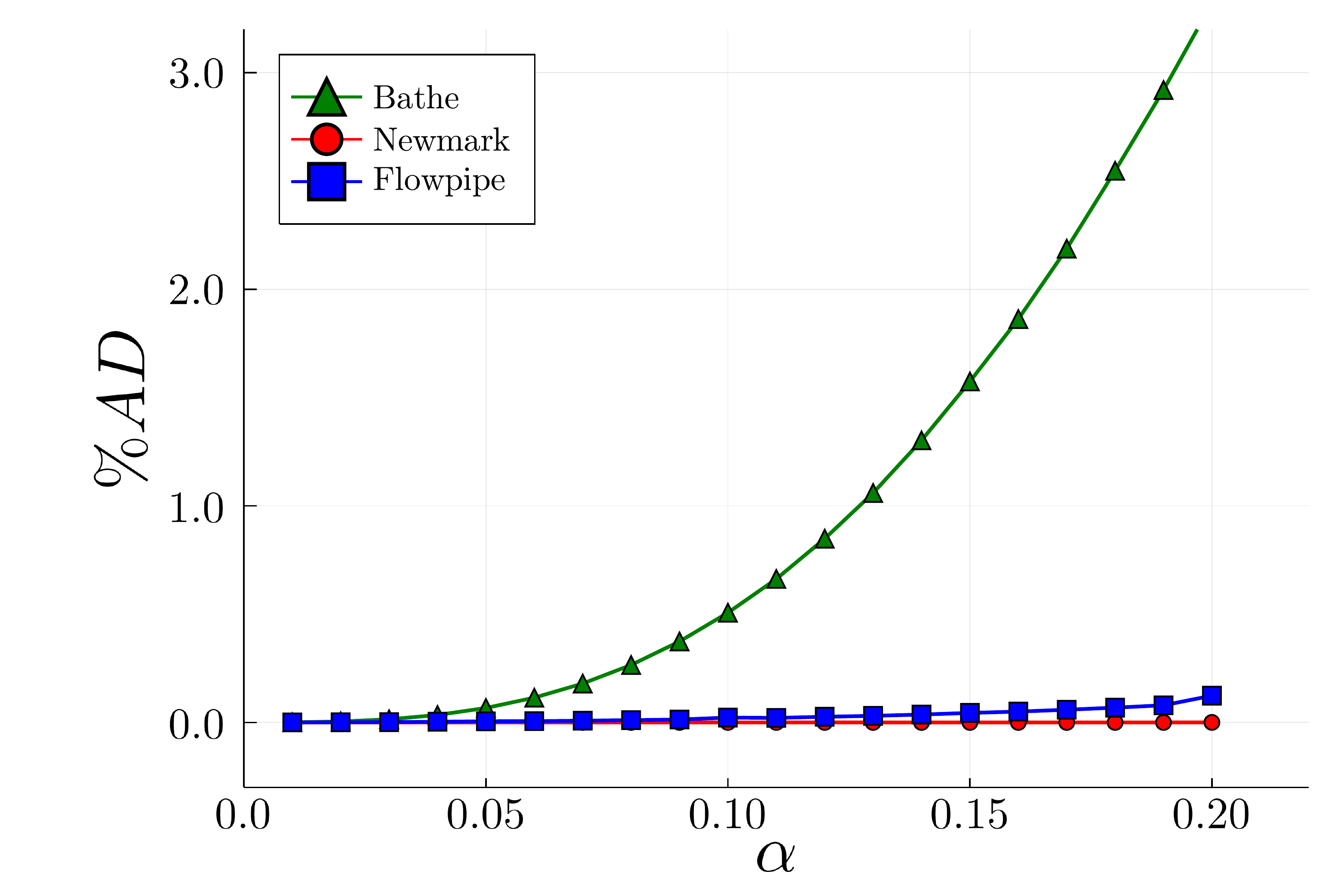}
		\caption{Percentage amplitude decay ($AD$).}
		\label{fig:AD_all}
	\end{subfigure}
	\caption{Example 1: PE and AD as a function of the relative step size $\alpha = \delta / T$. In red (green) color, the results for Newmark (Bathe) methods. In blue, the results obtained with the set propagation approach.}
	\label{fig:sdof_PE_AD}
\end{figure}

In Fig.~\ref{fig:PE_all} the PEs computed for the Bathe, Newmark and set propagation methods are presented, where the latter is shown in blue squares and the Bathe and Newmark methods using triangles and circles, respectively. %
The values of $\alpha$ considered are intended to reproduce the results shown in \citep{Bathe2014}. %
It can be seen that Newmark and Bathe methods present PE values, that match with literature results, while the set propagation results do not present PE.

The numerical amplitude $A_{num}$ after $n_A$ periods is $A_{num} = A (1 - AD)^{n_A}$, hence the Amplitude Decay $AD$ is:
\begin{equation}
	\centering
	AD = 1 - \left(\dfrac{A_{num}}{A}\right)^{\frac{1}{n_A}}.
	\label{eq:AD2}
\end{equation}

In Fig.~\ref{fig:AD_all} the computed values of AD are shown. As it is expected, the Newmark method does not present AD while the Bathe method presents a curve that matches the literature results. The set propagation method has not significant AD, with values below $0.15 \%$ for any $\alpha$.

The results obtained show that the set propagation method presents advantages with respect to the numerical integration methods, in terms of PE, and competitive results in terms of AD.
The set propagation method presents good results for single degree of freedom problems. In the following example, the problem is extended to increase the number of unknowns.

\clearpage
\subsection{Example 2 - Clamped-Free Bar} \label{sec:Example2_Clamped}

We consider a uni-dimensional wave propagation problem from \citep{Malakiyeh2019}. 
The purpose of this example is to evaluate the conservative time integration scheme for single initial conditions in a high-dimensional problem.
Subsequent examples consider variations in the initial conditions and input parameters.

\subsubsection{Problem definition}

\begin{figure}[htb]
	\centering
	\def\svgwidth{0.55\textwidth}
\begingroup%
  \makeatletter%
  \providecommand\color[2][]{%
    \errmessage{(Inkscape) Color is used for the text in Inkscape, but the package 'color.sty' is not loaded}%
    \renewcommand\color[2][]{}%
  }%
  \providecommand\transparent[1]{%
    \errmessage{(Inkscape) Transparency is used (non-zero) for the text in Inkscape, but the package 'transparent.sty' is not loaded}%
    \renewcommand\transparent[1]{}%
  }%
  \providecommand\rotatebox[2]{#2}%
  \newcommand*\fsize{\dimexpr\f@size pt\relax}%
  \newcommand*\lineheight[1]{\fontsize{\fsize}{#1\fsize}\selectfont}%
  \ifx\svgwidth\undefined%
    \setlength{\unitlength}{395.30601229bp}%
    \ifx\svgscale\undefined%
      \relax%
    \else%
      \setlength{\unitlength}{\unitlength * \real{\svgscale}}%
    \fi%
  \else%
    \setlength{\unitlength}{\svgwidth}%
  \fi%
  \global\let\svgwidth\undefined%
  \global\let\svgscale\undefined%
  \makeatother%
  \begin{picture}(1,0.24316497)%
    \lineheight{1}%
    \setlength\tabcolsep{0pt}%
    \put(0,0){\includegraphics[width=\unitlength,page=1]{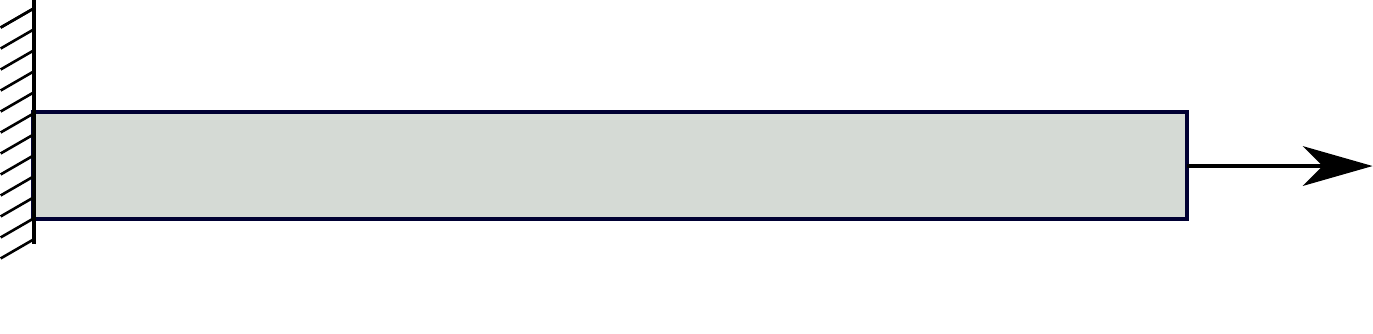}}%
    \put(0.90201673,0.17736682){\color[rgb]{0,0,0}\makebox(0,0)[lt]{\lineheight{1.25}\smash{\begin{tabular}[t]{l}$F(t)$\end{tabular}}}}%
    \put(0,0){\includegraphics[width=\unitlength,page=2]{clamped.pdf}}%
    \put(0.05027397,0.2127857){\color[rgb]{0,0,0}\makebox(0,0)[lt]{\lineheight{1.25}\smash{\begin{tabular}[t]{l}$x, u$\end{tabular}}}}%
    \put(0,0){\includegraphics[width=\unitlength,page=3]{clamped.pdf}}%
    \put(0.38079087,0.03924555){\color[rgb]{0,0,0}\makebox(0,0)[lt]{\lineheight{1.25}\smash{\begin{tabular}[t]{l}$L$\end{tabular}}}}%
  \end{picture}%
\endgroup%

	\caption{Example 2: diagram of the clamped-free bar excited by end load.}
	\label{fig:clampedDiagram}
\end{figure}

The domain consists of a bar of length $L=200$ and cross-section area $A=1$, formed by a linear elastic material with Young modulus $E = 30\times 10^6$ and density $\rho = 7.3 \times 10^{-4}$.
The free end of the bar is submitted to a step force $F(t) = 10,000 H(t)$, where $H(t)$ is the Heaviside function. %
The bar is governed by the partial differential equation
\begin{equation}
	EA \dfrac{\partial^2 u} {\partial x^2}(x,t) - \rho A \frac{\partial^2 u}{\partial t^2}(x,t) = 0, \label{eq:clamped_bar}
\end{equation}
where $u(x,t)$ is the displacement of the point in position $x$ at time $t$, considering the axis shown in Fig.~\ref{fig:clampedDiagram}.

The bar is considered to be initially at rest, with $u(x, 0) = 0$ and $\dot{u}(x, 0) = 0$ for all $x \in [0, L]$. %
The boundary conditions are $u(0, t) = 0$, corresponding to the fixed end, and  $\sigma(L, t) A = F(t)$, for the free end, where $\sigma(x,t) = E \dfrac{\partial u}{\partial x}(x,t)$ is the stress function at point $x$ and time $t$.

The analytical solution of this problem in the continuum domain can be obtained using mode superposition \citep{geradin2014mechanical}, and it is given by:
\begin{equation}
	u(x, t) = \dfrac{8FL}{\pi^2 E A } \sum_{s = 1}^{\infty} \left\{ \dfrac{(-1)^{s-1}}{(2s-1)^2}\sin \dfrac{(2s-1)\pi x}{2L}\left(1 - \cos \dfrac{(2s-1)\pi \mu t}{2L} \right)
	\right\},
	\label{eq:clamped_solution}
\end{equation}
where $\mu = \sqrt{E/\rho}$.

\subsubsection{Numerical resolution}

The numerical space discretization is done considering $N=1000$ two-node finite elements with linear interpolation. The following $N\times N$ mass and stiffness matrices are obtained:
\begin{equation*}
	\bfK = \dfrac{EA}{\ell}
	\begin{bmatrix}
		2 & -1 &  &  &  & 0  \\
		-1  & 2 & -1 &   &  &  \\
		& -1 & 2 & \ddots &  &  \\
		&  & \ddots & \ddots & -1 &  \\
		&  &  & -1 & 2 & -1 \\
		0 &  &  &  & -1 & 1 \\
	\end{bmatrix},~~ 
	\bfM = \dfrac{\rho A \ell}{2}
	\begin{bmatrix}
		2 &  &  &  &  & 0  \\
		& 2 &  &   &  &  \\
		&  & 2 &  &  &  \\
		&  &  & \ddots &  &  \\
		&  &  & & 2 & \\
		0 &  &  &  &  & 1 \\
	\end{bmatrix},
	\label{eq:clampedMatrices}
\end{equation*}
where $\ell = L / N$ is the length of each element. %
The damping matrix $\bfC$ in Eq.~\eqref{eqn:dynamicsfem} is a zero matrix of dimension $N\times N$. %

\begin{figure}[htb]
	\centering
	\includegraphics[width=0.75\textwidth]{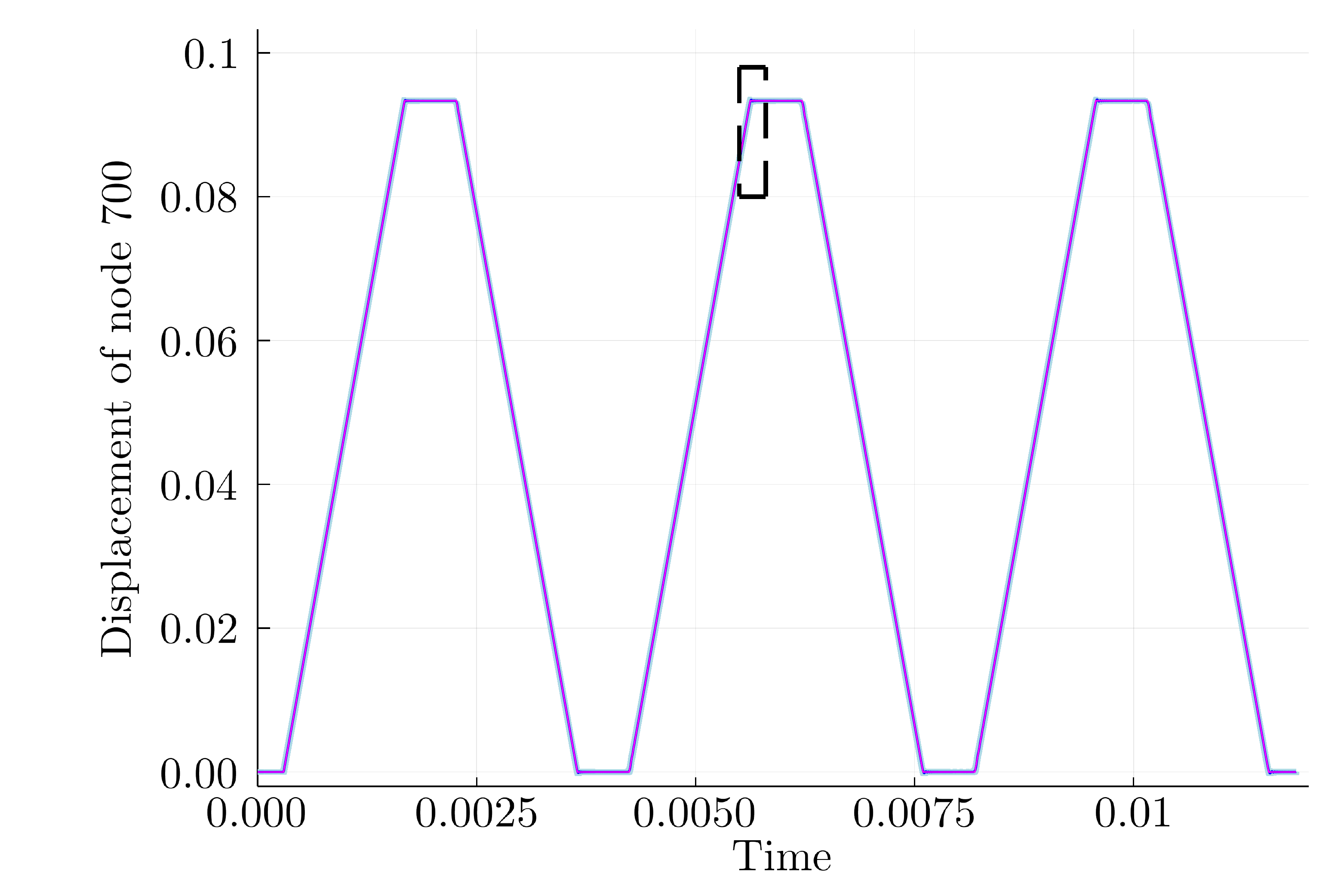}
	\caption{Example 2: Displacement of node 700 as a function of time corresponding to the analytic solution of the PDE (magenta), the analytic solution of the ODE (blue) and the flowpipe (light-blue). The rectangle shows the zoomed region.}
	\label{fig:ex2disp700}
\end{figure}

\begin{figure}[htb]
	\centering
	\includegraphics[width=0.75\textwidth]{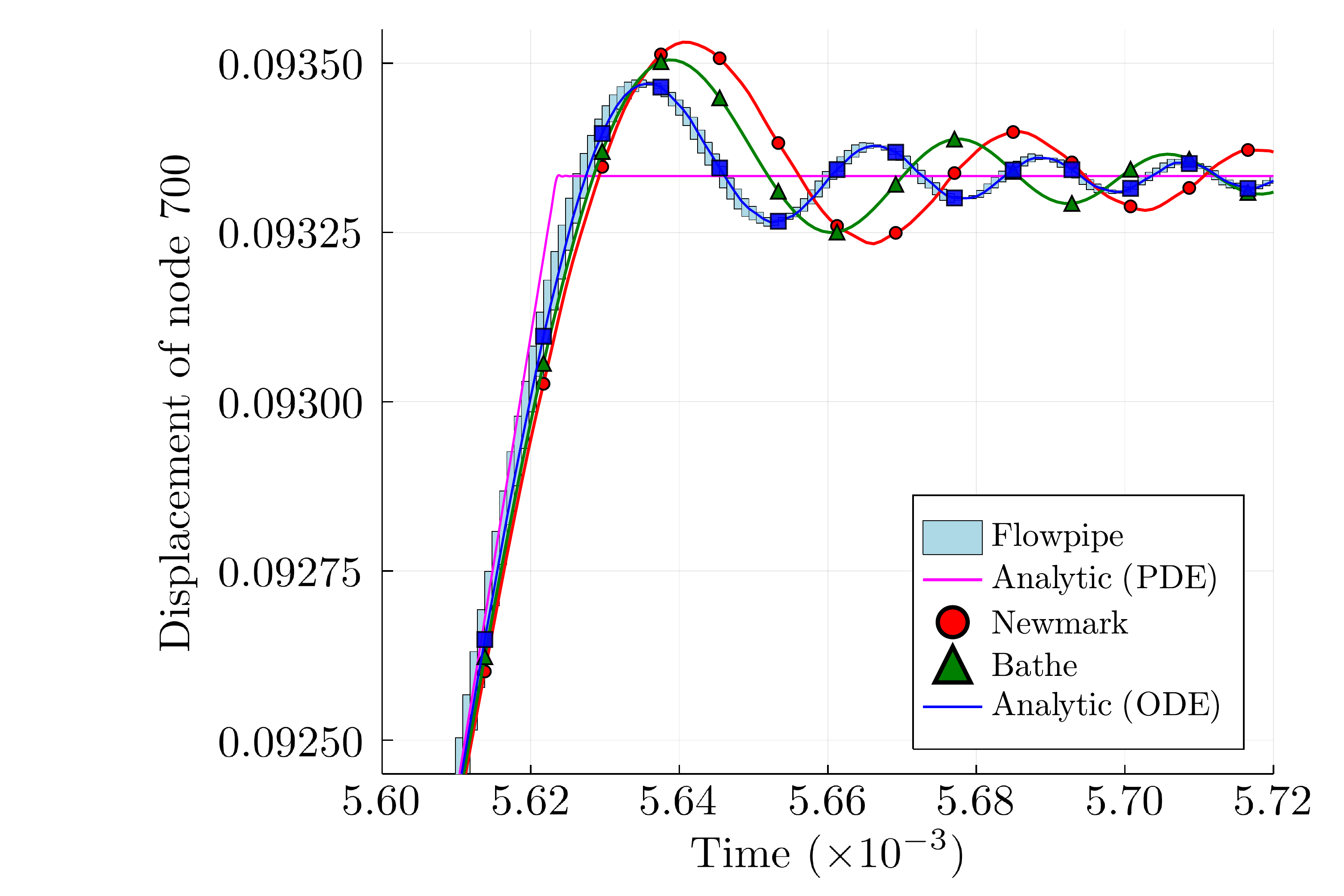}
	\caption{Example 2: Displacement of node 700 of the clamped beam model as a function of time.}
	\label{fig:ex2dispZoom}
\end{figure}

The displacements obtained for node 700 are shown in Fig.~\ref{fig:ex2disp700} and in Fig.~\ref{fig:ex2dispZoom} a zoom at time $5.6\times 10^{-3}$ is shown. %
In this case we present both the analytical solution in the continuum and the analytical solution in the discrete domain. %
The flowpipe computing using set propagation includes the analytic solution in the space-discretized domain, since that is the exact solution of the system of ODEs.

\begin{figure}[htb]
	\centering
	\includegraphics[width=0.75\textwidth]{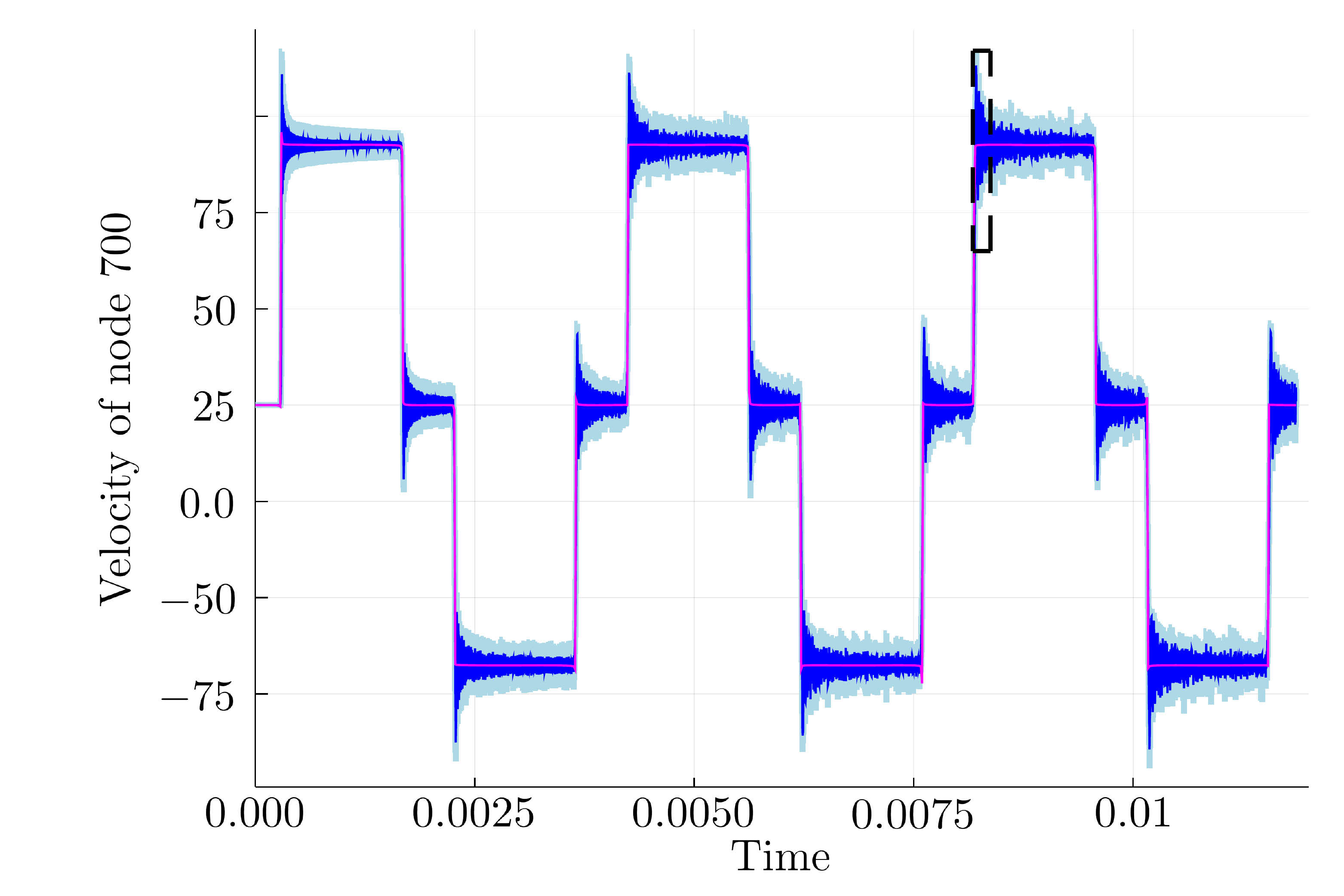}
	\caption{Example 2: velocity of node 700 as a function of time corresponding to the analytic solution of the PDE (magenta), the analytic solution of the ODE (blue) and the flowpipe (light-blue). The rectangle shows the zoomed region.}
	\label{fig:ex2vel}
\end{figure}

\begin{figure}[htb]
	\centering
	\includegraphics[width=0.75\textwidth]{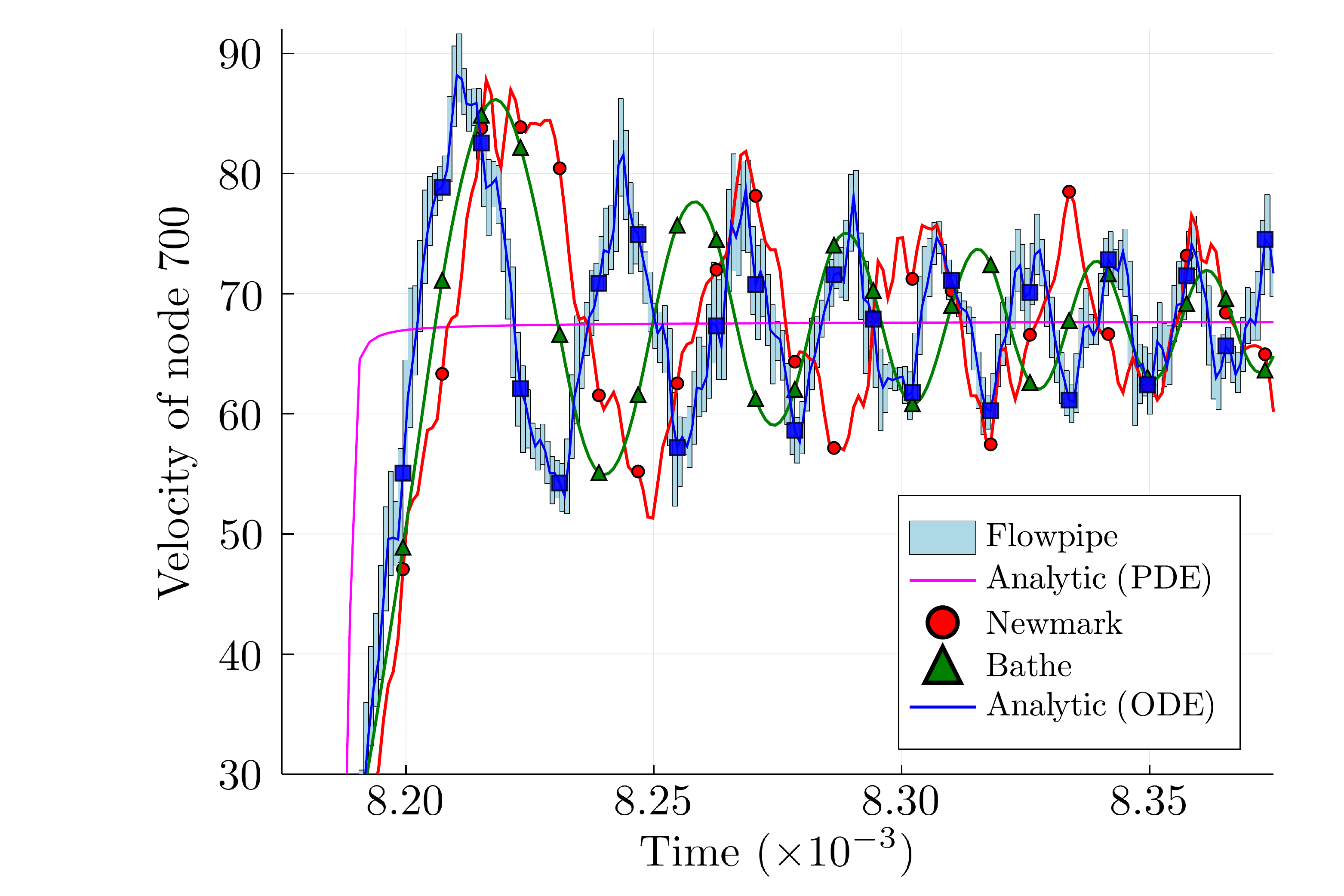}
	\caption{Example 2: zoom for the velocity of node 700.}
	\label{fig:ex2velzoom}
\end{figure}

In Fig.~\ref{fig:ex2vel} the velocities obtained in node 700 are shown, and in Fig.~\ref{fig:ex2velzoom} a zoom is shown. The analytical solution is again within the flowpipe, confirming the theoretically expected results. In this problem the flowpipe has been obtained using support functions along the node of interest. Execution times are reported in Table~\ref{tab:runtimes_clamped}.

\begin{table}[htb]
	\centering
	{\rev
		\begin{tabular}{@{}clll@{}}
			\cmidrule(l){2-4}
			& Set Propagation & Newmark & Bathe \\ \midrule
			& 3.3 s (0.6 GB)  & 0.7 s (2.5 GB)  & 1.3 s (4.3 GB)
		\end{tabular}
		\caption{Execution times and memory usage for the clamped-free bar example.}
		\label{tab:runtimes_clamped}
	}
\end{table}

\clearpage
{\rev 
\subsection{Example 3 - Two-dimensional wave propagation problem} 
}

In this example a two-dimensional wave propagation problem is considered. The geometry and load applied are inspired on the Lamb's Problem considered in \citep{kim2021}. %
The goals of this example are twofold: comparing the results regarding accuracy for a 2D wave propagation problem, and obtaining results for problems with distributed initial conditions.

\subsubsection{Problem definition}

\begin{figure}[htb]
	\centering
	\begin{subfigure}[b]{0.52\textwidth}
		\centering
		\def\svgwidth{\textwidth}
		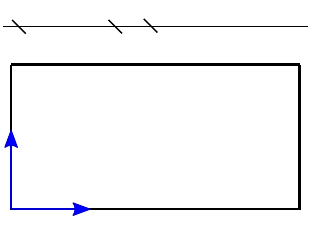 
		\caption{Diagram of domain and boundary conditions considered.}
		\label{fig:ex3Diagram}
	\end{subfigure}
	~~
	\begin{subfigure}[b]{0.38\textwidth}
		\includegraphics[width=.92\textwidth]{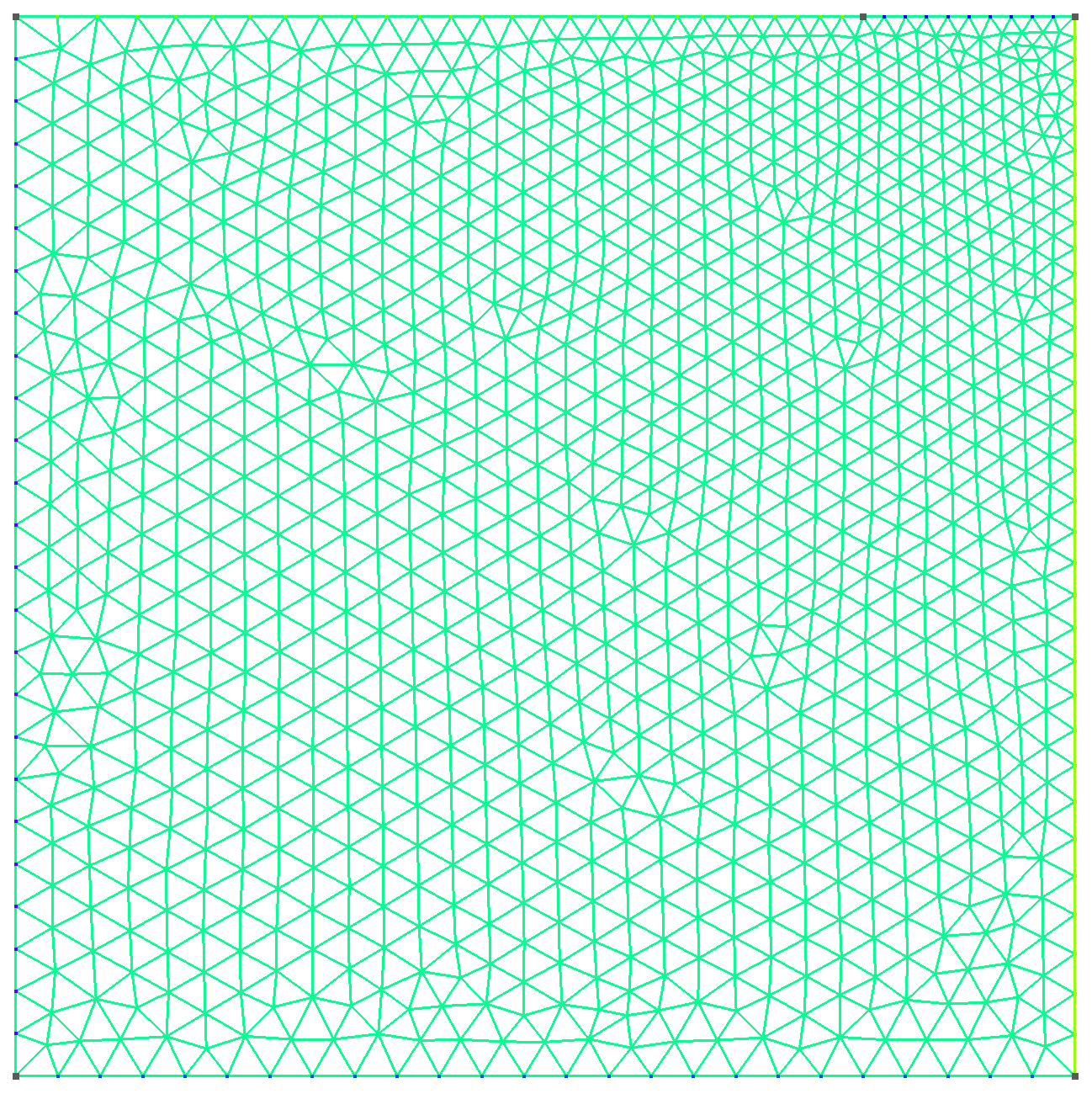}
		\caption{Finite Element Method mesh used, formed by triangular elements.}
		\label{fig:ex3Mesh}
	\end{subfigure}
	\caption{Example 3: diagram and mesh of two-dimensional wave propagation problem considered.}
\end{figure}

The domain geometry and boundary conditions considered are shown in Figure~\ref{fig:ex3Diagram}, where the load is assumed to be given by a step function $P(t) = 2 \times 10^6 H(t)$ (in Newtons), with $H(t)$ being the Heaviside function.
Given the symmetry of the problem, half of the domain is considered. %
The FEM mesh, shown in Figure~\ref{fig:ex3Mesh}, was generated using GMSH and it is formed by 2235 linear triangular elements and 1180 nodes. 

The solid is assumed to be submitted to a plane strain state and formed by a linear elastic isotropic material, with Young modulus $E = 1.8773 \times 10^{10}$ Pa, Poisson ratio $\nu = 0.25$ and  density $\rho = 2200 $ kg/m$^3$. %

The goal of the problem is to compute the horizontal displacement and velocity of the control node $A$, located at $x_A = 2560$~m and $y_A = 3200$~m, as shown in Fig.~\ref{fig:ex3Diagram}.
For the initial conditions, two cases are considered: a single initial condition and distributed initial conditions.

\subsubsection{Single initial condition analysis}

For the first analysis zero displacement and velocity initial conditions are considered, i.e. $\bfu(0) = \bszer$ and $\dot{	\bfu}(0) = \bfv(0) = \bszer$.
For the set propagation method such choice corresponds to $\mcX_0 = \{\bszer\}$ and $\mcC_0 = \{1 \times 10^6\}$, 
with $q = 1$ and $\alpha = 0$ in Eq.~\eqref{eq:input_function_cases}.

\begin{figure}[htb]
	\centering
	\begin{subfigure}[b]{0.48\textwidth} 
		\includegraphics[width=\textwidth]{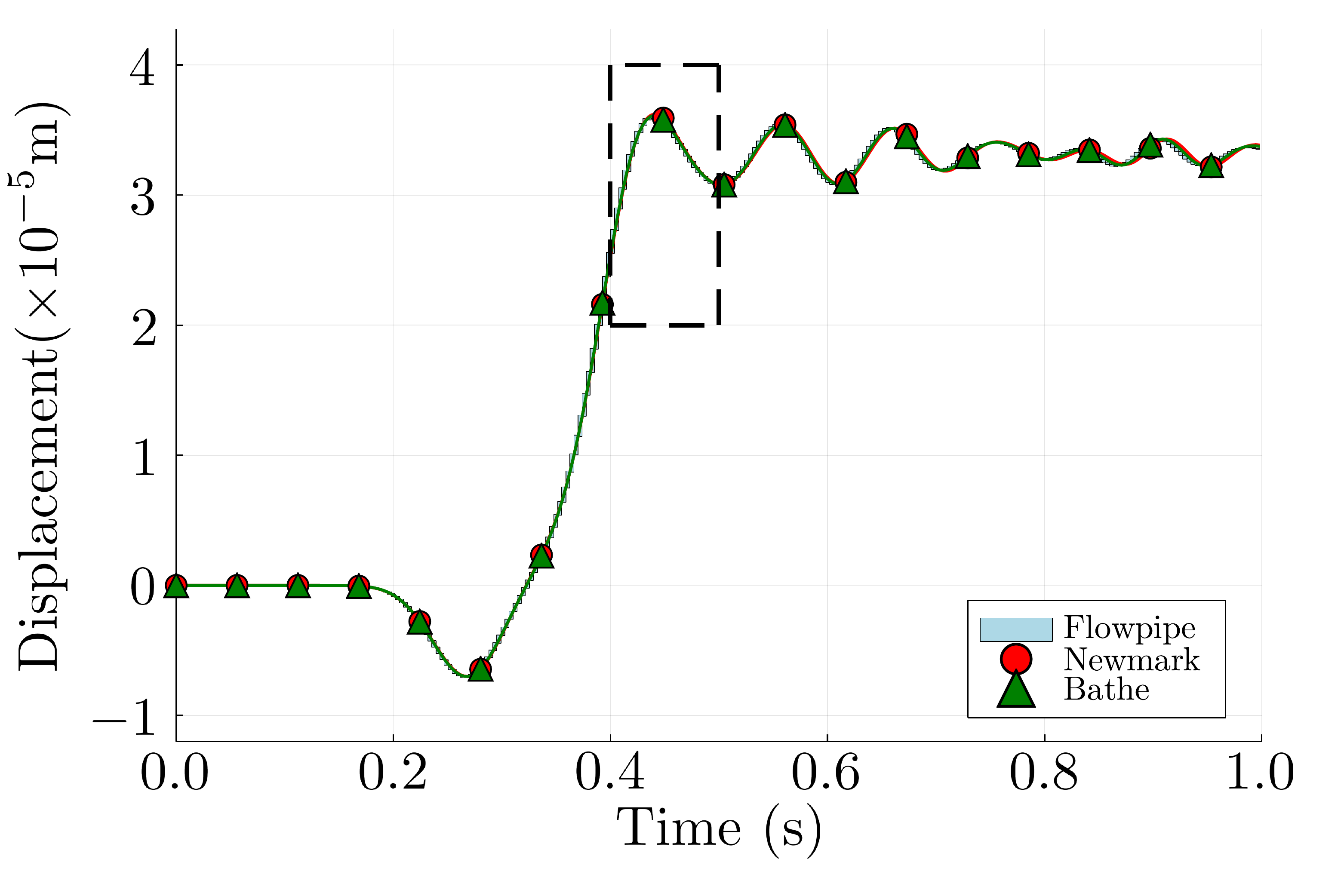}
		\caption{Displacements for $t \in [0, 1]$.}
		\label{fig:wave_singleton_horizontal}
	\end{subfigure} 
~~
	\begin{subfigure}[b]{0.48\textwidth} 
	\includegraphics[width=\textwidth]{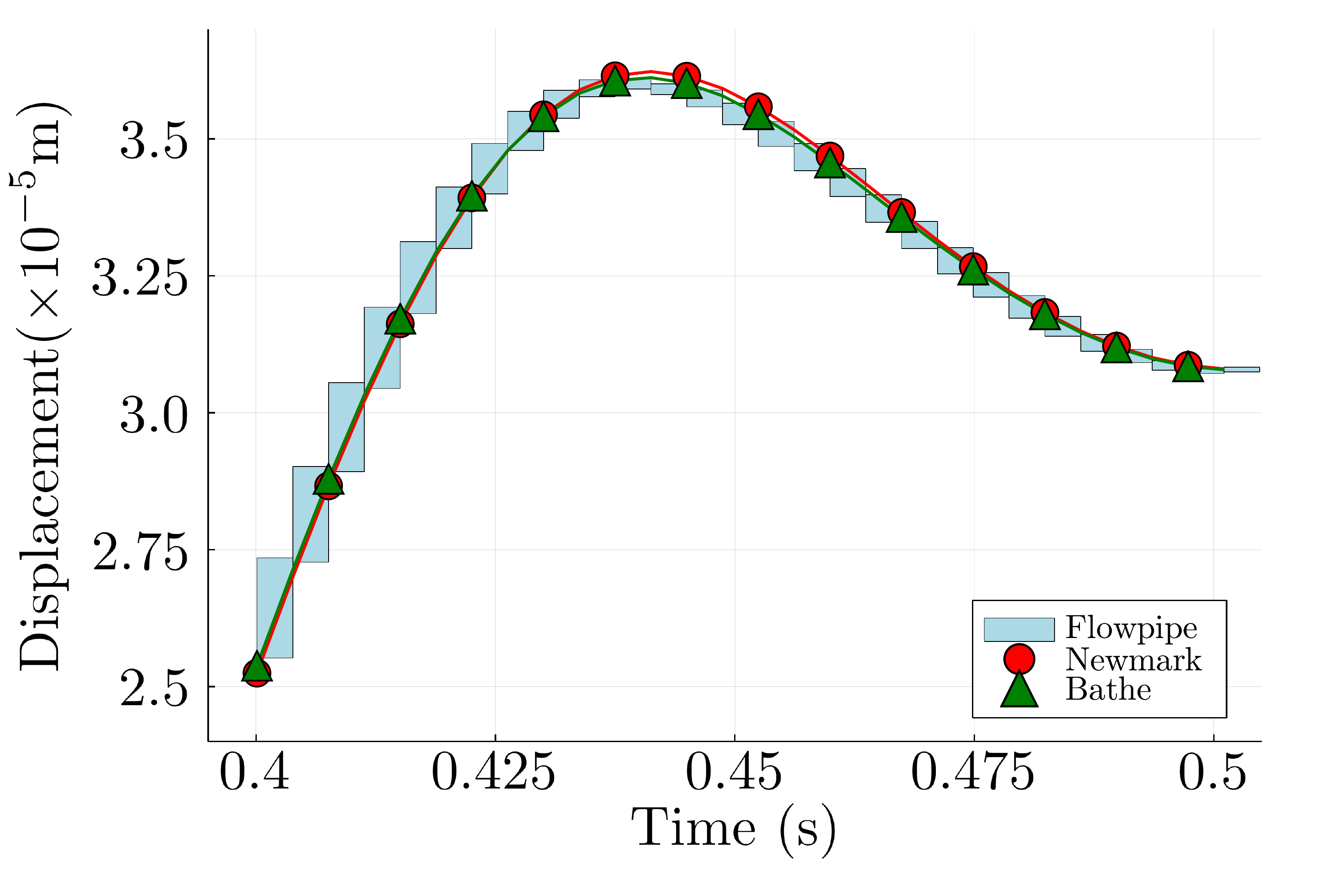}
	\caption{Zoom for $t \in [0.4, 0.5]$.}
	\label{fig:wave_singleton_horizontal_zoom}
\end{subfigure}

\caption{Example 3: Horizontal displacement at control node A for the wave propagation problem subject subject to a step-load and a single initial condition.}
\end{figure}

\begin{figure}[htb]
	\centering
	\begin{subfigure}[b]{0.48\textwidth}
		\includegraphics[width=\textwidth]{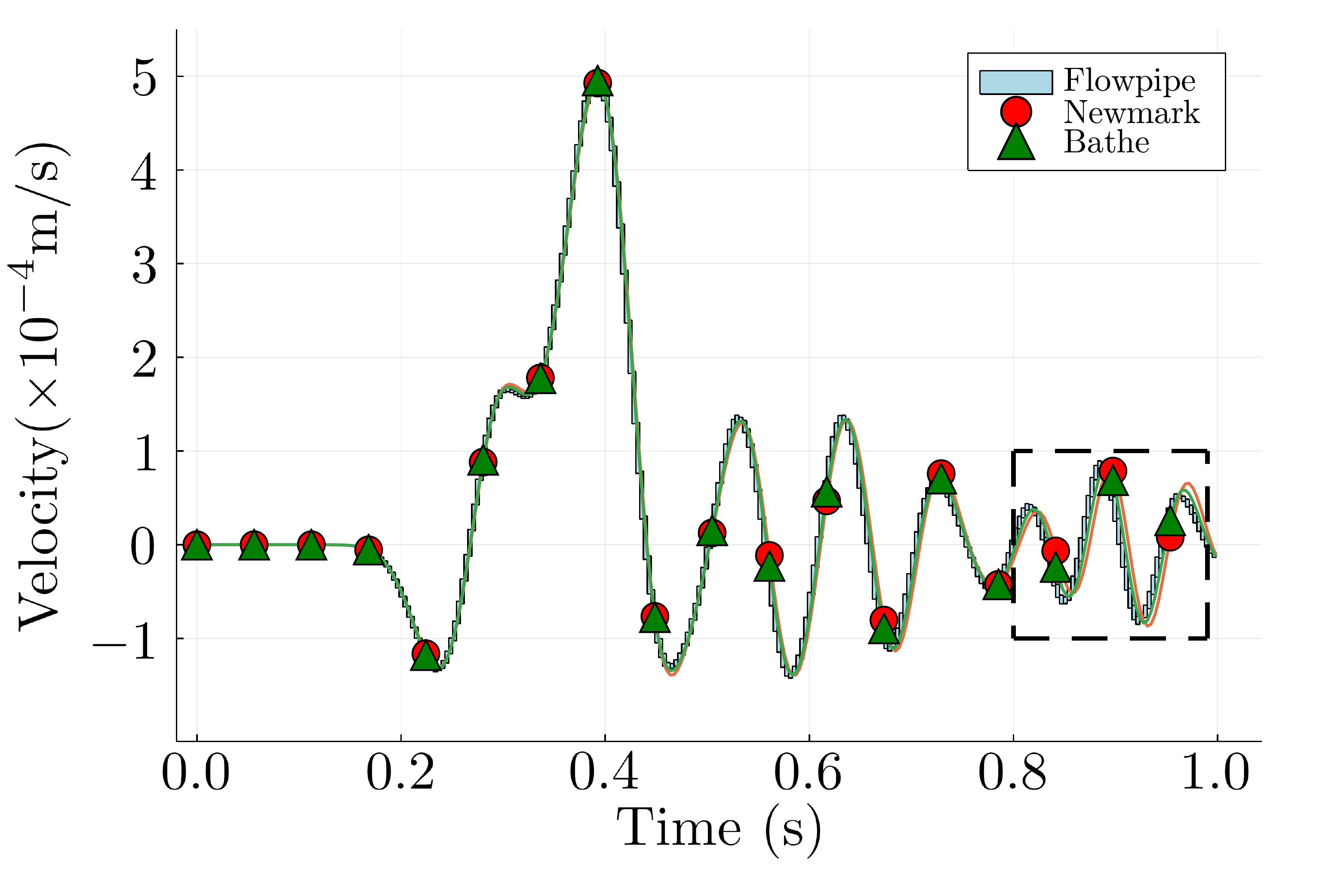}
		\caption{Velocities for $t \in [0, 1]$.}
		\label{fig:wave_singleton_horizontal_vel}
	\end{subfigure}
	~~
	\begin{subfigure}[b]{0.48\textwidth}
		\includegraphics[width=\textwidth]{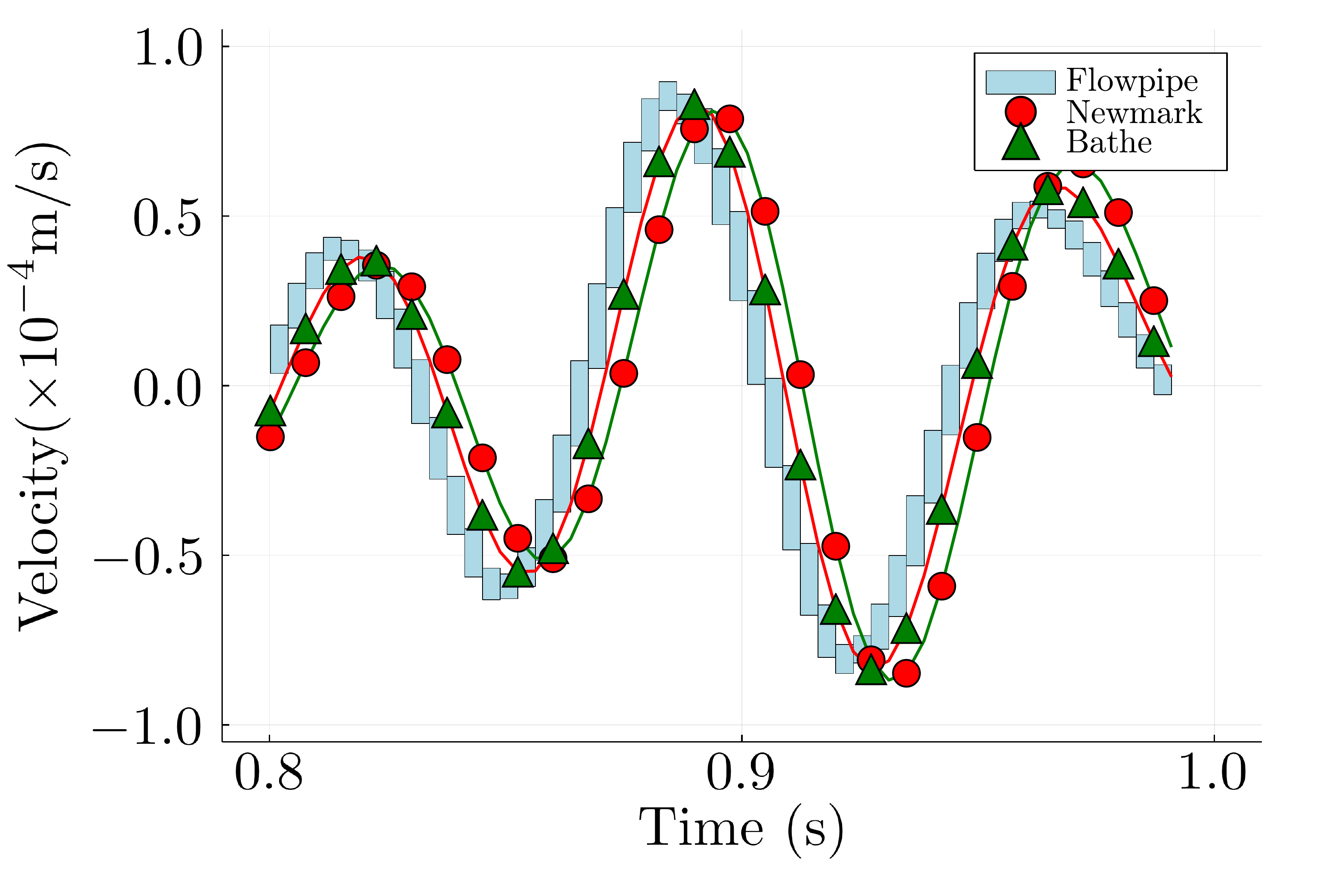}
		\caption{Zoom for $t \in [0.8, 1.0]$.}
		\label{fig:wave_singleton_horizontal_vel_zoom}
	\end{subfigure}
	
    \caption{Example 3: Horizontal velocity at control node A for the wave propagation problem subject to a step-load and a single initial condition.}
\end{figure}

Given that the goal of the problem is to compute the evolution of a few state variables of the system, support functions are chosen for the set propagation. %
The time step size was chosen using the formula proposed in \citep{kim2021}, then
$\delta = \sqrt{\frac{2\cdot 3200^2}{2235}} \cdot \frac{0.125}{3200} \approx 0.0037$. The set propagation result was obtained using support functions along the displacements and velocity components of node $A$.
The obtained horizontal displacements of node A for the Set Propagation, Newmark and Bathe methods are shown in Fig.~\ref{fig:wave_singleton_horizontal}, where in Fig.~\ref{fig:wave_singleton_horizontal_zoom} a zoom for time interval $[0.4,0.5]$ is shown. The corresponding plots for the velocity of node A are shown in Fig.~\ref{fig:wave_singleton_horizontal_vel}. In Fig.~\ref{fig:wave_singleton_horizontal_vel_zoom} the zoom for the velocities at time interval $[0.8,1.0]$ is  shown, where it can be observed that numerical methods, and in particular Newmark, provide solutions outside of the Flowpipes (as in the previous example).

The runtimes (and memory usage) of the methods are: Newmark 3.7 s (0.6 GB), Bathe 4.7 s (1.0 GB) and Set Propagation 5.6 s (1.1 GB).

\subsubsection{Distributed initial conditions analysis}

Let us now consider a case in which the values of initial displacement and velocity of each node are assumed to be included in the intervals $[-\delta_u,\delta_u]$ and $[-\delta_v,\delta_v]$, respectively. %
This is equivalent to $\bfu(0) \in \mcU_0 = [-\delta_u,\delta_u]^{N_{dof}}$ and $\dot{\bfu}(0) = \bfv(0) \in \mcV_0 = [-\delta_v,\delta_v]^{N_{dof}}$, where $N_{dof}=2322$ is the number of degrees of freedom of the system. %
In this problem, the radii of the uncertainty intervals, $\delta_u$ and $\delta_v$, are considered as the 10 \% of the maximum magnitude of displacement and velocity of node A for the previous singleton case, i.e. $3.5 \times 10^{-7}$ m and $5 \times 10^{-6}$ m/s.

Each initial condition $(\bfu_0,\bfv_0)$ has a corresponding solution for the horizontal velocity at node A, denoted as $v_{x,A}(\bfu_0,\bfv_0,t)$. %
The maximum envelope and minimum envelope functions are defined as:
\begin{equation}
v_{env}^+(t) =
\max_{\bfu_0 \in \mcU_0, \, \bfv_0 \in \mcV_0} v_{x,A}(\bfu_0,\bfv_0,t),
\quad
\text{and}
\quad 
v_{env}^-(t) =
\min_{\bfu_0 \in \mcU_0, \, \bfv_0 \in \mcV_0} v_{x,A}(\bfu_0,\bfv_0,t),
\end{equation}
respectively, and the absolute envelope is defined as:
\begin{equation}
v_{env}(t) =
\max \{ |v_{env}^+(t)|, |v_{env}^-(t)| \}.
\end{equation}
We define the following metrics to evaluate the numerical results:
\begin{equation}
\left\| 
v_{env}
\right\|_{L_\infty}
= \max_{t \in [0,1]} v_{env}(t),
\qquad
\left\| v_{env} \right\|_{L_1} 
= \int_{0}^{1} \left| v_{env}(t) \right| dt.
\end{equation}

The Set Propagation method propagates all the initial conditions, thus the computation of $v_{env}^+$, $v_{env}^-$ and $v_{env}$ can be done executing a single analysis for $\mcX_0=\mcU_0\times\mcV_0$. %
For numerical integration methods, the most direct approach to compute this would be to perform the simulations for a given number of random initial points in the initial conditions set. %
A more efficient approach is to consider only the vertices of the initial states. %
It can be demonstrated (and it was empirically verified) that this selection strategy provides solutions with higher $v_{env}$ values, which seems the best scenario for numerical integration methods. %
The number of different initial conditions considering this vertex constraint is $2^{2322} \approx 10^{698}$, hence only an small fraction of such solutions can be practically evaluated.

Only the Newmark method is used as numerical integration method, considering that similar results would be obtained using Bathe's method. %
The time-step used is the same as in the previous case. %
The numerical values of $v_{env}^+$ and $v_{env}^-$ are computed for the first: 1, 10, 100, 1000, 10000, and 100000 intial conditions, sampling uniformly from the vertices of $\mcU_0 \times \mcV_0$.

\begin{figure}[htb]
	\centering
	\includegraphics[width=0.8\textwidth]{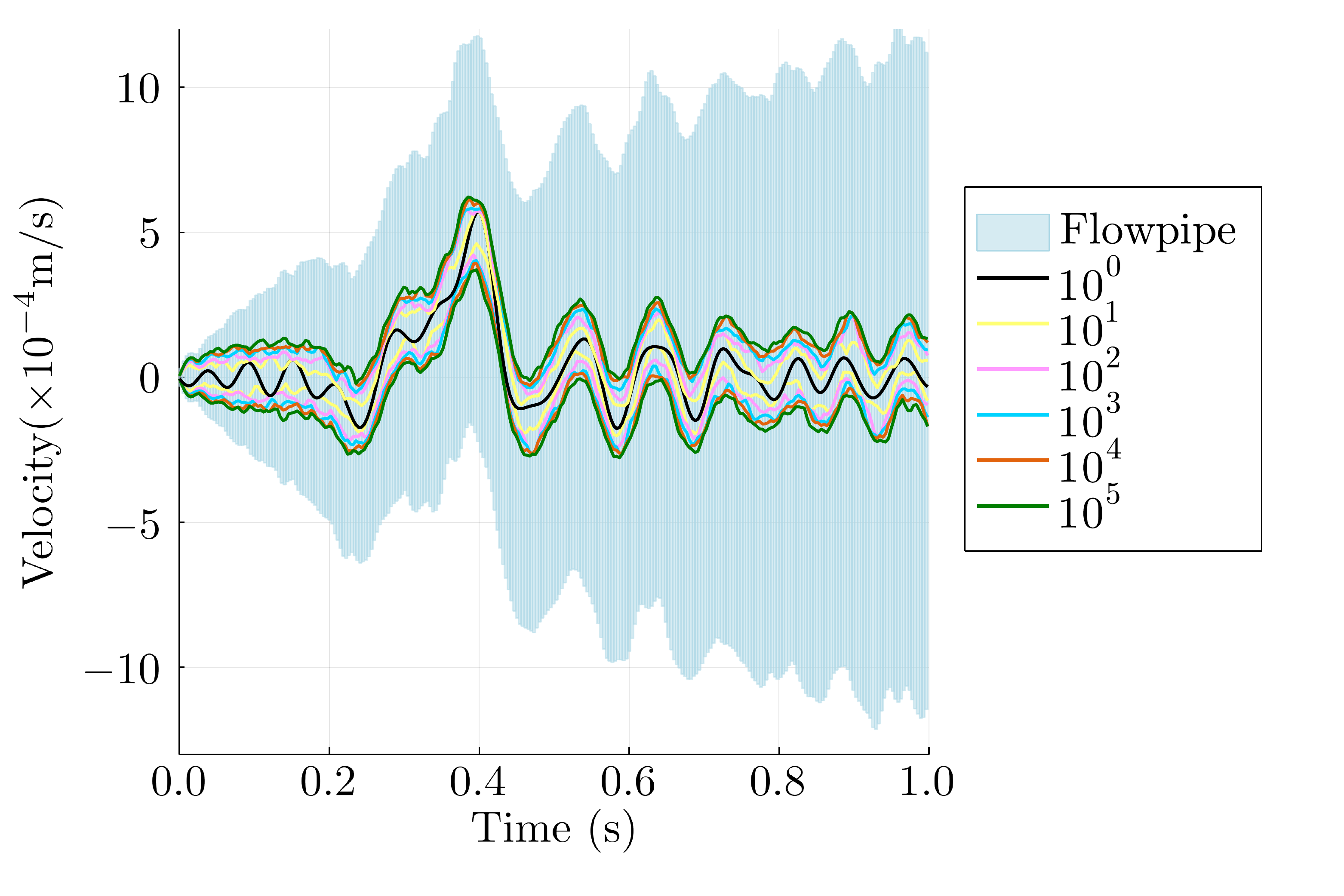}
	\caption{Example 3: plot of the flowpipe obtained by Set Propagation and the envelope functions $v_{env}^+$ and $v_{env}^-$ for  1, 10, 100, 1000, 10000, and 100000 random intial conditions using Newmark's method.}
	\label{fig:ex3Envelopes}
\end{figure}

In Fig.~\ref{fig:ex3Envelopes}, the functions $v_{env}^+$ and $v_{env}^-$ obtained for the Newmark method and the flowpipe obtained using the Set Propagation method are shown.

The results obtained for the norms of $v_{env}$ and the execution times are shown in Table~\ref{tab:waveprop}. %
It can be observed that, as the number of trajectories increases, the norms of the solutions obtained by Newmark, also increase with no clear convergence. %
On the other hand, the Set Propagation method provides an upper bound that covers all the admissible trajectories.

\begin{table}[htb]										
	\centering										
\begin{tabular}{lrccc}										\hline
	Method	&	\# Trajectories	&	Time (s)	&	$\left\| v_{env}\right\|_{L_1}$ (10$^{-5}$)	&	$\left\| v_{env}\right\|_{L_\infty}$ (10$^{-5}$)	\\	\hline
	Newmark	&	1	&	0.3	&	9.27	&	56.98	\\	\hline
	Newmark	&	10	&	2.0	&	13.52	&	57.53	\\	\hline
	Newmark	&	100	&	17.7	&	16.61	&	57.59	\\	\hline
	Newmark	&	1000	&	175.5	&	18.52	&	58.22	\\	\hline
	Newmark	&	10000	&	1771.4	&	19.98	&	61.18	\\	\hline
	Newmark	&	100000	&	17796.1	&	21.42	&	62.21	\\	\hline
	Set Propagation	&	-	&	8.5	&	81.33	&	122.25		
\end{tabular}										
	\caption{Example 3: results of execution time and memory usage for distributed initial conditions using Newmark method and Set Propagation method.}\label{tab:waveprop}										
\end{table}

\subsection{Example {\rev4 } - One-dimensional heat {\rev transfer problem} } \label{sec:Example3_Heat}

In this example, a one-dimensional heat transfer problem is considered and interval sets are introduced for the initial values.
The main goal of this problem is to evaluate the set propagation method as a tool to analyze the maximum temperature of a heat transfer problem under input uncertainty.

\subsubsection{Problem definition}

The domain consists of the interval $\Omega=[0,1]$ and homogeneous Dirichlet boundary conditions are considered on both ends. %
The domain is assumed to have $\kappa = 1$ as thermal conductivity, $\rho = 1$ as density and $c = 1$ as specific heat capacity, all assumed uniform and constant.

The governing equations can be written as:
\begin{equation}
	\left\{
	\begin{array}{lr}
		\dfrac{\partial T(x,t)}{\partial t} = \dfrac{\kappa}{\rho c} \dfrac{\partial^2 T(x,t)}{\partial x^2}, & x \in \Omega\\
		T(x, t) = 0, & x\in \{0,1\},
	\end{array}
	\right.
	\label{eq:heat1d}
\end{equation}
%
where $T(x, t)$ is the temperature in position $x$ at time $t$. %

For the initial conditions a set of possible temperature profiles is considered. %
The expression of the initial condition function is given by:
\begin{equation}
	\displaystyle T(x,0) = (1+\varepsilon) \left( \sin(\pi x) + \dfrac{1}{2} \sin(3\pi x) \right), \qquad \varepsilon \in [-0.1,0.1],
	\label{eq:heat1d_initial}
\end{equation}
where $\varepsilon $ is a relative difference parameter associated with a given level of uncertainty in the initial values used for the simulations.

The Finite Element Method is used to discretize the governing equations considering a 100 two-node elements mesh with linear interpolation. %
It is important to remark that although the FEM discretization is used, the set of initial temperatures contains an infinite number of functions since $\varepsilon$ belongs to an interval that is dense in $\mathbb{R}$.

\subsubsection{Maximum temperature analysis}

The results obtained using the set propagation and the Backward Euler methods are compared with the analytic (exact) solution of Eq.~\eqref{eq:heat1d}.
The time-step used for both resolution methods is $\delta = 10^{-5}$. %
In the case where set propagation is used, the entire set of parameters $\varepsilon$ from Eq.~\eqref{eq:heat1d_initial} is considered, while for the Backward Euler method the analysis is only done for the two extreme cases: $\varepsilon = - 0.1$ and $\varepsilon =  0.1$.

\begin{figure}
	\centering
	\begin{subfigure}[b]{0.45\textwidth}
		\centering
		\includegraphics[width=\textwidth]{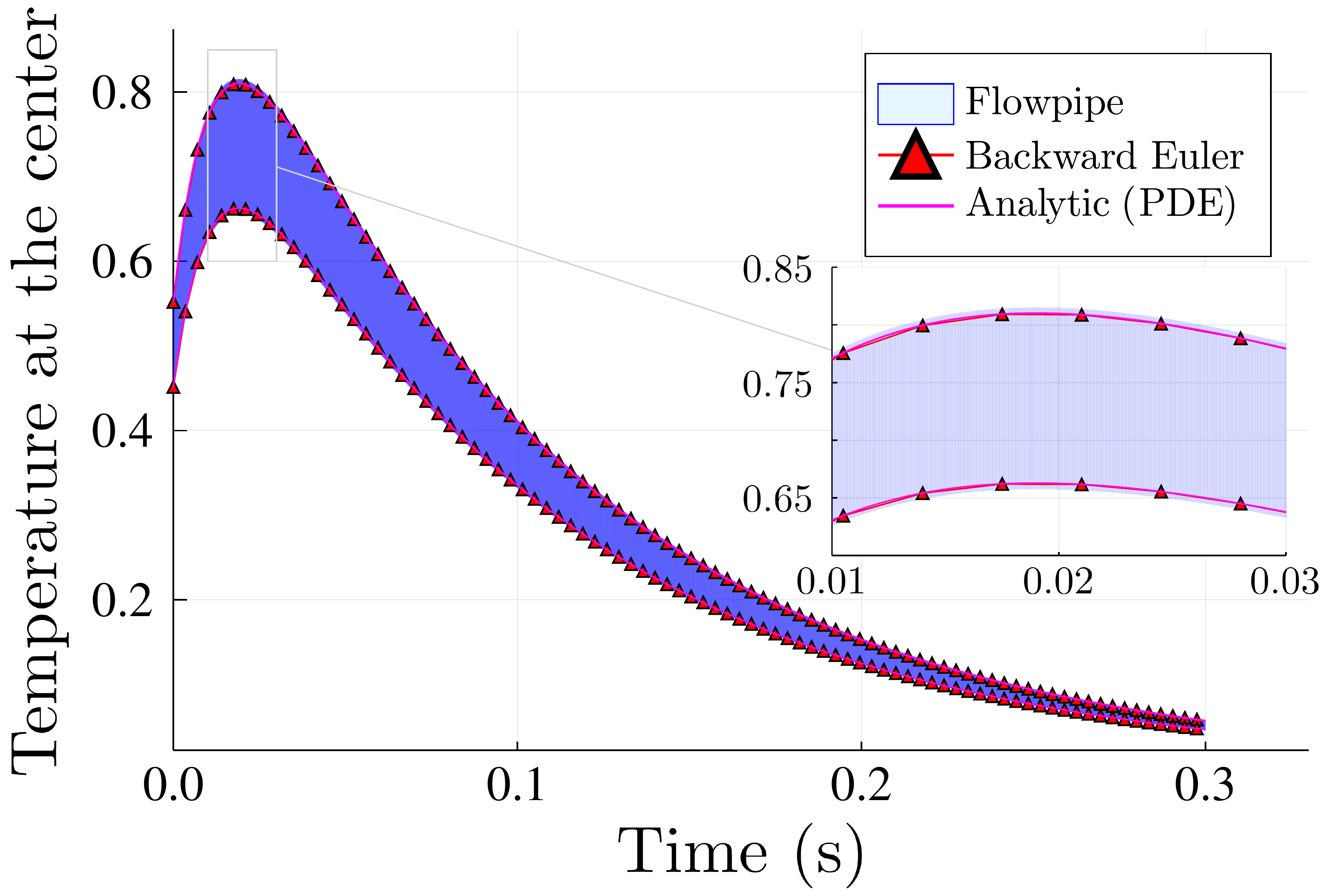}
		\caption{Temperature evolution obtained at the center node of the mesh ($x=0.5$) using the analytic solution (magenta), the Backward Euler method for extreme cases (red triangle), and the set propagation method (blue flowpipe). 
		}
		\label{fig:heat1d_center}
	\end{subfigure}
	~~
	\begin{subfigure}[b]{0.45\textwidth}
		\centering
		\includegraphics[width=\textwidth]{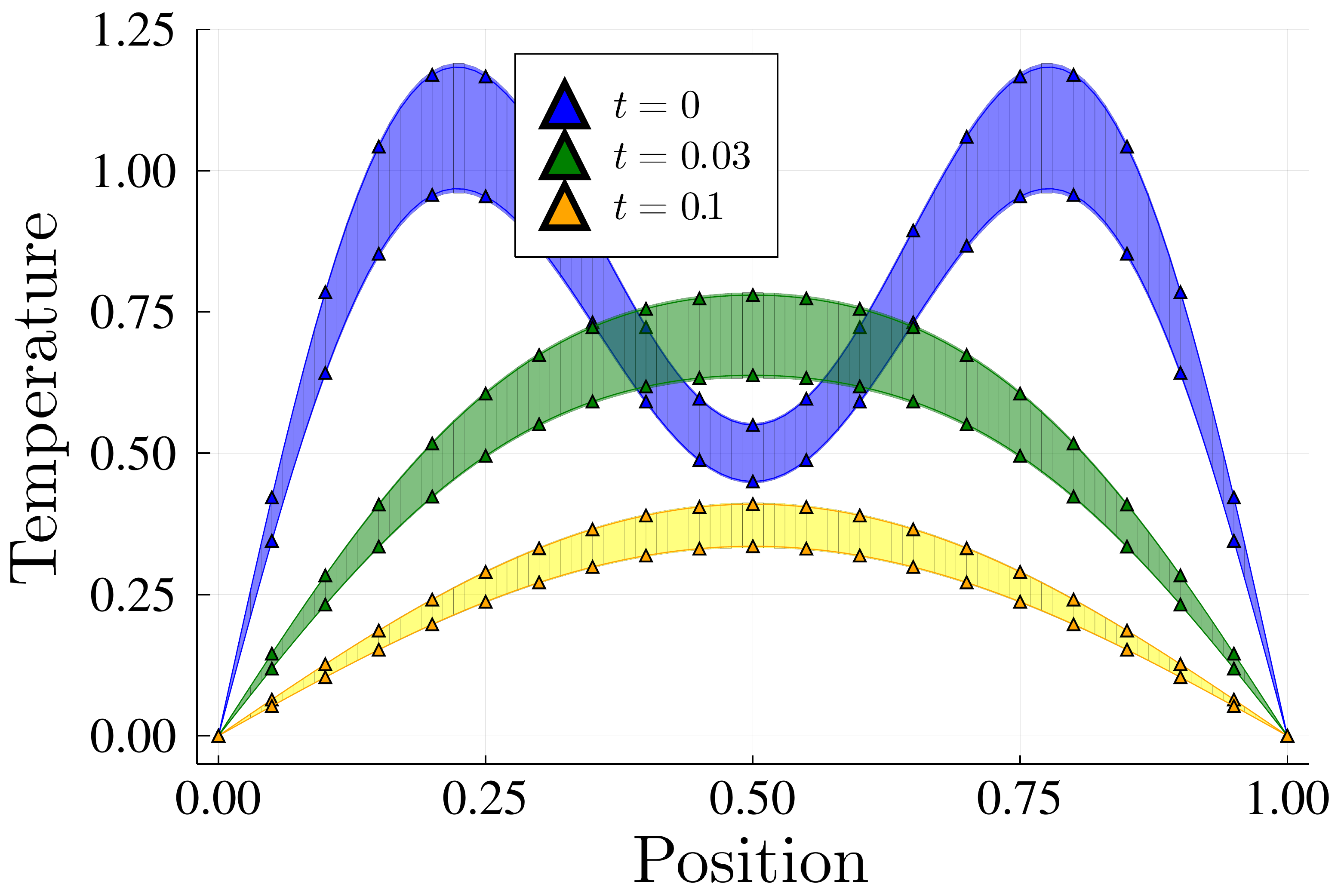}
		\caption{Temperature profiles computed at different time points: $t=0$ (blue), $t=0.03$ (green) and $t=0.1$ (yellow). Results from the Backward Euler method at extreme cases are shown using triangle markers.}
		\label{fig:heat1d_nodes}
	\end{subfigure}
	\caption{Example 4: Results for the one-dimensional heat problem with set initial conditions obtained with the set propagation method and with the Backward Euler method.}
	\label{fig:heat1d}
\end{figure}

The results for the temperature in the node $x=0.5$ are shown Fig.~\ref{fig:heat1d_center} and temperature profiles at times $t = 0$, $0.03$ and $t = 0.10$ are shown in Fig.~\ref{fig:heat1d_nodes}.
We compute the flowpipe using a hyperrectangular representation. Each direction corresponds to a different node.
We observe that the flowpipe borders contain and tightly match the numerical solutions from the Backward Euler method and the analytical solutions starting from the two extreme initial conditions. 

\subsubsection{Maximum temperature gradient analysis}

In this case, the goal of the analysis is to obtain the maximum value of the temperature gradients among all the nodes of the domain for all the admissible initial condition profiles.

Since it is not clear how to choose the initial temperature profile to obtain an extreme temperature gradient at a given time, here we repeated the calculations with Backward Euler using several different initial temperature profiles (a thousand) chosen randomly from Eq.~\eqref{eq:heat1d_initial}. %
The Set Propagation method allows to compute all the solutions for the gradients, in this case, using support functions.

\begin{figure}[htb]
	\centering
	\begin{subfigure}[b]{0.45\textwidth}
		\centering
		\includegraphics[width=\textwidth]{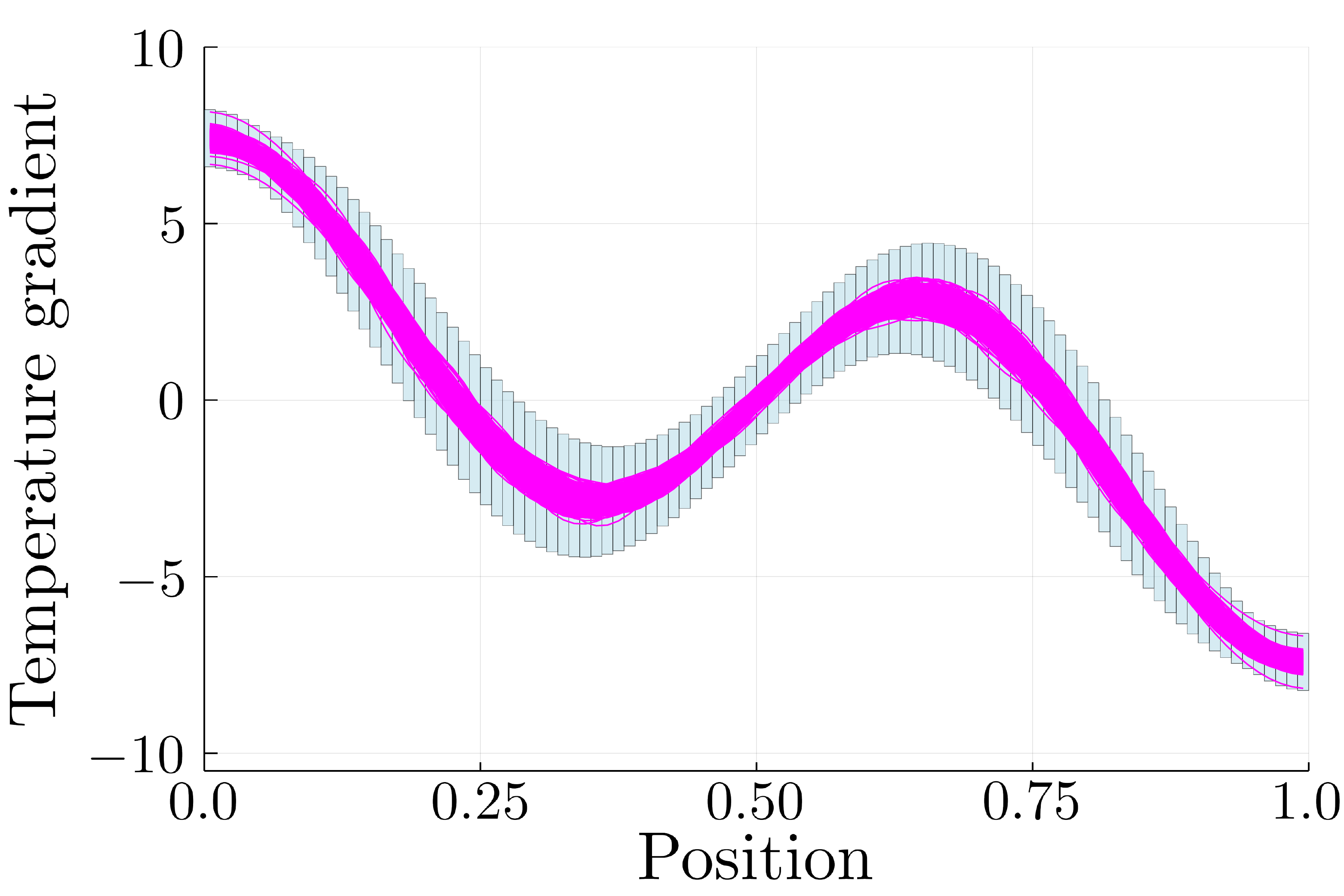}
		\caption{Temperature gradient profile at time $t = 0.001$.}
		\label{fig:heat1d_perfil_grad_temp}
	\end{subfigure}
	~~
	\begin{subfigure}[b]{0.45\textwidth}
		\centering
		\includegraphics[width=\textwidth]{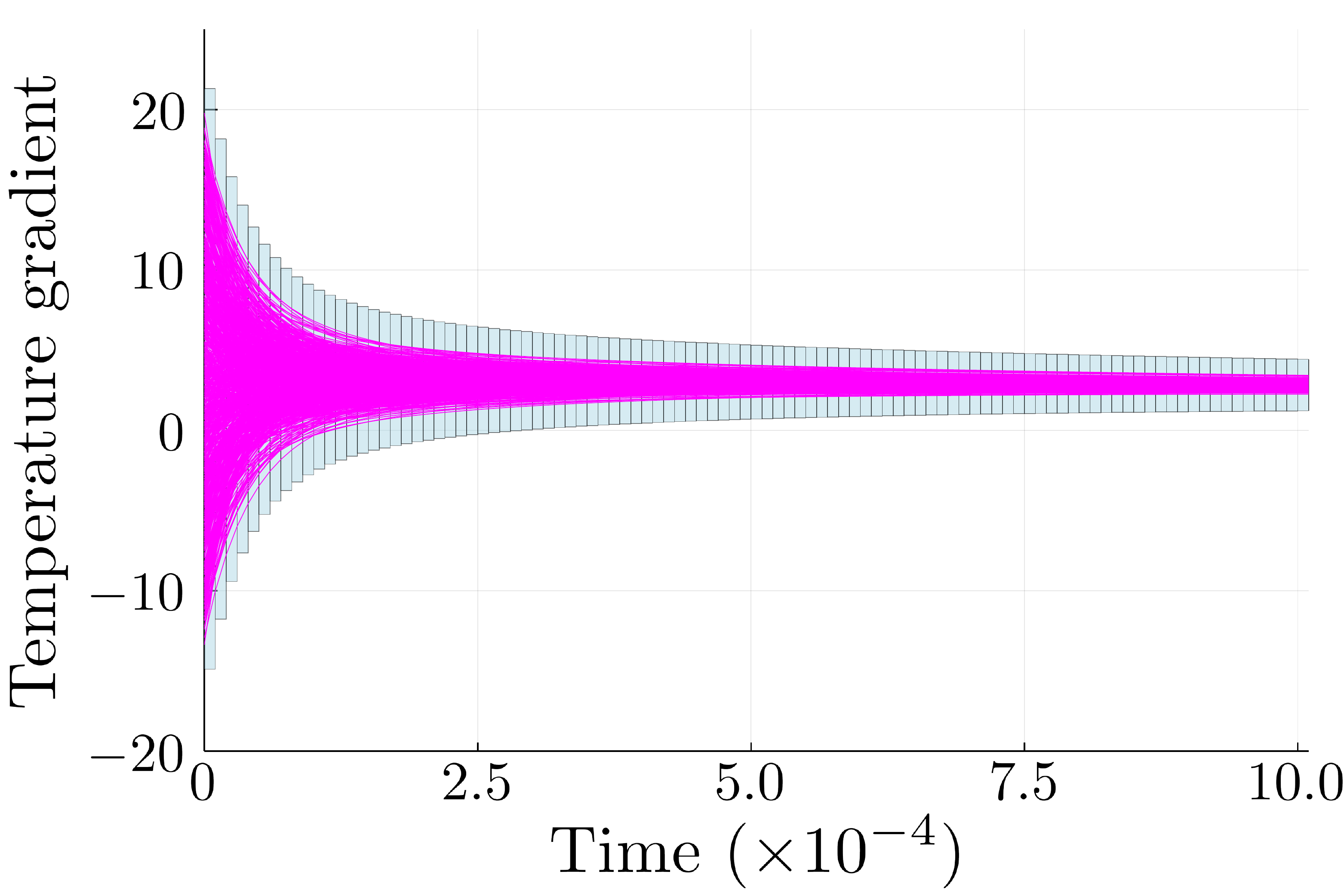}
		\caption{Temperature gradient vs. time at $x = 2/3$.}
		\label{fig:heat1d_grad_tiempo}
	\end{subfigure}
	\caption{Example 4: Flowpipe (light blue) calculation of the temperature gradient for the same case study from Fig.~\ref{fig:heat1d}. Superposed are shown several temperature gradients computed with the Backward Euler method (magenta) with different initial states chosen randomly within the given uncertainty in the initial temperature profile.}
	\label{fig:heat1d_gradientes}
\end{figure}

The results obtained for the temperature gradient profile at time $t=0.001$ are shown in Fig.~\ref{fig:heat1d_perfil_grad_temp} and the time evolution of the temperature gradient at $x = 0.66$ are shown in Fig.~\ref{fig:heat1d_grad_tiempo}. %
We observe that the flowpipe contains all the trajectories calculated with Backward Euler. We remark that statistical strategies are needed to estimate the temperature gradient bounds if we use the Backward Euler method. %
On the other hand, the flowpipe provides safe bounds for this gradient with only one integration. We conclude that the Set Propagation method is an effective tool to evaluate the maximum temperature under uncertainty.

The execution times are reported in Table~\ref{tab:runtimes_heat1d}.

\begin{table}[htb]
	\centering
	{\rev 
		\begin{tabular}{@{}clll@{}}
			\cmidrule(l){2-4}
			& Set Propagation & B. Euler (1 trajectory) & B. Euler (1000 trajectories) \\ \midrule
			Profile & 2.6 (113 MB)  & 0.42 s (107 MB)  & - \\ \midrule
			Gradient & 1.1 s (22 MB)  & - & 17 s (7.1GB)
		\end{tabular}
		\caption{Execution times and memory usage for the one-dimensional heat problem.}
		\label{tab:runtimes_heat1d}
	}
\end{table}

\subsection{Example 5 - concrete casting heat of hydration} \label{sec:Example4_Heat3D}

In this example a model for heat of hydration during casting of a massive concrete structure is considered. %
The development of reliable computational models for this kind of processes is a greatly challenging task due to the uncertainty in parameters such as the internal heat generation. %
In \citep{Wilson1974}, for instance, different levels of heat of hydration (associated with different sizes of aggregate) are considered. %
The parameter values and the modeling assumptions considered are based on \citep{Bofang2014,Tahersima2017}. %

\subsubsection{Problem definition}

\begin{figure}[htb]
	\centering
	\begin{subfigure}[b]{0.55\textwidth}
		\def\svgwidth{\textwidth}
		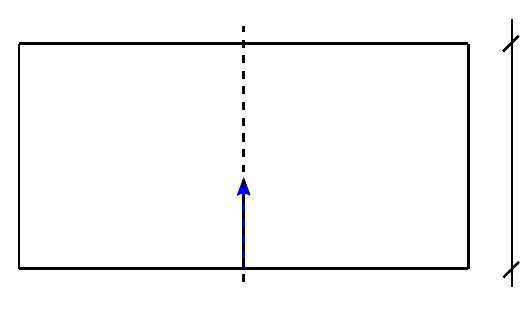
		\caption{Lateral view, control points: A located at $(0,0,0.6 \, \text{m})$, B located at $(0,0,0.9 \, \text{m})$.}
	\end{subfigure}
	\hfill
	\begin{subfigure}[b]{0.35\textwidth}
		\def\svgwidth{\textwidth}
		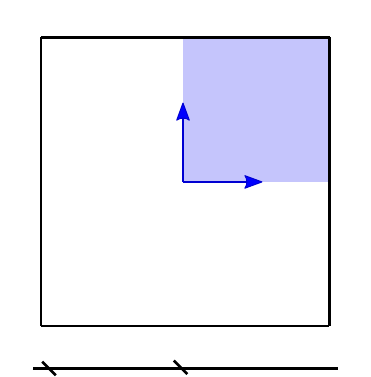
		\caption{Floor plan}
	\end{subfigure}
	\caption{Example 5: diagram of domain, boundary conditions and control points.}
	\label{fig:ex4diag}
\end{figure}

The domain considered consists in a square cuboid $\Omega = [-1,1]\times[-1,1]\times[0,1]$, with $1$ m height (in the $z$ direction) and $2$~m width in the other directions ($x$ and $y$). %
In Figure~\ref{fig:ex4diag} the geometry, boundary conditions and control points of the problem are shown.

The concrete density, the specific heat and the conductivity are considered $\rho=2485 $ kg/m$^3$, $c= 0.967$ kJ/(kg $^\circ$C) and $\kappa = 9.37 $ kJ/(m h $^\circ$C), respectively.

Regarding the boundary conditions, on the upper face ($z=1$) a convection condition is assumed with value $h_{air} = 40$ kJ/(m$^2$ h $^\circ$C), while on the other faces it is $h_{timb}= 500$ kJ/(m$^2$ h $^\circ$C). %
The external temperature $T_\infty(t)$ dependence in time is assumed to be given by the following sinusoidal law:
\begin{equation}
	T_\infty (t) = T_{min} + T_{var} \left( \frac{1}{2} + \frac{1}{2} \sin \left( \omega t - \frac{\pi}{2} \right) \right),
\end{equation}
where $T_{min} = 17 \, ^\circ $C is the minimum temperature, $T_{var}$ is the variation or difference between the minimum and maximum temperatures and $\omega = \pi/12~\textrm{h}^{-1}$ is the angular frequency.
The initial temperature of all the points in the concrete is assumed $T(\bfx,0) = T_{min}$, which is equal to the initial external temperature $T_{\infty}(0)$. %

The accumulated heat of hydration is assumed to be given by the exponential law:
\begin{equation}
	Q_{AH}(t) = Q_{FH} (1-e^{-mt}),
\end{equation}
where $Q_{FH}$ is the final heat of hydration in kJ/kg and $m = 7.95 \times 10^{-3}~\textrm{h}^{-1}$.
%
This quantity is multiplied by the density, derived and included in the thermal-work Eq.~\eqref{eqn:heatfem_force} as:
\begin{equation}
	Q_{int}(t) = Q_{FH} \cdot \rho \cdot m e^{-mt}.
\end{equation}

\subsubsection{Numerical resolution with fixed parameters} \label{sec:example4_fixed_params}

\begin{figure}[htb]
	\centering
	\begin{subfigure}[b]{0.45\textwidth}
		\centering
		\includegraphics[width=\textwidth]{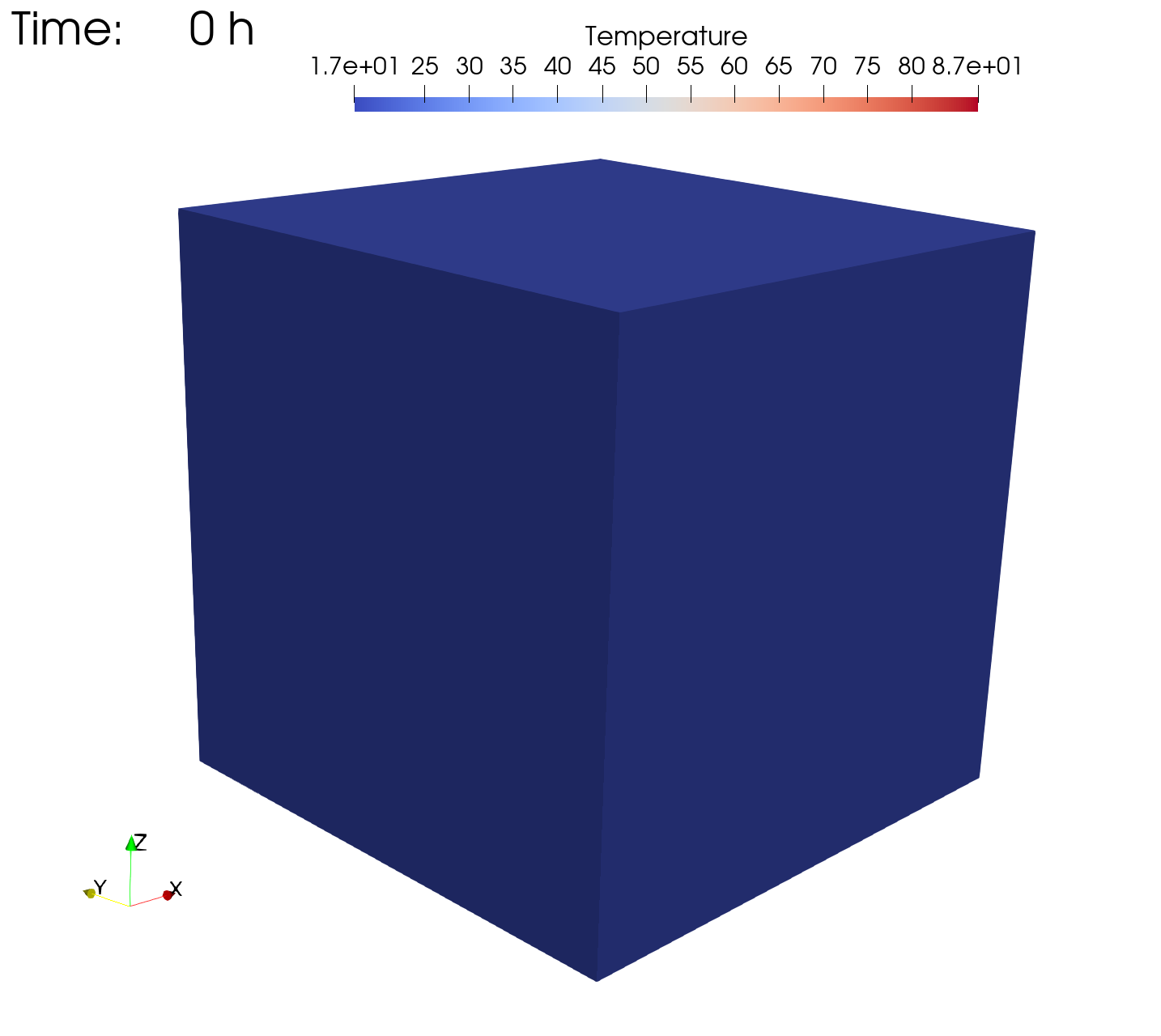}
	\end{subfigure}
	\hfill
	\begin{subfigure}[b]{0.45\textwidth}
		\includegraphics[width=\textwidth]{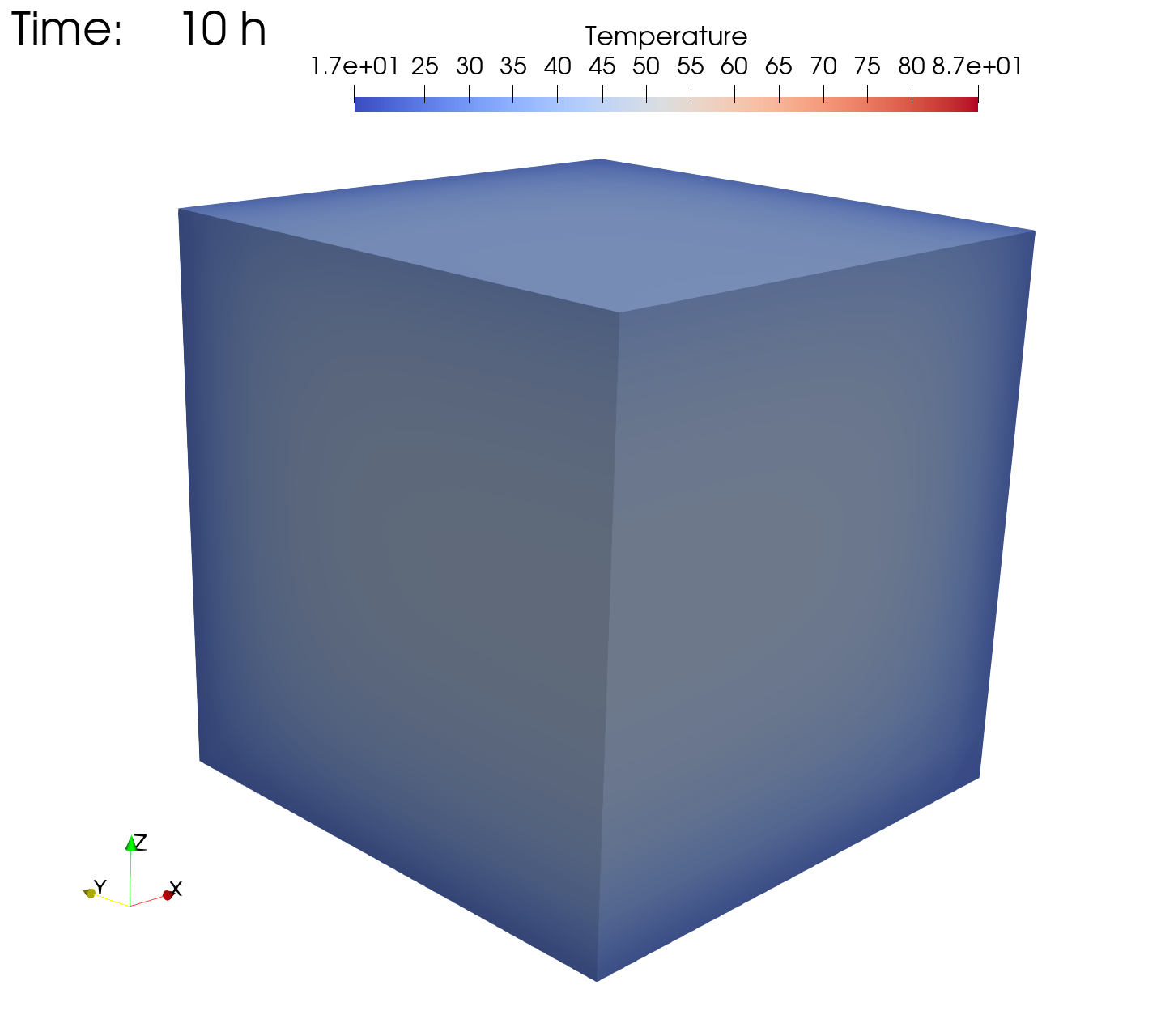}	\end{subfigure}\\
	\begin{subfigure}[b]{0.45\textwidth}
		\centering
		\includegraphics[width=\textwidth]{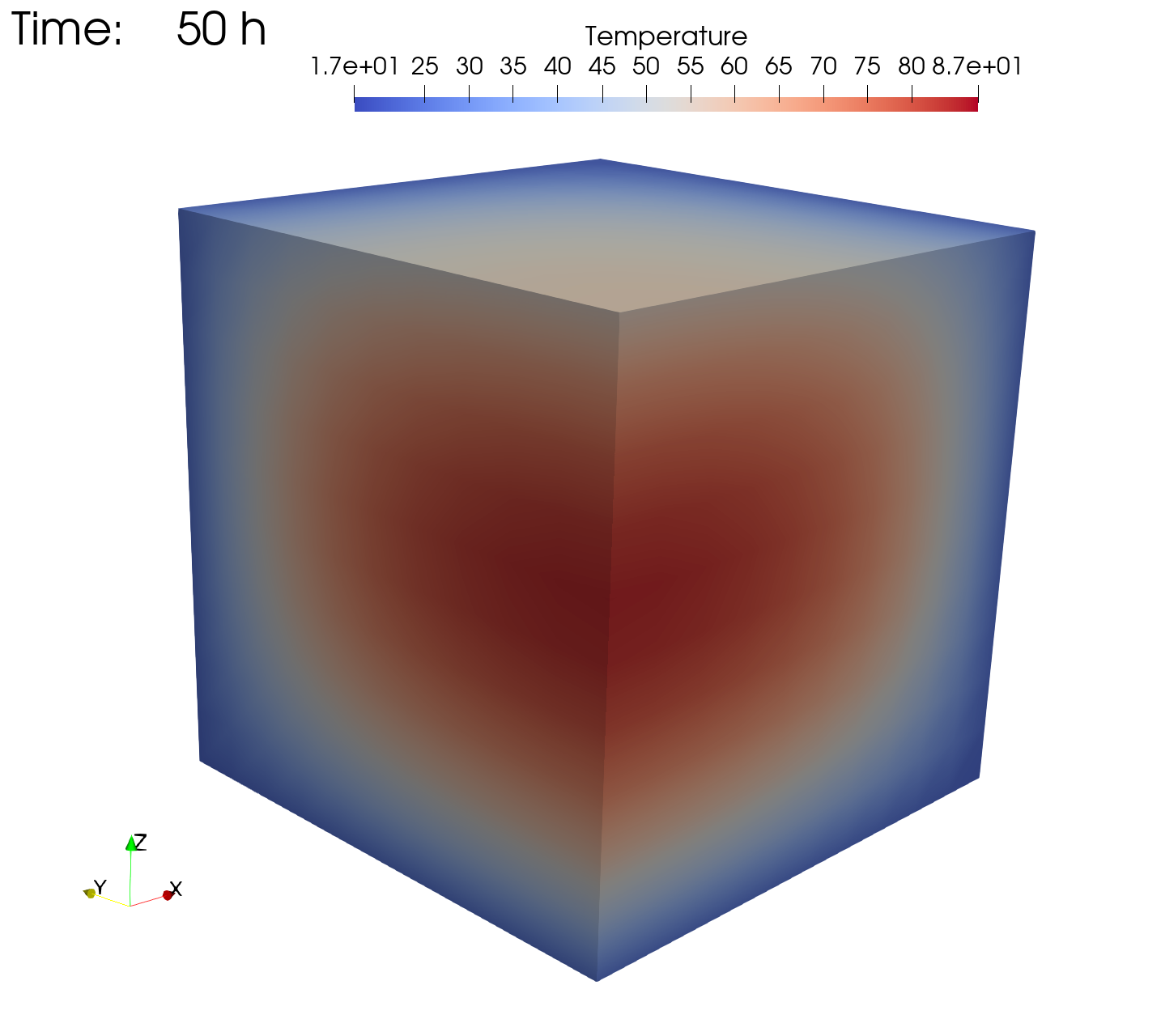}
	\end{subfigure}
	\hfill
	\begin{subfigure}[b]{0.45\textwidth}
		\includegraphics[width=\textwidth]{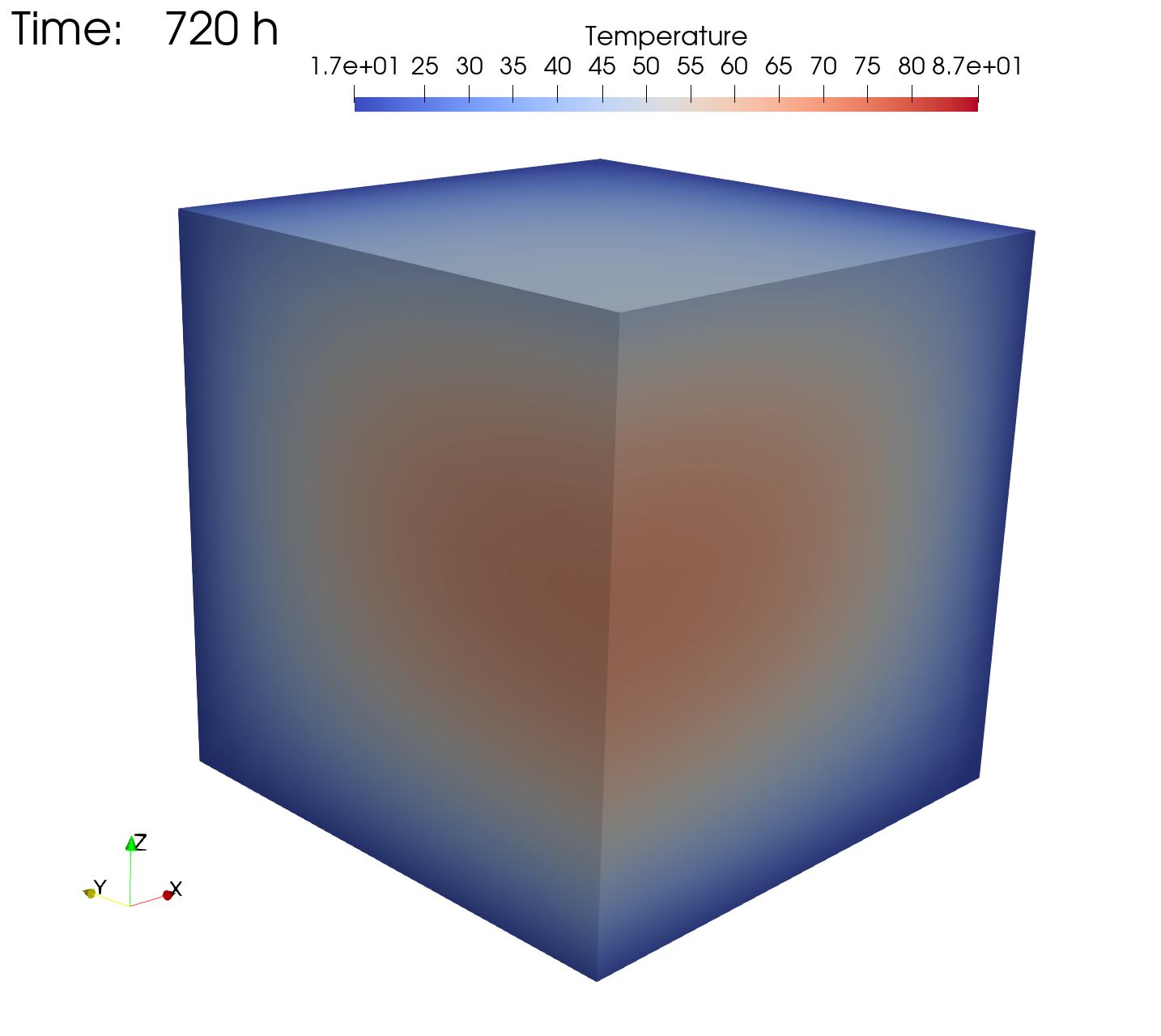}	\end{subfigure}
	\caption{Example 5: Evolution of the temperature field in the meshed domain at four different observation times: $0$ h (top-left), $10$ h (top-right), $50$ h (bottom-left) and $720$ h (bottom-right).}
	\label{fig:ex4TempFields}
\end{figure}

Let us consider the following parameters: $Q_{FH} = 330$ kJ/kg and $T_{var} = 6$ $^\circ$C. We use a time step $\Delta t = 1/3$ h and the final time of the simulation is 240 hours (10 days).
Numerical time integration is done using Backward Euler. 
The temperature fields obtained at times $0$, $10$, $50$ and $720$h are shown in Figure~\ref{fig:ex4TempFields}.
Let us remark that by symmetry, only one quarter of the domain is meshed, and a regular structured mesh of 1331 nodes and 6000 four-node tetrahedron elements is considered.

In order to apply the set propagation method, the source term is factorized using the formulation described in Section~\ref{sec:transformation_firstorder}. %
The source term for this model at time $t$ can be written as:
\begin{equation}\label{eqn:ex4ft}
	\bff_\theta(t)
	=
	\bff_{Q} Q_{FH} \rho m e^{-mt}
	+ \bff_{T_\infty}  T_{var}\frac{1}{2} \sin\left( \omega t - \frac{\pi}{2} \right)
	+ \bff_{T_\infty} \left( T_{min}+T_{var}\frac{1}{2}\right),
\end{equation}
where $\bff_{Q}$ and $\bff_{T_\infty}$ are vectors given by the finite element discretization integration. %
The terms in Eq.~\eqref{eqn:ex4ft} can be factorized to the form of Eq.~\eqref{eqn:decompEta} using:
\begin{equation}
	\left\{
	\begin{array}{ll}
		\bff_0^{(1)} = \bff_{T_\infty} \qquad &\eta^{(1)}(t) = T_{min} + T_{var} \frac{1}{2} \\
		\bff_0^{(2)} = \bff_{Q} \cdot \rho \cdot m \qquad &\eta^{(2)}(t) = Q_{FH} e^{-m t} \\
		\bff_0^{(3)} = \bff_{T_\infty} \qquad  &\eta^{(3)}(t) =  T_{var}\frac{1}{2} \sin\left( \omega t - \frac{\pi}{2} \right).
	\end{array}
	\right.
\end{equation}
The function $\eta^{(1)}(t)$ can be considered as the solution of the differential equation $\dot{y}=0$ with initial condition $T_{min} + T_{var} \frac{1}{2}$, the function $\eta^{(2)}$ is the solution of the differential equation $\dot{y} = -m y$ with initial condition $Q_{FH}$ and the function $\eta^{(3)}$ is the solution of the differential equation $\ddot{y}+\omega^2 y = 0$ with initial conditions $y(0) = -T_{var} \frac{1}{2}$ and $\dot{y}(0) = 0$.

This can be written in matrix form as:
\begin{equation}
	\left[ 
	\begin{array}{c}
		\dot{\bfx} \\
		\dot{\xi_1^{(1)}} \\
		\dot{\xi_1^{(2)}} \\
		\dot{\xi_1^{(3)}} \\ 
		\dot{\xi_2^{(3)}} 
	\end{array} \right]
	=
	\left[ 
	\begin{array}{ccccc}
		-\bfC^{-1}\bfK & -\bfC^{-1}\bff_{T_\infty} & -\bfC^{-1} \bff_{Q} \rho m & -\bfC^{-1} \bff_{T_\infty} & \bf0 \\
		\bszer  & 0 & 0 & 0 & 0 \\
		\bszer  & 0 & -m & 0 & 0\\
		\bszer & 0 & 0 & 0 & 1 \\
		\bszer & 0 & 0 & -\omega^2 & 0 \\
	\end{array} \right]
	\left[ 
	\begin{array}{c}
		\bfx \\
		\xi_1^{(1)} \\
		\xi_1^{(2)} \\
		\xi_1^{(3)} \\ 
		\xi_2^{(3)} 
	\end{array} \right],
\end{equation}

\begin{figure}[htb]
	\centering
	\includegraphics[width=0.8\textwidth]{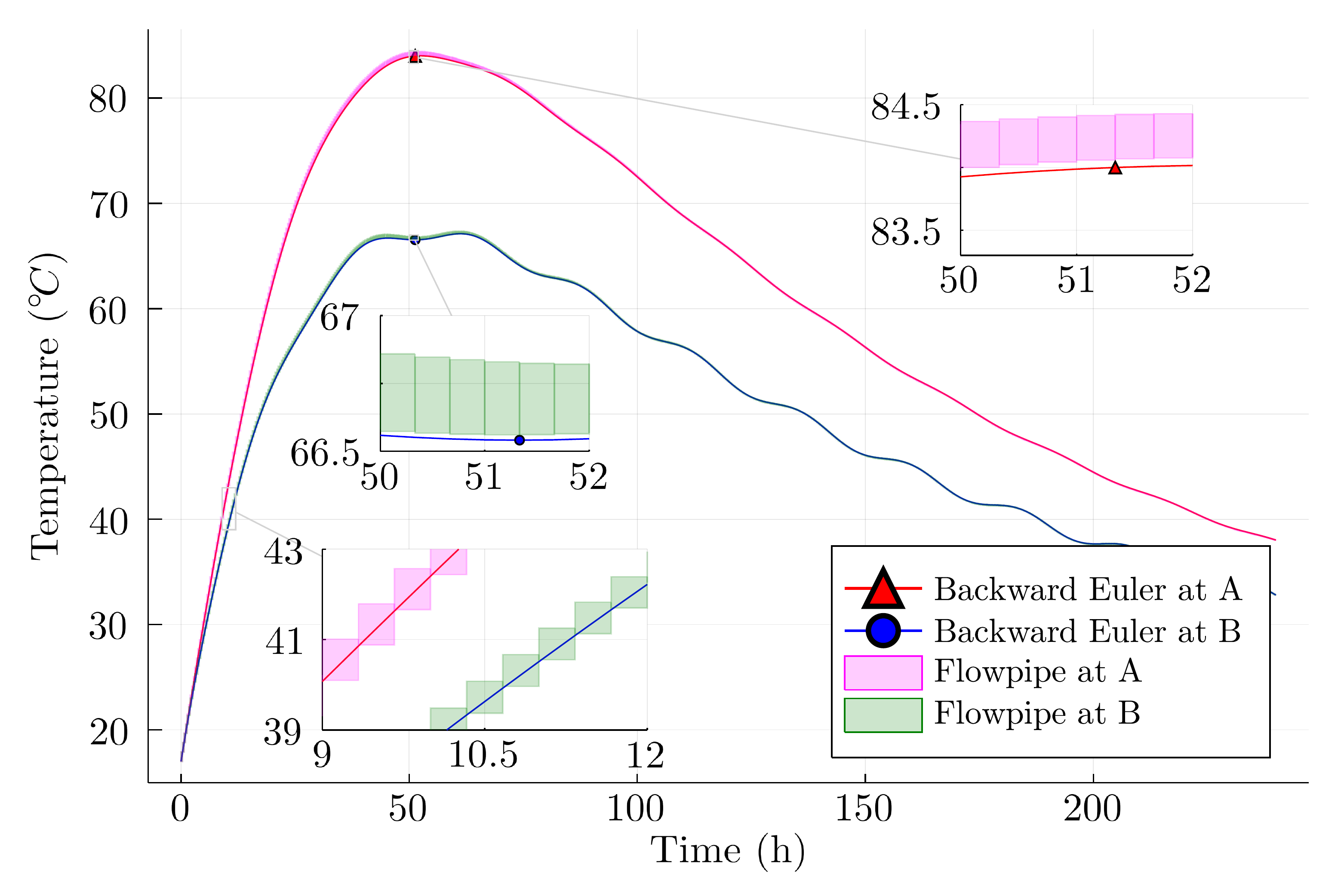}
	\caption{Example 5: temperature evolution at control points A and B.}
	\label{fig:ex4Case1Evol}
\end{figure}

Let us solve the problem considering fixed values for the parameters $Q_{FH}$ and $T_{var}$. In that case, the initial values for the input variables $\xi$ are single fixed values. %
The evolution of the temperatures obtained in control points A and B is shown in Figure~\ref{fig:ex4Case1Evol}. %

We observe that the Backward Euler method may underestimate the temperature values, as shown in the zoomed region in Fig.~\ref{fig:ex4Case1Evol} for times within the interval $[50,52]\textrm{h}$.
The flowpipe is computed using support functions along the canonical directions associated to the control points $A$ and $B$.

\subsubsection{Set propagation under variation in heat sources}

\begin{figure}[htb]
	\centering
	\begin{subfigure}[b]{0.4\textwidth}
		\centering
		\includegraphics[width=\textwidth]{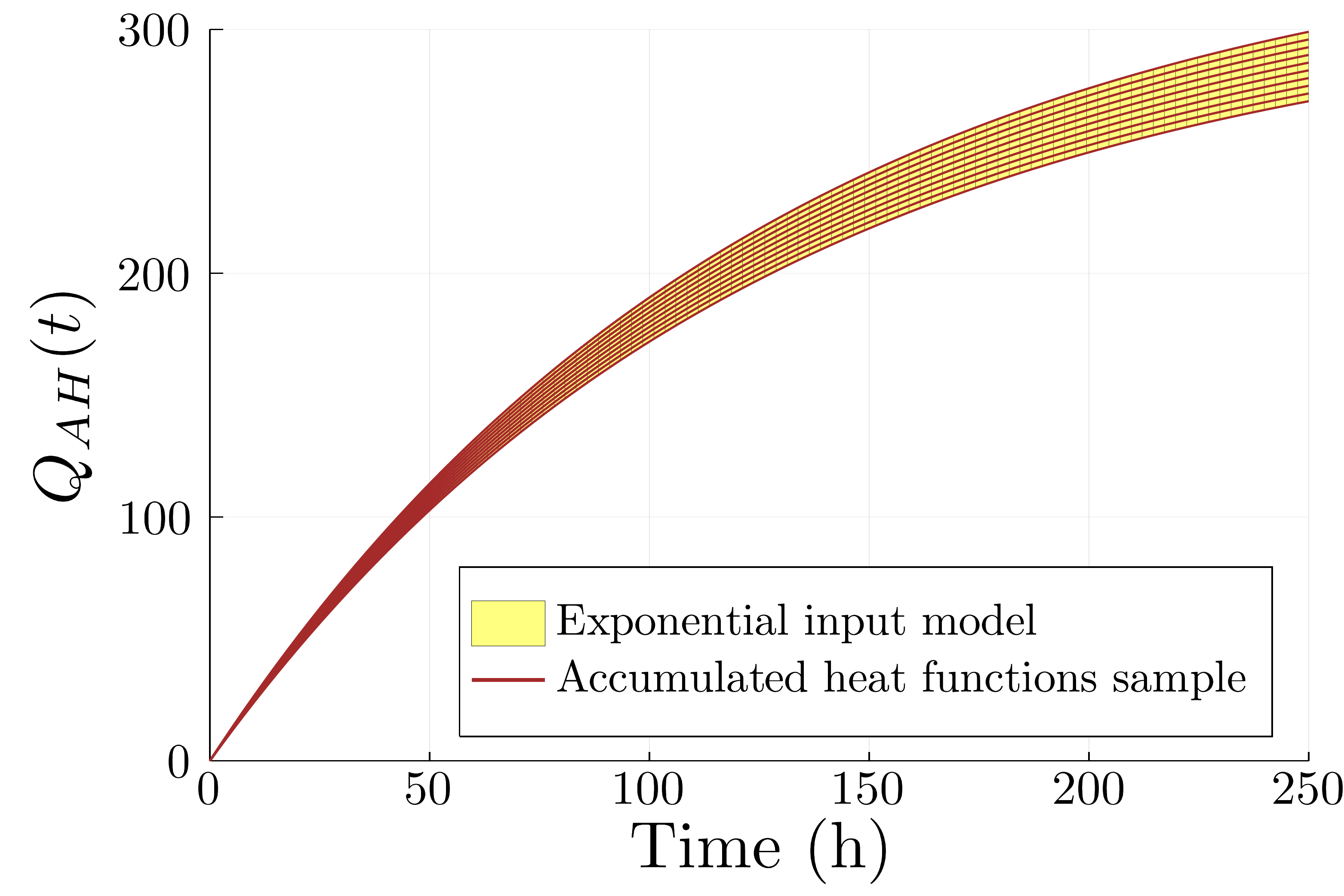}
		\caption{Accumulated heat functions considered.}
		\label{fig:acumheatRA}
	\end{subfigure}
	~~
	\begin{subfigure}[b]{0.5\textwidth}
		\centering
		\includegraphics[width=\textwidth]{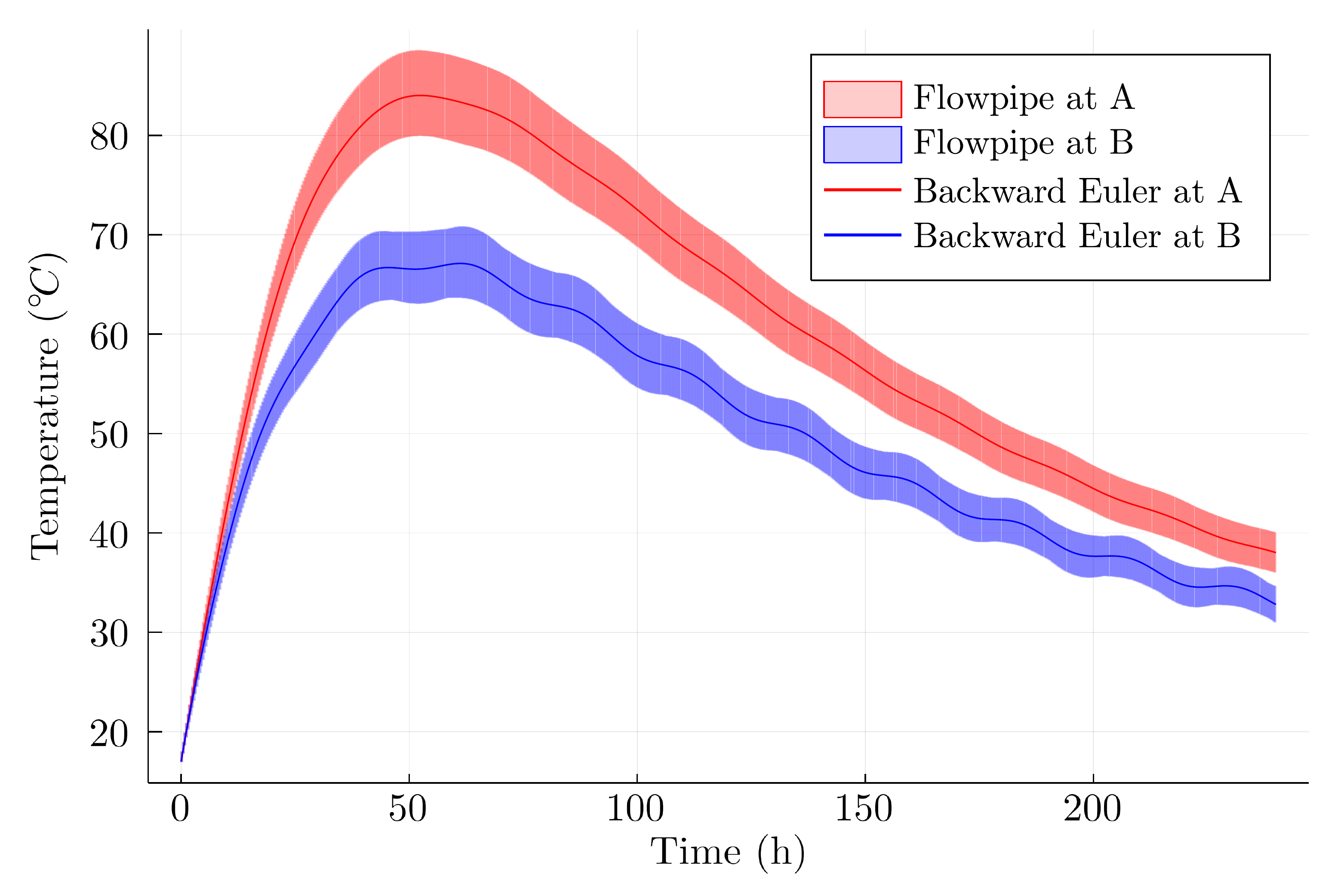}
		\caption{Results considering variations.}
		\label{fig:ex4case2}
	\end{subfigure}
	\caption{Example 5: set propagation results for the heat of hydration model under variation in heat sources.}
	\label{fig:ex5}
\end{figure}

Let us consider a variation in the parameters $Q_{FH}$ and $T_{var}$. For $Q_{FH}$, a $5 \%$ error interval is used, then $Q_{FH}\in [313.5, 346.5]~\text{kJ/kg}$, while for $T_{amb}$, a $2 ^\circ C$ variation is considered, then $T_{amb}\in[4, 8]~^\circ C$. %
These intervals are used to construct the set $\mcC_0$ of initial values of the input variables. %
A plot of the set of accumulated heat functions is shown in Fig.~\ref{fig:acumheatRA}.
The flowpipes corresponding to the evolution of the temperature in the control nodes are shown in Fig.~\ref{fig:ex4case2}.
The solution obtained with Backward Euler for the mean values of the parameters (as in Section~\ref{sec:example4_fixed_params}) is represented with solid lines.

\begin{table}[htb]
	\centering
	{\rev
		\begin{tabular}{@{}llll@{}}
			\cmidrule(l){2-4}
			& Set Propagation & Backward Euler \\ \midrule
			Case 1 (fixed parameters) & 2.05 s (373 MB)  & 3.7 s (78 MB)  \\ \midrule
			Case 2 (variation of heat sources) & 3.1 s (555 MB)  & 3.8 s (78 MB) 
		\end{tabular}
		\caption{Execution times and memory usage for the concrete casting heat of hydration example. The column for Backward Euler corresponds to solving for one initial condition.}
		\label{tab:runtimes_heat3d}
	}
\end{table}

Execution times are reported in Table~\ref{tab:runtimes_heat3d}. We remark that the computational cost of the set propagation is approximately the same whether we use fixed parameters or parameters with uncertainty. Therefore, reachability analysis becomes an efficient option for safety verification in the design of massive concrete structures under uncertainty.

\section{Conclusions} \label{sec:conclusions}

In this article a new approach for numerically solving the differential equations that govern linear heat transfer and structural dynamics problems is presented.
The approach is based on recent advances in reachability analysis for high-dimensional linear systems.
Our work is an attempt to bridge the gap between spatially discretized equations using the Finite Element Method and the incipient field of set-based numerical integration.
The results of our approach are compared with reference numerical time integration methods through the resolution of {\rev five} examples.

{\rev Given a system of differential equations with a set of admissible initial conditions, the} set propagation method consists {\rev in} computing an initial set that contains all possible trajectories within {\rev a time interval $[0, \delta]$}, and {\rev propagating} that set using set-based techniques.
The union of {\rev the obtained} sets is called \textit{flowpipe}.
Let us remark that, while standard numerical integration provides single values for each time point, our approach returns a \textit{set} for each time interval that is guaranteed to contain all possible solutions of the differential equation for the given initial state(s) and input(s).
We formalized the idea of homogeneization of the governing equations subject to set-based initial conditions and inputs. 
We also presented a simple modification of an existing theoretical result for conservative time discretizaton, which ensures more accurate results at a marginal increase in the computational cost.

In Example 1 a fundamental test problem of a single degree of freedom oscillator {\rev was considered}. 
{\rev This problem is} typically used to evaluate numerical properties such as Period Elongation (PE) and Amplitude Decay (AD). %
The results obtained let us conclude that the set propagation method has neither PE nor AD, {\rev therefore it provides accurate results}. %

In Example 2 a reference uni-dimensional wave propagation problem {\rev was} considered. %
In this example we found that the {\rev proposed method can handle a relatively high number of degrees of freedom. %
Moreover it allows to accurately compute the solution at a specific node with execution times in the same order and considerably less memory usage than Bathe and Newmark methods}.

{\rev
In Example 3 a two-dimensional wave propagation problem was considered. %
For single initial conditions the three compared methods provide accurate results with similar execution times and memory usage. %
For distributed initial conditions, our approach is able to compute, in a single execution, an envelope of the solutions at a specified control node for all admissible conditions. %
On the other hand, the results obtained for the numerical integration methods, let us infer that a prohibitively large number of simulations would be required to provide a fair estimation of the envelope. %
}

In Example {\rev 4} the Set Propagation and Backward Euler methods are compared in a unidimensional heat transfer problem.
A dense set of spatially distributed initial condition functions is considered.
The obtained flowpipe accurately encloses the analytic solution for all initial functions considered.
Moreover, we computed the gradient of the solution using a support function approach. 
While numerical integration approaches such as Backward Euler require to perform several simulations for different initial functions, the proposed approach provides a robust enclosure with a single set-based integration.
Let us remark that computing a high number of simulations to ensure a very high degree of coverage is computationally demanding due to the curse of dimensionality. {\rev  However, the set propagation method provides robust results in a single analysis,} hence the cost can be drastically reduced.

Finally, the results obtained in Example {\rev 5} let us conclude that the set propagation method can be applied to a realistic three-dimensional massive concrete hydration problem.
The problem includes a combination of different types of heat sources.
The obtained flowpipe is accurate in comparison with numerical solutions, {\rev  showing that the homogeneization strategy can correctly handle} various sources of input uncertainty.
We observed that for some time intervals the solution obtained using the Backward Euler method underestimates the maximum temperature. 
Given the theoretical basis of reachability analysis, and the accuracy of our results, we conclude that set propagation may be a useful tool for the design of massive concrete structures.

{\rev
	
	It is important to recall that the solution of the FEM equations, or spatially discrete system of ODEs, is already in error when compared with the original exact problem in the continuum. Estimating and absorbing such error in a set-based fashion was not addressed in this article.
Moreover, it is well-known that high order modes are artifacts of the discretization process, thus numerical or algorithmic dissipation is often a favorable or even required mechanism to improve the accuracy of the results when compared to the exact problem. %
On the other hand, solutions obtained by reachability analysis are not numerically damped, therefore do not include artificial damping behavior. %
This, apart from being an advantage of the method itself, may be useful to quantify how numerical solutions depart from the true trajectories of the FEM equations.
}

Several research directions remain open.
It would be enlightening to apply the developed method in real engineering problems, in particular in the context of Example 5.
For that purpose the theoretical development of other, more general input models might be of interest, as well as  the evaluation of the computational cost when the size of the problems increase. 
The set-based approach might be extended to other scientific domains such as inverse problems and non-linear structural analysis. %
{\rev It would  also be interesting to extend this method to integrate the differential equations given by other finite element approaches, such as Overlapping Finite Elements or the Extended Finite Element Method. }

\section*{Acknowledgements}

The authors would like to thank the \textit{Comisión Sectorial de Investigación Científica} of \textit{Universidad de la República} and the \textit{Agencia Nacional de Investigación e Innovación} for the funding of the \textit{Timbó} portal. %
The authors gratefully acknowledge comments from two anonymous reviewers, to Christian Schilling for valuable discussions at various stages of this work, to Pedro Curto for useful comments in heat transfer analysis, and to Ander Gray for comments on non-probabilistic approaches.

\bibliography{references}

\clearpage
\appendix

\section{Proof of Proposition 1} \label{sec:prop1}

Using the forward-only transformations in \citep{bogomolov2018reach} applied to the special case when the system is homogeneous,
\begin{equation}
	\Omega_0^+ = CH(\mathcal{X}_0, \bfPhi \mathcal{X}_0 \oplus E^+(\bfA, \mathcal{X}_0, \delta)), \label{eq:discretization_step_intersect_pos_app}
\end{equation}
where $E^+(\bfA, \mathcal{X}_0, \delta) = \boxdot\left(\bfP(\vert \bfA \vert, \delta)\boxdot(\bfA^2 \mathcal{X}_0)\right)$
and $\bfP(\bfA, \delta) = \sum_{i=0}^\infty \bfA^i \delta^{i+2}/(i+2)!$, $\delta > 0$.

The exact reachable set at the time point $t = \delta$ is $\mathcal{R}^e(\mathcal{X}_0, \delta) = \bfPhi \mathcal{X}_0$. 
Under the transformation $t \mapsto -t$, the state transition matrix of the linear system \eqref{eq:linearODE} becomes $\bfPhi^{-1} = e^{-\bfA \delta}$. Therefore the reachable states for the time interval $[0, \delta]$ can be enclosed using Eq.~\eqref{eq:discretization_step_intersect_pos_app} backwards in time:
\begin{eqnarray}
	\Omega_0^- &=& CH(\bfPhi\mcX_0, \bfPhi^{-1} \bfPhi \mcX_0 \oplus E^+(-\bfA, \bfPhi \mcX_0, \delta)) \nonumber \\
	&=& CH(\bfPhi\mcX_0, \mcX_0 \oplus E^+(\bfA, \bfPhi \mcX_0, \delta)). \label{eq:discretization_step_intersect_neg_app}
\end{eqnarray}
Since both $\Omega_0^+$ and $\Omega_0^-$ are sound overapproximations of $\mathcal{R}^e(\mcX_0, [0, \delta])$, we conclude that
\begin{equation}
	\mathcal{R}^e(\mcX_0, [0, \delta]) \subseteq \Omega_0^+ \cap \Omega_0^-.
\end{equation}
The convergence property follows from Lemma 3 in \citep{frehse2011spaceex}.

\section{Proof of Proposition 2} \label{sec:prop2}

Let $H_0 = \langle \bfc_0, \bfr_0 \rangle_H$ be the initial hyperrectangle. It can be represented as a zonotope with a diagonal generators matrix, $Z_0 = \langle \bfc_0, \mathbf{D}_0\rangle_Z$ where $\mathbf{D}_0$ is the diagonal matrix whose entries are the elements of $\bfr_0$. Using Eq.~\eqref{eq:recurrence_sol_zonotope}, let $Z_k = \langle \bfc_k, \bfG_k \rangle_Z$ such that $\bfc_k = \bfPhi^k \bfc_0$ and $\bfG_k = \bfPhi^k \bfD_0$ for each $k \geq 0$.
From Prop. 2.2 in \cite{le2009reachability}, the support function of $Z = \langle \bfc, \bfG\rangle_Z$ along direction $\bfd \in \mathbb{R}^n$ is 
\begin{equation} \label{eq:support_function_zonotope}
	\rho(\bfd, Z) = \bfd^T \bfc + \Vert \bfG^T \bfd \Vert_1 = \sum_{i=1}^n d_i c_i + \sum_{i=1}^n \left\vert \sum_{j=1}^n (\bfG)_{ji} d_j\right\vert.
\end{equation}
Let $\bfd = \pm \bfe_i$ for any $i = 1, \ldots, n$ be a canonical direction (and its opposite). Then,
\begin{eqnarray}
	\rho(\pm \bfe_i, Z_k) &=& (\pm\bfc_k)_i + \sum_{i'=1}^n \left\vert  (\bfG_k)_{ii'} \right\vert = (\pm\bfPhi^k \bfc_0)_i + \sum_{i'=1}^n \left\vert  (\bfPhi^k \mathbf{D}_0)_{ii'} \right\vert \nonumber \\ 
	&=& (\pm \bfPhi^k \bfc_0)_i + \left(\vert  \bfPhi^k \vert \bfr_0\right)_{i}.
\end{eqnarray}
Therefore, the radius of the tight hyperrectangular overapproximation of $Z_k$, $H_k = \langle \bfc_k, \bfr_k \rangle_H$, along direction $\bfe_i$ is
\begin{equation}
	(\bfr_k)_i = \dfrac{\rho(\bfe_i, Z_k) + \rho(-\bfe_i, Z_k)}{2} = \left(\vert \bfPhi^k \vert \bfr_0 \right)_i,
\end{equation}
which concludes the proof.

\section{Estimation of PE and AD in Flowpipes} \label{sec:apePEAD}

In this section we describe a numerical procedure for the computation of Period Elongation (PE) and Amplitude Decay (AD) for the set propagation method results presented in Example 1. %

This problem can be reduced to the computation of $T_{num}$ for Eq.~\eqref{eq:PE} and $A_{num}$ for Eq.~\eqref{eq:AD2}. %
In order to compute amplitudes and periods, maximums of the flowpipe solution have to be identified.

Let us consider a flowpipe $\mcF$ formed by the union of a finite number of reach-sets $\{X_k\}_k$. %
The relative maximum of the displacement happens when the conditions $u(t)>0$ and $v(t)=0$ are met. The intersection of $\mcF$ with the set $U^+=\{(u, v) \in \mathbb{R}^2: (u > 0) \wedge (v=0) \}$ results in one or several reach-sets
such that their intersection with $U^+$ is non-empty, i.e. they include the behavior that satisfies the maximum displacement conditions (if more than one reach-set intersects $U^+$, we take their set union). %

To estimate a period $T_{num}$ to the flowpipe $\mcF$ we proceed as follows.
Let $\mcT_m$ be the time interval associated with the $m$-th intersection of $\mcF$ with $U^+$. %
Then, we have $m \cdot T_{num} \in \mcT_m$ for any period $m$, therefore we can define the following time interval:
\begin{equation}
	\centering
	S_T = \bigcap_{m} \mcT_m \frac{1}{m}.
	\label{eq:SDOF_STdef}
\end{equation}
Since $S_T$ is an intersection of intervals, all containing $T_{num}$, then 
\begin{equation}
	T_{num} \in S_T.
\end{equation}
We repeated the same example for different values of $\alpha$ and calculated $S_T$ in each case using $1000$ periods. %
The value of $PE$ reported in Fig.~\ref{fig:PE_all} corresponds to the maximum of $T_{num}$ among all periods.


Regarding the estimation of the amplitude, using the time intersection in Eq.~\eqref{eq:SDOF_STdef} we compute an interval of displacements. The lower value of displacement represents a safe estimation of the amplitude, and also the worst possible case for the Amplitude Decay. This is the value considered for the computation of $AD$ in Fig.~\ref{fig:AD_all}.

\end{document}